\documentclass{gtpart}

\usepackage[all]{xy}
\usepackage{graphicx}
\usepackage{hyperref}
\usepackage{rotating}
\usepackage{pinlabel}
\usepackage{amscd}

\author{Tim Perutz}
\givenname{Tim}
\surname{Perutz}
\address{DPMMS, Centre for Mathematical Sciences, University of Cambridge, Wilberforce Road, Cambridge CB3 0WB, UK.}
\email{tim.perutz@cantab.net}

\title{Lagrangian matching invariants for fibred four-manifolds: II}

\subject{primary}{msc2000}{53D40}
\subject{primary}{msc2000}{57R57}
\subject{secondary}{msc2000}{57R15}

\keyword{Four--manifolds}
\keyword{Lefschetz fibrations}
\keyword{Seiberg--Witten invariants} 
\keyword{pseudo-holomorphic curves}
\keyword{Lagrangian correspondences}

\theoremstyle{plain}
\newtheorem{Thm}{Theorem}[section]

\newtheorem{Lem}[Thm]{Lemma}
\newtheorem{Prop}[Thm]{Proposition}
\newtheorem{Cor}[Thm]{Corollary}

\newtheorem{Conj}[Thm]{Conjecture}

\newtheorem{Rk}[Thm]{Remark}
\newtheorem{Ex}[Thm]{Example}

\newtheorem*{TheoremD}{Theorem D}
\newtheorem*{TheoremE}{Theorem E}
\newtheorem*{TheoremF}{Theorem F}

\theoremstyle{definition}
\newtheorem{Defn}[Thm]{Definition}

\newenvironment{pf}{ \begin{proof} }{ \end{proof} }

\DeclareMathAlphabet\EuScript{U}{eus}{m}{n}
\SetMathAlphabet\EuScript{bold}{U}{eus}{b}{n}

\DeclareFontFamily{U}{eus}{\skewchar\font'60}%
\DeclareFontShape{U}{eus}{m}{n}{<-6>eusm5<6-8>eusm7<8->eusm10}{}%
\DeclareFontShape{U}{eus}{b}{n}{<-6>eusb5<6-8>eusb7<8->eusb10}{}%

\newcommand{\PS}{\mathbb{P}}

\DeclareMathOperator*{\SO}{\mathrm{SO}}

\DeclareMathOperator*{\U}{\mathrm{U}}

\DeclareMathOperator*{\Sp}{\mathrm{Sp}}
\DeclareMathOperator*{\Spin}{\mathrm{Spin}}
\DeclareMathOperator*{\pin}{\mathrm{Pin}}
\DeclareMathOperator*{\orth}{\mathrm{O}}

\DeclareMathOperator{\im}{im}
\DeclareMathOperator{\End}{\mathrm{End}}
\DeclareMathOperator{\rank}{\mathrm{rank}}
\DeclareMathOperator{\Hom}{\mathrm{Hom}}
\DeclareMathOperator{\coker}{\mathrm{coker}}

\DeclareMathOperator{\cone}{\mathrm{cone}}

\DeclareMathOperator{\sym}{Sym}

\DeclareMathOperator{\hilb}{Hilb}
\DeclareMathOperator{\dbar}{\bar{\partial}}

\DeclareMathOperator{\aut}{Aut}
\DeclareMathOperator{\diff}{Diff}
\DeclareMathOperator{\ham}{Ham}

\DeclareMathOperator{\supp}{Supp}
\DeclareMathOperator{\spinc}{\mathrm{Spin}^c}

\DeclareMathOperator{\ind}{ind}

\DeclareMathOperator{\torus}{T}

\DeclareMathOperator{\sect}{sect}
\DeclareMathOperator{\hor}{\mathcal{H}}
\DeclareMathOperator{\CZ}{CZ}
\DeclareMathOperator{\Div}{div}

\DeclareMathOperator{\Gr}{Gr}
\DeclareMathOperator{\Det}{Det}
\DeclareMathOperator{\Lag}{Lag}

\newcommand{\ii}{\mathrm{i}}
\newcommand{\id}{\mathrm{id}}
\newcommand{\Tv}{T^{\mathrm{v}}}
\newcommand{\Th}{T^{\mathrm{h}}}
\newcommand{\ev}{\mathrm{ev}}
\newcommand{\reg}{\mathrm{reg}}
\newcommand{\univ}{\mathrm{univ}}
\newcommand{\crit}{\mathrm{crit}}

\begin{document}

\begin{abstract}
In the second of a pair of papers, we complete our geometric construction of  `Lagrangian matching invariants' for smooth four-manifolds equipped with broken fibrations. We prove an index formula, a vanishing theorem for connected sums and an analogue of the Meng--Taubes formula. These results lend support to the conjecture that the invariants coincide with 
Seiberg--Witten invariants of the underlying four-manifold, and are in particular independent of the broken fibration. 
\end{abstract}
\maketitle

\section{Introduction}
This paper is the sequel to \cite{Pe2}, which we shall refer to as Part I. Our overall aims and main theorems were stated in the introduction to that paper, but briefly, the goal is to understand 
Seiberg--Witten invariants of a four-manifold equipped with a `broken fibration'---that is, a singular Lefschetz fibration in the sense of Auroux--Donaldson--Katzarkov \cite{ADK}---using 
constructions in symplectic topology directly reflecting the geometry of the fibration. We shall not achieve quite as much as that: what we shall manage is to construct `Lagrangian matching invariants', as invariants of the fibration, and demonstrate, in a number of ways, their resemblance to Seiberg--Witten invariants of the four-manifold. However, equality of the Lagrangian matching invariants with the Seiberg--Witten invariants remains conjectural, as does their independence of the fibration. 

\subsection{Fibred coisotropic hypersurfaces}
An important role in the programme is taken by certain  Lagrangian correspondences between symmetric products of surfaces. In Part I we constructed these correspondences
as vanishing cycles for certain degenerations and studied their properties.  They arise from \emph{fibred coisotropic hypersurfaces} in symmetric products. One way to  describe their significance, which we did not do explicitly in Part I, is by means of a theorem which we now proceed to state.  

Let $\Sigma$ be a closed Riemann surface equipped with an area form $\alpha$. Fix $n\geq 2$ and a K\"ahler form $\omega$ on $\sym^n(\Sigma)$ representing a cohomology
class (in $H^2(\sym^n(\Sigma);\R)$) invariant under action of the mapping class group of $\Sigma$.  When $g(\Sigma)>0$ the space of such invariant classes is two-dimensional, and there is a standard basis $(\eta,\theta)$ for it:
\begin{itemize}
\item
$\eta$ is dual to any of the smooth divisors $ \sym^{n-1}(\Sigma)\subset \sym^n(\Sigma)$, embedded by the map $D\mapsto x+ D$ for a fixed $x\in \Sigma$.
\item
$\theta$ corresponds to the intersection form on $H_1(\Sigma;\R)$ under the natural  map $\Lambda^2 H^1(\Sigma)\to H^2(\sym^n(\Sigma))$ (i.e., the usual isomorphism $\Lambda^2 H^1(\Sigma)\to H^2(\mathrm{Jac}(\Sigma))$ followed by $\mathrm{AJ}^*$, where
$\mathrm{AJ}\colon \sym^n(\Sigma)\to \mathrm{Jac}(\Sigma)$ is the Abel--Jacobi map).  
\end{itemize}
When $g(\Sigma)=0$, the class $\theta$ vanishes, but $\eta$ spans the invariant part of the second cohomology.  The classes $\eta$ and $\theta$ will come up repeatedly in this paper, as they did in Part I.

After rescaling $\omega$, we may write its cohomology class as $[\omega]= \eta + \lambda \theta$ where $\lambda$ lies in (some bounded-below interval in) $\R$. 

One obtains a homomorphism between symplectic mapping class groups
\[ \kappa_{n,\lambda} \colon  \frac{\aut(\Sigma,\alpha)}{\ham(\Sigma,\alpha)} \to  \frac{\aut(\sym^n(\Sigma),\omega)}{\ham(\sym^n(\Sigma),\omega)}   \]
by the following procedure. Given $\phi\in \aut(\Sigma,\alpha)$, its mapping torus $Y=\torus(\phi)\to S^1$ is a bundle with fibre $\Sigma$ carrying a natural cohomology class $\alpha_\phi\in H^2(\torus(\phi);\R)$ extending  $\alpha$. A choice of complex structure on its vertical tangent bundle makes its relative symmetric products $\sym^n_{S^1}(Y) = \{ (t,D): t\in S^1, D\in \sym^n(Y_t)\} $ into smooth manifolds. They fibre over $S^1$ with complex fibres. A point which arose in Part I is that there are natural maps 
\[ H^2(Y)\to H^2(\sym^n_{S^1}(Y)), \quad  c\mapsto c^{[1]} \]
and 
\[H^0(Y)\to H^2(\sym^n_{S^1}(Y)),\quad c\mapsto c^{[2]}, \]
defined by means of the universal divisor 
\[ \Delta = \{( x, D) : \, x\in \supp(D)\} \subset Y\times _{S^1}\sym^n_{S^1}(Y)\] 
and its Poincar\'e dual cohomology class $\delta$.
The second projection $\mathrm{pr}_2\colon Y \times _{S^1}\sym^n_{S^1}(Y) \to \sym^n_{S^1}(\Sigma)$ is a $\Sigma$--bundle. Using also the first projection $\mathrm{pr}_1$, one puts
\[ c^{[1]} = (\mathrm{pr}_2)_!(\mathrm{pr}_1^*(c)\smile \delta),\quad c^{[2]}=   (\mathrm{pr}_2)_!(\mathrm{pr}_1^*(c)\smile \delta\smile \delta) .\]
There is a non-empty, convex set of closed two-forms on $\sym^n_{S^1}(Y)$ which are K\"ahler on the fibres and which globally represent the unique cohomology class in the linear span of $(\alpha_\phi)^{[1]}$ and $1^{[2]}$ which extends the class $\eta + \lambda \theta$ on a reference fibre. One verifies that $1^{[2]}$ restricts to the fibre $\sym^n(\Sigma)$ as $2(n\eta - \theta)$; this, in tandem with MacDonald's formula $c_1(\sym^n(\Sigma))=(n+1-g)\eta - \theta$, shows that the class extending $\eta + \lambda \theta$ is
\begin{equation}\label{linear comb}
\frac{1+\lambda n}{\langle \Sigma, \alpha \rangle} (\alpha_\phi)^{[1]} - \frac{\lambda}{2} 1^{[2]}.
\end{equation}
\emph{We then define $\kappa_{n,\lambda}([\phi])$ to be the Hamiltonian-isotopy class of the symplectic monodromy of $\sym^n_{S^1}(Y)$.} One checks that this procedure defines a homomorphism.

In Section 2 of Part I, we defined spherically-fibred coisotropic submanifolds 
and Dehn twists along them.

\begin{Thm}
Let  $\tau_\gamma$ denote the positive Dehn twist along the embedded circle $\gamma \subset \Sigma$. Then, assuming $\lambda>0$, $\kappa_{n,\lambda}(\tau_\gamma)$ coincides with the Hamiltonian isotopy class represented by a fibred Dehn twist along an $S^1$--fibred coisotropic submanifold $V_\gamma \subset \sym^n(\Sigma)$, itself determined by $\gamma$ up to Hamiltonian isotopy. The reduced space of $V_\gamma$ is diffeomorphic to $\sym^{n-1}(\bar{\Sigma})$, where $\bar{\Sigma}$ is the result of surgery along $\gamma$.
\end{Thm}
This theorem is essentially due to I Smith \cite[Proposition 3.7]{Smi}, drawing on ideas of Seidel, though to obtain the monodromy symplectically (up to Hamiltonian isotopy), and not merely smoothly, we invoke the monodromy theorem from Part I, Section 2. 

The theorem follows readily from Part I, in which the construction of $V_\gamma$ was given in detail (Theorem A). Beware that whilst it is easy to think of a hypersurface in $\sym^n(\Sigma)$ determined by $\gamma$, namely the image of the natural map $\gamma\times \sym^{n-1}(\Sigma)\to \sym^n(\Sigma)$, this is not  diffeomophic, let alone isotopic, to $V_\gamma$ (it is, however, homologous to $V_\gamma$). The significance of $V_\gamma$ is more subtle: it is the vanishing cycle of a relative Hilbert scheme of points on the fibres of an elementary Lefschetz fibration. The monodromy theorem then shows that the monodromy of an elementary symplectic Morse--Bott fibration, such as this Hilbert scheme, is the fibred Dehn twist about its vanishing cycle (modulo Hamiltonian isotopies).  

\begin{Cor}
When $\lambda>0$, $\im(\kappa_{n,\lambda})$ is generated by fibred Dehn twists. 
\end{Cor}
\begin{pf}
We claim that $\aut(\Sigma,\alpha)/ \ham(\Sigma,\alpha)$ is generated (as a group, not a monoid) by Dehn twists; the corollary then follows from the theorem.

Moser \cite{Mos} proved that $\aut(\Sigma,\alpha)$ has the same homotopy type as $\diff(\Sigma)$. Thus $\pi_0 \aut(\Sigma,\alpha)$ is the ordinary mapping class group $\pi_0 \diff(\Sigma)$, which is generated by (area-preserving) Dehn twists. Hence it suffices to show that the kernel of the quotient map $\aut(\Sigma,\alpha)/ \ham(\Sigma,\alpha)\to \pi_0  \aut(\Sigma,\alpha) $ is generated by Dehn twists. This kernel is isomorphic, via the flux homomorphism, to $H^1(\Sigma;\R)/\Gamma$, where $\Gamma$ is the flux group. By the Earle--Eells theorem \cite{EE}, $\Gamma$ is zero when $g(\Sigma)>1$ and  $\Gamma=H^1(\Sigma;\Z)$ when $g(\Sigma)=1$. By considering the Dehn twist along a curve $\gamma$, followed by the inverse twist along a parallel curve $\gamma'$, one sees that elements of  $H^1(\Sigma;\R)/\Gamma$ are indeed generated by Dehn twists.
\end{pf}
It is unknown whether the corollary holds when $\lambda\leq 0$. 

%Our aim is to put these vanishing cycles to work in the context not just of the Lefschetz fibrations from which they arise, but also of \emph{broken fibrations}.

\subsection{Plan of the paper}
In Section 4 of Part I we showed how families of Lagrangian correspondences, parametrised by $S^1$, could be cast as Lagrangian boundary conditions associated with elementary broken fibrations (Theorem B). In Part II we continue where we left off.  We begin, in Section \ref{Matching conditions}, by restating Theorem B, or rather, stating a corollary of it which is our jumping-off point for this part. In Section \ref{pseudohol} we set out what we need from pseudo-holomorphic curve theory. This is then applied to relative Hilbert schemes of points arising from broken fibrations, with the Lagrangian boundary conditions we have constructed. 

The definition of the Lagrangian matching invariants is given in Section \ref{definition}. It invokes a particular case of the theory of pseudo-holomorphic sections already described. However, rather than applying the theory in its `raw' form, we show that  the resulting invariants $\EuScript{L}_{(X,\pi)}$ for a broken fibration $\pi$ on a four-manifold $X$ can be organised in a way which is strongly reminiscent of Seiberg--Witten theory: 
\begin{itemize}
\item	
A `topological sector' for our theory determines (and is in many cases determined by) a $\spinc$-structure for $X$ (Section \ref{invariants}).
\item
The local expected dimension of the moduli space of pseudo-holomorphic sections defining the invariants is equal to that of the Seiberg--Witten moduli space (Theorem D).
\item
There are relative invariants in Floer homology groups $HF_*(Y,\mathfrak{t})$, which one can take to be finitely generated over $\Z$ for each $\spinc$-structure. These are constructed as summands (corresponding to components of the twisted loopspace) in the fixed-point Floer homology groups $HF_*(\kappa_{n,\lambda}(\phi))$, where $Y=\torus(\phi)$, and $\lambda$ is chosen so as to avoid periods. The group $HF_*(Y,\mathfrak{t})$ is graded by the $\Z$-set of homotopy classes of oriented two-plane fields underlying $\mathfrak{t}$ (Section \ref{Gradings}).  
Such geometric gradings appear in Kronheimer and Mrowka's monopole Floer homology theory for three-manifolds \cite{KM}.
\end{itemize} 
Computations are carried out in Section \ref{calculations}, including the proofs of Theorems E (`Meng--Taubes formula') and F (`vanishing for connected sums'), both of which were stated in the introduction to Part I.

\subsection{Acknowledgements} I thank Simon Donaldson, Paul Seidel, Ivan Smith and Michael Usher for their useful input. This work was supported by the EPSRC via Research Grant EP/C535995/1.

\section{Matching conditions}\label{Matching conditions}
In this section we state some of our conclusions from Part I in a language we shall find convenient for the pseudo-holomorphic curve theory still to come.
\begin{Defn}\label{matching}
Let $S$ be a compact surface with an even number of boundary components. A {\bf matching} for $S$ is a decomposition 
\[\partial S = \partial S_+ \amalg  \partial S_-  \] 
(defined by a continuous function $\partial S \to \{\pm 1\}$) together with an orientation-reversing diffeomorphism $\tau \colon  \partial S_+ \to \partial S_-$.
\end{Defn}

In Part I, Definition 2.1, a \emph{symplectic Morse--Bott fibration} was defined in terms of data $(E,\pi,\Omega,J_0,j_0)$, where $E$ is a manifold with boundary, $\pi\colon E \to S$ a smooth proper map to a surface, mapping $\partial E$ submersively to $\partial S$; $\Omega$ is a closed two-form, non-degenerate on the tangent distribution $\ker(D\pi)$; $J_0$ (resp. $j_0$) is the germ of an almost complex structure near the critical set $X^\crit$ (resp. near $\pi(X^\crit)$), making $D\pi$ holomorphic. $X^\crit$ is assumed to be a submanifold, and the complex Hessian form on its normal bundle non-degenerate. We additionally made certain integrability assumptions near $X^\crit$, which we do not repeat here.

\begin{Defn}\label{lmc}
Let $(E^{2n+2},\pi,\Omega, J_0,j_0)$ be a symplectic Morse--Bott fibration over a surface $S$ equipped with a matching. Let $(\partial E_{\pm}, \pi_\pm,\Omega_\pm)$ be the restriction of $E$ to $\partial S_{\pm}$.  A {\bf Lagrangian matching condition} for $(E,\pi)$ is a submanifold $Q^{n+1} \subset \partial E_+ \times_{\partial S_+} \tau^* (\partial E_-) $ such that (i) the projection $ Q\to \partial S_+$ is a proper, surjective submersion, and (ii) $((-\Omega_+) \oplus \tau^* \Omega_- )| Q=0$. 
\end{Defn}

Our main examples of matched surfaces arise from broken fibrations (Part I, Definition 1.1). Let $(X,\pi)$ be a broken fibration over a closed, oriented surface $S$. Let $Z$ be the one-dimensional part of $X^\crit$, and suppose that $\pi|X^\crit $ is injective. 

Let $N$ be a narrow, closed tubular neighbourhood of $\pi(Z)\subset S$, with retraction $r\colon N\to Z$. Then $r$ gives rise to a matching for $S' := S\setminus \mathrm{int}(N)$. Indeed, if $\Gamma$ is a component of $\pi(Z)$, there are two components of $\partial S'$ which retract to it under $r$. The one for which the fibres have lower Euler characteristic (higher genus, in the connected case) is designated as belonging to $\partial S'_+$, the other as belonging to $\partial S_-$. Clearly the matching $\tau$ is an invariant of $(X,\pi)$, up to isotopy. 

Let $X'=X|S'$ be the restricted fibration.

Now let $\nu \colon S' \setminus S^\crit \to \Z_{\geq 0}$ be a locally constant function with the property that
\[  2\nu(s) + \chi(X_s) = \mathrm{const}. \]
For example, $\nu$ could be the intersection number $s\mapsto \beta\cdot [X_s]$ for some relative homology class $\beta\in H_2(X,Z;\Z)$ with $\delta (\beta)= [Z]\in H_1(Z;\Z)$ (see Lemma 4.5).

Choose a positively oriented complex structure $J$ on the tangent distribution $\Tv X' $. Let 
\begin{equation}\label{Xnu}
X^{[\nu]}= \hilb^{\nu}_{S'} (X') \end{equation} 
be the relative Hilbert scheme, i.e., the disjoint union, over components $C\subset S'$, of the relative Hilbert schemes $\hilb^{\nu(C)}_{C}(X|C)$ 
(see Part I, Section 3). It comes with a map $\pi^{[\nu]}\colon X^{[\nu]} \to S'$, and a complex structure $J^{[\nu]}$ on its vertical tangent distribution $\ker D\pi^{[\nu]}$ (this distribution is not of constant rank).

Suppose that $W\in H^2(X;\R)$ is a class such that $\langle W, F \rangle >0 $ for every class $F\in H_2(X;\Z)$ represented by a component of a fibre of $\pi$. Such a $W$ exists by definition of broken fibrations.

\begin{Defn}\label{admissible}
A closed two-form $\Omega$ on $X^{[\nu]}$ is {\bf $W$--admissible} if it meets the following conditions:
\begin{itemize}
\item
$(X^{[\nu]},\pi^{[\nu]},\Omega,J_0,j_0)$ is a symplectic Morse--Bott fibration for some $(j_0,J_0)$, where $J_0$ agrees with $J^{[\nu]}$ on tangent vectors where both are defined.
\item  $\Omega$ is K\"ahler on each regular fibre $\sym^{\nu(s)}(X_s)$, and also on the components of $\crit(\pi^{[\nu]})$.
\item There is a $\lambda>0$ such that the restriction of $\Omega$ to each regular fibre $\sym^{\nu(s)}(X_s)$ represents the class $\eta_{X_s}+ \lambda \theta_{X_s}$.
\item  The restriction of $\Omega$ to $\partial X^{[\nu]}$ represents the cohomology class
\[ \frac{1+\lambda n}{\langle W , \mathrm{fibre} \rangle} W^{[1]}-\frac{\lambda}{2}1^{[2]}  \]
occurring in the construction of $\kappa_{n,\lambda}$ (cf. equation (\ref{linear comb})).
\end{itemize}
It is {\bf admissible} if it is $W$--admissible for some $W$.
\end{Defn}
\begin{Lem}\label{admissible existence}
$W$--admissible two-forms exist. Moreover, the set of admissible forms is convex, hence they are unique up to deformation.
\end{Lem}
\begin{pf}
For existence, one uses the Thurston--Gompf method: see \cite[Theorem 2.2]{Gom} and also \cite[Lemma 4.4]{DS}. One can use equation \ref{c1} to rewrite the class
$ \frac{1+\lambda n}{\langle W , \mathrm{fibre} \rangle } 
 (W|\partial X')^{[1]} -\frac{\lambda}{2} 1^{[2]} $ in the cohomology of the boundary as 
\[   \left(\frac{1+\lambda n}{\langle W , \mathrm{fibre} \rangle} (W|\partial X') 
+ \frac{\lambda}{2} c_1(\Tv X)|\partial X' \right)^{[1]} - \lambda c_1(\Tv \partial X^{[\nu]}). \]
One then looks for a closed two-form globally representing the class
\[ \left(\frac{1+\lambda n}{\langle W , \mathrm{fibre} \rangle} 
(W|\partial X') + \frac{\lambda}{2} c_1(\Tv X')  \right)^{[1]} - \lambda c_1(T X^{[\nu]}) . \] 
This is possible locally in the base: this is a non-trivial statement near a critical value, but it was established in Part I, Section 3. Hence, by Thurston--Gompf patching, it is also possible globally. The convexity statement is obvious.
\end{pf}
Admissibility depends on the choice of $J$, but the set of $J_t$-admissible forms varies continuously in a path $J_t$.

Our rephrasing of Theorem B from Part I is the following:
\begin{Thm}\label{Theorem B}
For any choice of $J$-admissible two-form $\Omega$, the symplectic Morse--Bott fibration $(X^{[\nu]},\pi^{[\nu]},\Omega, J_0, j_0)$ admits a canonical Lagrangian matching condition $\EuScript{Q}$, up to isotopy.  In a smooth path of complex structures $J_t$, and of $J_t$-admissible two-forms $\Omega_t$, there are smoothly-varying matching conditions $\EuScript{Q}_t$.
\end{Thm}

\section{Pseudo-holomorphic sections with Lagrangian matching conditions}\label{pseudohol}
In this section we set out a general framework for `Lagrangian matching invariants', not specifically tied to the Lagrangian matching conditions arising from broken fibrations on four-manifolds.

\subsection{Almost complex structures adapted to symplectic Morse--Bott fibrations}
Let $M$ be a smooth manifold and $V\to  M $ an oriented real vector bundle. 
Denote by $\EuScript{J}(V)$ the space of all orientation-compatible, $C^\infty$ complex structures on $V$, with its $C^\infty$ topology. If $(V,\sigma)$ is a symplectic vector bundle then the subspace $\EuScript{J}(V,\sigma)\subset\EuScript{J}(V)$ of compatible complex structures---those $J\in \EuScript{J}(V)$ such that the formula $g_x(v_1,v_2):=\sigma(v_1,J v_2)$ defines a metric on the vector bundle $V$---is contractible. 

\begin{Defn}
Let $(E,\pi,\Omega, J_0,j_0)$ be a symplectic Morse--Bott  fibration over $S$, and let $E^*=E\setminus \crit(\pi)$. A pair $(J,j)\in\EuScript{J}(TE)\times \EuScript{J}(TS)$ is {\bf adapted} to the fibration if 
\begin{enumerate}
\item           
$\pi$ is $(J,j)$-holomorphic, i.e. $D\pi \circ  J = j\circ D\pi$;
\item
$J|{\Tv E^*}\in \EuScript{J}(\Tv E^*, \Omega|{\Tv E^*})$, that is, $J$ is compatible with $\Omega$ on the smooth points of the fibres;
\item
$J$ extends the germ $J_0$, and $j$ extends $j_0$.
\end{enumerate}
\end{Defn}

Trying not to make the notation \emph{too} cumbersome, we write $\EuScript{J}(E,\pi)$ for the space of adapted pairs, and $\EuScript{J}(E,\pi,j) =\{ J: (J,j)\in \EuScript{J}(E,\pi)\}$. With respect to the splitting $TE^* = \Tv E^* \oplus \pi^* TS^*$ defined by the symplectic connection (or, in fact, any other connection), a pair $(J,j)\in \EuScript{J}(E,\pi)$ has a block decomposition over $E^*$ of shape
\begin{equation}  \label{J}
J = \left(
\begin{array}{cc}
J^{vv} & J^{vh} \\
0 & j,
\end{array}
\right), \quad  J^{vv} \circ J^{vh} + J^{vh} \circ j =0. 
\end{equation} 
With $J^{vv}$ and $j$ fixed, $J^{vh}$ is just a $\C$-antilinear homomorphism: 
\[J^{vh}\in \Gamma \Hom^{0,1}(\Th E^*, \Tv E^* ).\]

\begin{Lem}
The spaces $\EuScript{J}(E,\pi,j)$ are contractible.
\end{Lem}
\begin{pf}
To see that $\EuScript{J}(E,\pi,j) \neq \emptyset$, fix a closed neighbourhood $U$ of $\crit(\pi)$ on which $J_0$ can be defined; let $U^*=U \setminus E^{\crit}$. One can extend $J_0|{\Tv U^*}$ to an element $J^{vv}\in \EuScript{J}(\Tv E^*,\Omega|{\Tv E^*})$. Over $U^*$, $J^{vv} \oplus  j$ differs from $J_0$ by an antilinear homomorphism, which can be extended to one defined over $E^*$. This gives rise to a complex structure $J$ of the right sort.

Now consider the restriction map 
$ r\colon \EuScript{J}(E,\pi,j) \to \EuScript{J}(\Tv E^*,\Omega|{\Tv E^*})$. 
Its fibres are affine spaces modelled on the vector spaces $\Gamma_c \Hom^{0,1}(\Th E^*, \Tv E^*)$, where $\Gamma_c$ means sections supported outside $X^{\crit}$. The map $r$ admits a section $s$, and $\EuScript{J}(E,\pi,j)$ deformation-retracts to the contractible space $\im(s)$, hence is contractible.
\end{pf}

An easy matrix calculation, linearising (\ref{J}), shows that $\EuScript{J}(E,\pi)$ has formal tangent spaces $T_J \EuScript{J}(X,\pi,j)$ which fit into short exact sequences
\[ 0\to \Hom^{0,1}( (\pi^* TS,j), (\Tv E,J^{vv})) \to   T_J \EuScript{J}(E,\pi,j) \to T_{J^{vv}}\EuScript{J}(\Tv E)\to 0.  \]
Moreover, $\EuScript{J}(E,\pi,j)$ can be made into a smooth Fr\'echet manifold which fibres smoothly over $\EuScript{J}(\Tv E)$.

The drawback of working with spaces of almost complex structures of class $C^\infty$ is that they are not complete, which means that in establishing genericity results, one cannot 
directly invoke the implicit function and Baire category theorems. This is not a serious problem: one fix, due to Floer, is to observe that for any rapidly-decreasing sequence 
$\epsilon=(\epsilon_n)_{n\geq 1}$ of positive reals there is a dense Banach submanifold $\EuScript{J}(E,\pi,j)_\epsilon\subset \EuScript{J}(E,\pi,j)$ whose 
tangent vectors are bundle maps with finite `$C^\infty_\epsilon$-norm'. We refer to Schwarz \cite[Section 4.2.1]{Sch} for the details. 

\subsection{The moduli space} 
Now suppose that $(E,\pi,\Omega,J_0,j_0)$ is a symplectic Morse--Bott fibration over a matched surface, with Lagrangian matching condition $Q$ (Definitions \ref{matching} and \ref{lmc}).

\begin{itemize}
\item
Write $\sect(E)$ for the space of $C^\infty$ sections of $\pi$, and $\sect(E,Q)$ for the subspace of sections $u$ which `map $\partial S$ into $Q$.' This is an abbreviated way of saying that the image of the section $(u|\partial S_+,\tau^*u|\partial S_-)$ of the fibre product $\partial E_+ \times_{\partial S_+} \tau^*\partial E_-$ lies in $Q$.
\item
The tangent space of $\sect(E,Q)$ at $u$ is the subspace $C^\infty( u^*  \Tv E ; \Tv Q)$ of sections $\xi \in C^\infty(u^*\Tv E)$ such that $(\xi|\partial S_+,\tau^* \xi|\partial S_-)$ is a section of $\Tv Q$.
\item
We need a Sobolev completion $\sect_1^p(E,Q)$ of $\sect(E,Q)$, defined using a Riemannian metric and a fixed number $p>2$. It is the smallest subset of the continuous sections $\sect_{C^0}(E,Q)$ which contains $\sect(E,Q)$ and which also contains the section $x\mapsto \exp_{u(x)} (\xi(x))$ for any $u\in \sect(E,Q)$,  $\xi \in L_1^p(u^* \Tv E; \Tv Q)$. The space $\sect_1^p(E,Q)$ is a smooth Banach manifold with tangent spaces $L_1^p(u^* \Tv E;\Tv Q)$.
\item
Fix $(J,j)\in \EuScript{J}(E,\pi, \Omega)$. The moduli space of $(j,J)$-pseudo-holomorphic sections with boundary in $Q$ is the space
\begin{equation}\label{def moduli space}
 \EuScript{M}_{J,j}(E,Q)  =  \{ u \in \sect_1^p(E,Q) :  J \circ (D u) = (D u) \circ j\}\subset \sect_1^p(E,Q), \end{equation} 
and it is this that we wish to analyse. 
\item
Let $\EuScript{E}_{J,j} ^Q\to \sect(E,Q)$ be the natural infinite-rank vector bundle with fibres $(\EuScript{E}_{J,j}^Q )_u=L^p (\Hom^{0,1}_Q(TS, u^*\Tv E))$. Here $\Hom^{0,1}_Q(TS, u^* \Tv E)$ is the vector space of $(j,J^{vv})$-antilinear homomorphisms $TS\to u^* \Tv E$ which (in the obvious sense) carry $T(\partial S)$ to $u^* \Tv Q$.
\end{itemize}

The moduli space $\EuScript{M}_{J,j}(E,Q)$ is the zero-set of the section
\begin{equation}\label{dbar}
 \dbar_{J,j} =\frac{1}{2}(D + J\circ D \circ j) \in \Gamma( \EuScript{E}_{J,j}^Q). \end{equation}  
Allowing $J$ to vary within $\EuScript{J}(E,\pi,j)$ one gets a `universal' vector bundle 
\[\EuScript{E}_j^Q \to \sect(E)\times \EuScript{J}(E,\pi;j).\] 
This has a section $\dbar_j$, whose zero-set $\EuScript{M}_j(E,Q) = \bigcup_J{\EuScript{M}_{J,j}(E)}$ is the `universal moduli space'. 

When $u\in \EuScript{M}_{J,j}(E,Q)$, there is an intrinsically defined linearised operator 
\[D_u\colon L^p_1(u^*\Tv E;Q )\to (\EuScript{E}_{J,j}^Q)_u. \] 
When $u\in \sect_1^p(E,Q)$ is not holomorphic, the linearisation is not intrinsic, but can be defined by choosing a connection on $\EuScript{E}^Q_{J,j}$. Then one has 
\begin{equation} \label{CR}
D_u v = \nabla_u^{0,1} v + a_u v
\end{equation}
for connections $\nabla_u$ on the complex vector bundles $u^*\Tv E$, and bundle maps $a_u\colon  u^*\Tv E \to \Hom^{0,1}(TS, E)$.

\subsection{Transversality}
\emph{Apart from minor differences of context, this is a review of standard theory. We follow Seidel \cite{Se5} closely, since our set-up is almost the same as his.}

The first point to make is that a $C^1$ section of $\pi$ \emph{cannot intersect $E^{\crit}$}, because,  along the image of the section, its derivative provides a right inverse to $D\pi$. Consequently the presence of critical points makes no difference to transversality theory for moduli spaces of sections. It does, however, affect their compactifications. 

Introduce a torsion-free connection $\nabla$ on $T(E\setminus E^{\crit} )$ which restricts to a connection on the vertical subbundle. The linearisation of $\dbar_{J,j} $ (\ref{dbar}) at $u\in \EuScript{M}_{J,j}(E)$ is a linear map
\begin{equation}\label{lin dbar}
 D_u \colon \Gamma( u^* \Tv E) \to (\EuScript{E}_{J,j})_u = \Gamma (\Hom^{0,1}(TS, u^*\Tv E)). \end{equation}
Explicitly,
\[  D_u (v) = (u^*\nabla)^{0,1} (v)  + \frac{1}{2} (\nabla_ v J)\circ Du \circ j. \label{lin2}
\]
The linearisation of $\dbar_j$ is $D_{u,J}^{\univ} $, where 
\[D_{u,J}^{\univ}\colon  (v,Y)\mapsto D_u v + \frac{1}{2}Y\circ Du \circ j.\]

The linearised operators extend continuously to maps between the Banach completions of their domains and targets. Thus we have
\begin{align*}
 D_u \colon & L^p_1 (u^* \Tv E, \Tv Q ) \to   L^p ( \Hom_Q^{0,1} ( TS, u^*\Tv E), \\
   D_{u,J}^{\mathrm{univ}}  \colon &  L^p_1 (u^* \Tv E, \Tv Q )  \oplus C^\infty_\epsilon(T_J\EuScript{J}(E,\pi,j)) 
  \to  L^p ( \Hom_Q^{0,1} ( TS, u^*\Tv E). 
\end{align*} 

The crucial point is that $D_u $ is Fredholm (we will come back to this point shortly). It follows that so too is $ D_{u,J}^{\univ}$. The latter is also \emph{surjective}, whence $\EuScript{M}_j(E,Q)$ is a smooth Banach submanifold of $ \sect^p_1(E,Q)\times \EuScript{J}(E,\pi;j)_\epsilon$. This surjectivity statement is non-standard only in one respect, namely, that we are considering almost complex structures on $\partial E_+\times_{\partial S_+}\tau^*\partial E_- $ which respect the fibre product decomposition. However, the standard argument, involving unique continuation, is unaffected by this point (if $D_{u}^*\eta=0$, where $\eta$ is an $L^q$ section of the dual bundle to $\Hom_Q^{0,1} ( TS, u^*\Tv E)$, $q^{-1}+p^{-1}=1$, then $\eta$ must be supported in $\partial S$; but $L^q\cap \ker D_u^*$ contains a dense subspace of continuous sections, hence $\eta=0$). 

Write $\EuScript{J}^{\reg}(X,\pi;j)$ for the space of {\bf regular} almost complex structures: those with the property that $D_u$ is onto for every $u\in \EuScript{M}_{J,j}(X,Q)$. Their importance is  that when $J$ is regular, $\EuScript{M}_{J,j}(E,Q)$ is a smooth manifold of local dimension $\ind D_u$. The key transversality statement, following from surjectivity of $D_{u,J}^{\mathrm{univ}}$, is that $\EuScript{J}^{\reg}(E,\pi,  j)\cap \EuScript{J}(E,\pi,J)_\epsilon$  is a (dense) Baire subset of the complete space $\EuScript{J}(E,\pi, j)_\epsilon$ (as trailed earlier, we use Floer's $C^\infty_\epsilon$ spaces here).

\subsubsection{Intersection with cycles in fibres} Mark a finite set of points $\{s_i\}_{i\in I}$ in $\mathrm{int}(S)$, manifolds $Z_i$, and smooth maps $\zeta_i \colon Z_i\to X_{s_i}$. There is an evaluation map 
\begin{equation}\label{eval}
\ev_I^J =\prod_{i\in I} {\ev_i} \colon \EuScript{M}_{J,j}(X,Q) \to \prod_i {X_{s_i}}.\end{equation} 
The arguments that establish the surjectivity of $D_{u,J}^{\univ}$ extend easily to give a `transversality of evaluation' lemma (see \cite[Lemma 2.5]{Se5}) which says that, for any $J \in \EuScript{J}(X,\pi;j)$, and any neighbourhood $U$ of $\{s_i:i\in I\}$, there exist arbitrarily small perturbations $J_t$ of $J$, supported in $\pi^{-1}(U)$, such that $J_t\in \EuScript{J}^{\reg}(X,\pi;j)$ and $\ev_I^{J_t}$ is transverse to $\zeta = \prod_i{\zeta_i}$. Similarly, given marked points $s_k' \in \partial S_+$, and maps  $\zeta'_k\colon Z_k' \to Q_{s_k'}$, there exist small perturbations of $J$ supported in a chosen neighbourhood of $\{s_k'\} \cup \{\tau^{-1}(s_k')\}$ making $\zeta'_k$ transverse to the natural evaluation map. (Again, the argument is unaffected by our use of almost complex structures which respect the fibre product structure.)

\subsubsection{Fredholm theory} 
It is a simple matter to prove that pseudo-holomorphic sections with Lagrangian matching conditions have a Fredholm deformation theory, granted the standard Fredholm theory for Lagrangian boundary conditions. The argument is most direct in the case where $S= S_+ \amalg S_- $ and the matching diffeomorphism $\tau$ of the boundary components extends to a diffeomorphism $\tau'\colon S_+ \to S_- $: here, a Lagrangian matching condition is simply a Lagrangian boundary condition in a fibre product.

In general, the argument goes as follows. Near $\partial S$, we can extend $\tau$ to an diffeomorphism $ \widetilde{\tau}$ between collar neighbourhoods $C_+$ of $\partial S_+$ and $C_- $ of $\partial S_-$. 
Take $u\in \sect(E,Q)$, and let $u_\pm = u|C_\pm$. Then we have a section $(u_+, \widetilde{\tau}^* u_-)$ of the fibre bundle $E|C_+\times_{C_+} \widetilde{\tau}^* (E|C_-)\to C_+$. 

We claim that the map (\ref{lin dbar}),
\[D_u\colon L^p_1(u^*\Tv E;Q )\to (\EuScript{E}_{J,j}^Q)_u, \] 
is Fredholm. Over compact subsets of $\mathrm{int}(S)$, one has the usual elliptic estimate for $D_u$. Near the boundary, one can equivalently work with $D_u$ over $C_+\cup C_-$ or with the corresponding operator  $D_{(u_+, \widetilde{\tau}^* u_-)}$ over $C_+$. The theory developed by Floer---see \cite[Lemma 2.3]{Flo}---gives an elliptic estimate for the latter. Together, these estimates imply that $\ker(D_u)$ is finite-dimensional. On small open sets in $\mathrm{int}(S)$ or in $C_+$, the operator $D_u$ has a bounded right inverse, by elliptic theory, and these may be patched together to form a global parametrix (as in \cite[Prop. 3.6]{Do2}), which shows that $D_u$ has closed range and finite-dimensional cokernel.

To state the index formula for $D_u$, we need to define a `Maslov index' map  
\[\mu_Q \colon \sect(E,Q)\to \Z.\] 
For a component $\partial S_+^{(i)}$ of $\partial S_+$, let $v_+^{(i)} = u|\partial S_+^{(i)}$. Let $v_-^{(i)} = u|\tau (\partial S_+^{(i)})$. Choose symplectic trivialisations $t_+^{(i)}$ for $(v_+^{(i)})^* \Tv E$ and $t_-^{(i)}$ for $\tau^* (v_-^{(i)})^*\Tv E$. Then the pullback $(v_+^{(i)}, \tau \circ v_-^{(i)})^*\Tv Q $ defines (via $t_+^{(i)}\oplus t_-^{(i)}$) a loop in a Lagrangian Grassmannian, and so has a classical Maslov index $l_i$. We put
\[ \mu_{Q}(u)  = 2c_1^{\mathrm{rel}}(u^* \Tv E ) + \sum_i{l_i} , \] 
where $c_1^{\mathrm{rel}}(u^*\Tv E)$ is the Chern number of $u^* \Tv E\to S$ relative to the trivialisations $t_\pm^{(i)}$ over $\partial S$. 
\begin{Lem}\label{matched index}
The index of $D_u$ is 
\[ \ind(D_u) = \mu_Q(u) +  \sum_{S_i \in \pi_0 (S)}{\mathrm{rk}_{\C}(\Tv E|S_i)\chi(S_i)}. \]
\end{Lem}
\begin{pf}
We use the observation from \cite{EGH} that the index of Cauchy--Riemann operators is invariant under `cutting and pasting'. For a positive boundary component $\partial S_+^{(i)}$, let $\gamma_+^{(i)}$ be a loop in $S'$, parallel to $\partial S_+^{(i)}$, and let $V_+^{(i)}=(u^* \Tv E)|\gamma_+^{(i)}$. Let $\delta\subset S^2$ be the equator, and let $U\to S^2$ be a complex vector bundle with $c_1(U)=0$ and $U|\delta=V_+^{(i)}$. Cut $S'$ along $\gamma_+^{(i)}$, and $S^2$ along $\delta$ to obtain a Riemann surface with boundary, then re-glue the four boundary components so as to obtain a new Riemann surface (the result of surgery on $\gamma_+^{(i)}$) which inherits a vector bundle. By \cite[1.8, `Axiom C3']{EGH}, this vector bundle also inherits a Cauchy--Riemann operator $D_{u,+}^{(i)}$, and $\ind(D_{u,+}^{(i)})=\ind (D_u) + 2 \rank_{\C}(V_+^{(i)})$. (According to \cite{EGH} one should consider real analytic data here; however, ultimately this will make no difference.)

Perform this surgery operation for each $i$; also perform similar surgeries on a loop $\gamma_-^{(i)}$ parallel to $\tau(\partial S_+^{(i)})$, for each $i$. The result is (i) a matched Riemann surface $C$ which is the union of a closed Riemann surface $C'$ and $k$ pairs of discs $(D_+^{(i)}, D_-^{(i)})$ with matched boundaries; (ii) a Cauchy--Riemann operator $D_u'$ over it. The ordinary Riemann--Roch theorem computes the index of $D_u'$ over $C'$. Riemann--Roch for surfaces with boundary computes the index over $D_+^{(i)}\cup D_-^{(i)}$ as $ l_i + \rank(\Tv E)|D_+^{(i)} + \rank(\Tv E)|D_-^{(i)}$, where $l_i$ is a Maslov index as in the definition of $\mu_Q(u)$. Putting these things together gives the result.
\end{pf}

\subsection{Compactness}
The {\bf action} of a section $u\in \sect(E,Q)$ is 
\[\EuScript{A}(u) = \int_S{u^*\Omega}. \]
Though $\Omega$ itself need not be symplectic, $\Omega+c\pi^*\beta$ is symplectic (tamed by $J$) for any positive area form $\beta\in \Omega^2_S$ and $c \gg 0$.
The symplectic area for sections then differs from the action by a constant, $c\int_S{\beta}$. Gromov's compactness theorem (the relevant version is that of Ye \cite{Ye}) says that any sequence $u_1,u_2,\dots$ in $\EuScript{M}_{J,j}(X,Q)$ of bounded symplectic area---equivalently, bounded action---has a subsequence which converges, in Gromov's topology, to a pseudo-holomorphic curve $v$ with bubbles and boundary bubbles. 

The principal component of a Gromov-limit is a differentiable map $S\to E$ which is a section over the complement of a discrete set in $S$ (the roots of the bubble-trees), hence is globally a section.
Hence to compactify $\{u\in \EuScript{M}_{J,j}(E,Q): \EuScript{A}(u)\leq \lambda\}$ one need only consider curves whose components are (i) sections; (ii) bubbles in regular fibres; (iii) bubbles in singular fibres; (iv) boundary bubbles.

\subsubsection{Fibred monotonicity} One does not expect moduli spaces of  pseudo-holomorphic sections with Lagrangian boundary conditions to resemble  
closed manifolds in general, even to the extent of carrying fundamental homology classes.  The reason is that there is typically a codimension 1 boundary, corresponding to bubbling off of discs. 
Various hypotheses can be imposed to make sure that there is in fact no such boundary in codimension 1. Ours will be fibred weak monotonicity; indeed, we will mostly be concerned with fibred monotonicity. In our applications, which will involve symmetric products $\sym^n(\Sigma)$, fibred monotonicity holds when $n\geq g(\Sigma)$. 

\begin{Defn}
(a) A closed Lagrangian submanifold $L$ in a compact symplectic manifold $(M,\omega)$ is called {\bf monotone} if there exists $c>0$ such that 
\[  \int_{\bar{D}} u^*\omega =  c\,\mu(u)\]
for every smooth map $u \colon (\bar{D},\partial\bar{ D})\to (M,L)$.

(b) A symplectic Morse--Bott fibration $(E,\pi,\Omega)$ over a matched surface, with Lagrangian matching condition $Q$, is {\bf fibre-monotone} if, for each $ s\in \partial S_+$, $Q_s$ is monotone in $(\partial E_+ \times_{\partial S_+}\tau^*\partial E_-)_s$.
\end{Defn}

For a Lagrangian submanifold $L\subset M$, let $\mu_L\colon \pi_2(M,L)\to \Z$ be the Maslov index homomorphism and $\mu_{\mathrm{min}}(L)$ the minimal Maslov index:
\begin{equation}\label{mumin}  
 \im(\mu_L) = \mu_{\min}(L) \Z,\quad \mu_{\mathrm{min}}(L)\geq 0. 
\end{equation}
Let $c_{\min}(M)$ be the minimal Chern number of $M$:
\begin{equation}\label{cmin}  
  \im(c_1\colon  \pi_2(M) \to \Z) = c_{\min}(M) \Z,\quad  c_{\min}\geq 0.  
  \end{equation}
Since the Maslov index of a disc whose boundary is mapped to a point is twice its Chern number, $\mu_{\min}(L)\leq 2c_{\min}(M) $. For a Lagrangian boundary condition $Q$, let $\mu_{\mathrm{min}}(Q)=\gcd_{s\in \partial S}{\mu_{\min}(Q_s)}$.

In the following proposition, we suppose that $E^{\crit}=\emptyset$, so that $(E,\pi,\Omega)$ is actually a locally Hamiltonian fibration (LHF). 
\begin{Prop}\label{mod sp}
Let $Q$ be a Lagrangian matching condition for the LHF $(E^{2n+2},\pi,\Omega)$ over a matched surface $S$, and fix $h\in \pi_0 \sect(E,Q)$. 
Fix also finite sets $I \subset \mathrm{int}(S)$, $ J\subset \partial S_+$; and smooth cycles $\zeta_s \colon Z_s \to E_{s}$ for each $s\in S$, and $\zeta' \colon Z_{s'} \to  Q_{s'}$ for $s'\in J$. Let $Z=\prod{Z_s}$, $\zeta = \prod{\zeta_s} \colon Z\to \prod{E_s}$; and  $Z'=\prod{Z_s'}$, $\zeta' = \prod{\zeta_s'} \colon Z\to \prod{Q_{s'}}$. Suppose that
\begin{enumerate}
\item
Both $E$ and $Q$ are fibre-monotone;
\item
the virtual dimension $d(h)=\ind(D_u)$ ($u\in h$) satisfies 
\begin{equation}\label{dim cond}
 d(h) - \sum_{s\in I}{(2n-\dim(Z_s))}-\sum_{s'\in J}{(n-\dim(Z_{s'}'))} - \mu_{\min}(Q) < 0;
\end{equation}
\end{enumerate}
Then for any $j\in \EuScript{J}(TS)$ and for a dense set of $J\in \EuScript{J}^{\reg}(X,\pi,j)$, the fibre product
\[ \EuScript{M}_{J,j}(E,Q) \times_{(\ev ,\zeta,\times\zeta')}(Z\times Z').  \]
is a compact, smooth manifold of of dimension $d(h) - \sum_I{(2n-\dim(Z_s))}-\sum_J{(n-\dim(Z_{s'}'))}$ when this number is non-negative, and empty when it is negative. 
\end{Prop}
 
\begin{pf}
By the remarks on transversality of evaluation (\ref{eval}), there is a dense subset of $\EuScript{J}^{\reg}(X,\pi,j)$ whose members $J$ have the property that the evaluation map 
\[ \ev= \ev_A \times \ev_B' \colon \EuScript{M}_{J,j}(X,Q)\to \prod_a E_{s_a}\times\prod_b{ Q_{s_b'}}\] 
is transverse to $\zeta \times \zeta'$. This means that there is a smooth moduli space 
\[  \EuScript{M} = \EuScript{M}_{J,j}(E,Q) \times_{(\ev ,\zeta,\times\zeta')}(Z\times Z').  \]
Any sequence in $\EuScript{M}$, for which the sections all lie in $h$, has a subsequence with a Gromov limit comprising a principal component $u_{\mathrm{prin}}$, non-constant bubbles $\{u_l\}_{l\in L}$, and non-constant boundary bubbles $\{u_k'\}_{k\in K}$. These satisfy
\begin{equation}\label{index sum}
 d(h) = \ind(u_{\mathrm{prin}}) + 2\sum_{l\in L}{c_1(u_l)}+ \sum_{k\in K}{\mu(u_k')}. 
 \end{equation}
By fibred monotonicity, the Maslov indices $\mu(u_k')$ are positive multiples of $\mu_{\min}(Q)$, hence $\geq \mu_{\min}(Q)$. The contribution of each term $c_1(\beta_l)$ from a bubble in a fibre is also positive, at least $\mu_{\min}(Q)$. Then if $L \cup K \neq \emptyset$ we will have $\ind(u_{\mathrm{prin}})- \sum_{a}{(2n-\dim(Z_a))}-\sum_{b}{(n-\dim(Z_b'))}<0$. But for regular almost complex structures, this number gives the dimension near $u_{\mathrm{prin}}$
of $\EuScript{M}$, which cannot be negative.
\end{pf}

There are many possible hypotheses to make singularities in the fibres permissible, and we do not attempt an axiomatisation. In the case of relative Hilbert schemes of points on Lefschetz fibrations, they can be admitted because of the following observation of Donaldson and Smith:
\begin{Lem}\cite[Lemma A.11]{DS}
Let $\hilb^n_S(E)$ be the relative Hilbert scheme of $n$ points on a Lefschetz fibration $E\to S$ with irreducible fibres, and $\hilb^n(E_s)$ a singular fibre. Then any holomorphic sphere in $\hilb^n(E_s)$ which arises as a bubble-component in a Gromov limit of pseudo-holomorphic sections is homotopic to a sphere in a regular fibre.
\end{Lem}

\subsubsection{Invariance} We now give the parametric version of Proposition \ref{mod sp}. All the data involved in the construction can be allowed to move in smooth families. To begin with, \emph{fix} the symplectic Morse--Bott fibration $(E,\pi,\Omega,J_0,j_0)$, Lagrangian boundary condition $Q$, and homotopy class $h$. 

Each $j$ and each $J\in\EuScript{J}^{\reg}(E,\pi; j )$ gives a moduli space $\EuScript{M}_{J,j}(E,Q;)_h$ which we know to be a smooth, oriented manifold; after a small perturbation of $J$, the fibre product $\EuScript{M}_{J,j}(E,Q;)_h\times_{\ev,\zeta\times\zeta'} Z\times Z'$ with a family of smooth oriented cycles in the fibres (as considered above) is also smooth and oriented. It carries a smooth boundary-evaluation map 
\[ \ev(J,j) \colon \EuScript{M}_{J,j}(E,Q)_h\to \sect(Q) \]
to the space of sections of $Q$. 

\begin{Prop} 
Assume that the hypotheses of Proposition \ref{mod sp} hold, and that the left-hand side of (\ref{dim cond}) is $< -(k+1)$, for some $k\geq 0$. Then the oriented bordism class of the map $\ev(J,j)$ is independent of $(J,j)$. Moreover, any smooth map $S^k \to \EuScript{J}^{\reg}(X,\pi)$ induces a bundle $\EuScript{M}_{S^k}$ over $S^k$ and a map $\ev \colon \EuScript{M}_{S_k}\to \sect(Q)$ which extends to a map to $\sect(Q)$ from a bundle over $B^{k+1}$. 
\end{Prop}

\begin{pf}
One considers parametrised moduli spaces. The space $\EuScript{J}(E,\pi)$ is contractible; thus any map $f\colon S^k\to  \EuScript{J}^{\mathrm{reg}}(X,\pi)$ extends to a map $F\colon B^{k+1}\to \EuScript{J}(X,\pi)$. The transversality theory extends in a straightforward way, showing that after a small homotopy of $F$, fixing $f$, one gets a smooth moduli space of the appropriate dimension which inherits an orientation from that on $B^{k+1}$. The numerical conditions imply that compactness goes through as before, and this gives the result. 
\end{pf}

Likewise, the bordism class of the evaluation map is unchanged under deformations of $\Omega$, $J_0$, $j_0$, the points $s_a$, $s_b'$ and the cycles $\zeta_a$ and $\zeta_b'$. 

\subsection{Fibred \emph{k}--negativity}

In the framework of fibred monotonicity, bubbling in low-dimensional moduli spaces (with regularity assumptions in force) is ruled out on the grounds that the principal component of a limit curve must have non-negative index. A more refined method incorporates transversality for pseudo-holomorphic spheres in the fibres and discs in the boundary fibres. One can then determine numerical conditions under which moduli spaces of sections are generically disjoint from those of fibrewise spheres and discs.
In the context of Hamiltonian Floer homology, these are the `weak monotonicity' conditions of Hofer and Salamon \cite{HoS} (a symplectic manifold $(M^{2n},\omega)$ is weakly monotone if it is monotone, \emph{or} $c_{\min}=0$, \emph{or} $c_{\min} \geq n-2$). To handle Lagrangian boundary conditions effectively, somewhat more conservative assumptions are required. We have seen that monotone Lagrangians in monotone symplectic manifolds are acceptable, and we will consider one other hypothesis, drawing on technical results of Lazzarini.\footnote{There is, unfortunately, no straightforward general theory for Lagrangians with $\mu_{\min}=0$ in symplectic manifolds with $c_{\min}=0$ (these are the characteristics of vanishing cycles in $\sym^{g-1}(\Sigma)\times \sym^{g-2}(\bar{\Sigma})$). This is a borderline situation in which there is just enough leeway to obtain compact zero-dimensional moduli spaces, but not enough to handle one-parameter families. We shall not venture into this hazardous territory.}

\begin{Defn}\label{kneg}
For any $k\geq 0$, a compact symplectic manifold $(M^{2n},\omega) $ is {\bf $k$--negative} if, for any sphere $f\colon S^2\to M$ with $\int_{S^2}{f^*\omega}>0$, one has $\langle f^* c_1(TM), [S^2]\rangle \leq 3-n-k$. A Lagrangian $L\subset M$ is called {\bf $k$--negative} if any disc $f\colon (D,\partial D)\to (M,L)$ with $\int_{D}{f^*\omega}>0$ has Maslov index $\leq 2-n-k$. 
\end{Defn} 
We shall need to consider $3$--negative manifolds and $2$--negative Lagrangians.

\subsubsection{Transversality in the fibres} For $k$--negative manifolds ($k\geq 0$) generic almost complex structures $J\in \EuScript{J}(M,\omega)$ are regular in the sense that for every element $\beta\in \EuScript{M}_J^{\mathrm{si}}(M)$ of the moduli space of parametrised, \emph{simple} $J$-holomorphic spheres, the relevant linearised operator $D_\beta$ is surjective (a pseudo-holomorphic curve $u\colon C \to M$ is \emph{simple} if it does not factor through a branched covering $C\to C'$ of degree $>1$; by a lemma of McDuff (see \cite[Proposition 2.51]{MS}) this is equivalent to the set of injective points, $\{z\in S^2: \beta^{-1}(\beta(z))=\{z\}\}$, being non-empty, or indeed open and dense). For such $J$, $\EuScript{M}_J^{\mathrm{si}}(M)$ is a manifold of local dimension $2(  c_1(\beta) + n)$, and its quotient by the free action of $\aut(S^2) = \mathrm{PSL}(2,\C)$ has dimension $2( c_1(\beta) + n-3)$. Thus there can be no simple  holomorphic spheres. Since every pseudo-holomorphic sphere factors through a simple one, there are no non-constant spheres at all. By the same argument, one sees that \emph{in a $k$--negative symplectic manifold, there are are no non-constant spheres in a generic $k$--parameter family of almost complex structures.} 

As with spheres, transversality theory for pseudo-holomorphic discs in $(M,L)$ works in a straightforward way for \emph{simple} discs---those for which the set of injective points is non-empty, or equivalently open and dense---but the complication is that not every non-constant disc factors through a simple one. One can, however, use Lazzarini's lemma \cite{Laz}: \emph{if $\delta$ is a non-constant $J$--disc then there is a simple $J$--disc $\delta'$ such that $\delta'(\partial \bar{D})\subset \delta(\partial \bar{D})$.} The index formula says that the moduli space of parametrised simple discs $\EuScript{M}^{\mathrm{si}}_J(M,L)$ has virtual dimension $\mu_L+ n$. Since $\aut(\bar{D})$ is 3-dimensional and acts freely, the unparametrised moduli space has dimension $\mu_L+n-3$. Thus, for generic $k$--parameter families of almost complex structures, the moduli space of simple discs is empty. Lazzarini's lemma allows us to remove `simple' from this conclusion.

Since we are interested in fibrations over surfaces, with Lagrangian boundary conditions, we should consider $3$--negative fibres and $2$--negative Lagrangians in them. (Remember that we need a spare parameter to prove invariance!)

\subsubsection{Singular fibres} Symplectic Morse--Bott fibrations can be handled using the method developed for Lefschetz fibrations in Seidel's thesis \cite{Se2}. One needs to make some hypothesis about the normalisation $\widetilde{M}_0$ of a singular fibre $M_0$. The almost complex structure $J_0=J|M_0$ on $M_0$ lifts to one $\widetilde{J}_0 $ on $\widetilde{M}_0$, integrable near the preimage of the normal crossing divisor; it is sufficient to assume that there $\widetilde{J}_0$ is compatible with a $1$--negative symplectic form on $\widetilde{M}_0$. The point is that any $J_0$-sphere in $M_0$ lifts uniquely to a $\widetilde{J}_0$-sphere in $\widetilde{M}_0$. One cannot expect $\widetilde{J}_0$ to be regular; however, for generic $J$, lifts of non-constant spheres in $M_0$ \emph{are} regular, hence have non-negative Chern number: here we use a standard argument which shows that one can achieve regularity for all spheres passing through an open set $U$ by making a perturbation of the almost complex structure supported in $U$. 

Let us summarise our conclusion:
\begin{Prop}
Let $Q$ be a Lagrangian matching condition for a symplectic Morse--Bott fibration $(E^{2n+2},\pi,\Omega)$ over a matched surface $S$, and fix $h\in \pi_0 \sect(E,Q)$. Fix also finite sets of regular values $I \subset \mathrm{int}(S)$ and $ J\subset \partial S_+$; and smooth cycles $\zeta_s \colon Z_s \to E_{s}$ for each $s\in S$, and $\zeta'_{s'} \colon Z_{s'} \to  Q_{s'}$ for $s'\in J$. Let $Z=\prod{Z_s}$, $\zeta = \prod{\zeta_s} \colon Z\to \prod{E_s}$; and  $Z'=\prod{Z_s'}$, $\zeta' = \prod{\zeta_s'} \colon Z\to \prod{Q_{s'}}$. Suppose that
\begin{enumerate}
\item
the fibres of $E$ are 3--negative; those of $Q$ are 2--negative;
\item
the components of $E^\crit$ are 1--negative;
\item
the number
\[ D(h ; Z_i,  Z_j ):=n\chi(S) + \mu_Q(h) - \sum_{s\in I}{(2n-d_s(Z))}-\sum_{s'\in J}{(n-d_{s'}'(Z'))} \]
is non-positive.
\end{enumerate}
Here $d_s(Z)=\dim Z_s$, $d'_{s'}(Z') =\dim Z_{s'}$. Then, for any $j\in \EuScript{J}(TS)$ and for a dense set of $J\in \EuScript{J}^{\reg}(E,\pi,j)$, the fibre product
\[ \EuScript{M}_{J,j}(E,Q) \times_{(\ev ,\zeta\times\zeta')}(Z\times Z').  \]
is a compact, smooth manifold of of dimension $D(h; Z,  Z')$ when this number is negative, and empty when it is zero. 
Moreover, for generic paths in $ \EuScript{J}^{\reg}(E,\pi,j)$, the one-parameter moduli space is a compact 1--manifold with boundary when $D(h; Z_i,  Z_j' )=0$, and empty when $D(h; Z,  Z' )<0$.
\end{Prop}

\begin{Rk} The same transversality arguments also apply when the fibres, Lagrangians, and critical manifolds are all monotone; this is a more complicated approach to the one discussed above, but gives a slightly sharper result: one need only assume that $D(h)\leq 0$. Note also that this technique handles the singular fibres of $\hilb^n_S(E)$ for $n\geq g$ in the case previously neglected, when $E$ has reducible singular fibres. 
\end{Rk} 

\subsection{Orientations} \label{orientations}

Since the first version of this paper was written, a draft of Seidel's book on Fukaya categories \cite{Se6} has become available. It includes a thorough  account of orientations in Lagrangian Floer theory incorporating non-orientable Lagrangians. It seems churlish to ignore this useful reference, particularly as the the draft of this section dealt unconvincingly with one or two points. The account we now present should be regarded, for the most part, as expository;  besides Seidel's book, its sources are De Silva's thesis \cite{Sil} and the book of Fukaya--Oh--Ohta--Ono \cite{FOOO}.

Consider a space $\EuScript{M}$ defined as the zero-set of a non-linear Fredholm map $f$ between Banach spaces. Suppose, moroever, that the linear Fredholm maps $F_x$ linearising $f$  are all surjective (`regularity'). Then $\EuScript{M}$ is a smooth manifold with tangent spaces $\ker(F_x)$. Orientability of this manifold means triviality of $\Lambda^{\max}T \EuScript{M}$, the line bundle with fibres $\Lambda^{\max} \ker(F_x)$. In practice, $\EuScript{M}$ sits inside a larger space $\EuScript{B}$ parametrising a family of Fredholm operators $\EuScript{F} =\{ F_x\}_{x\in \EuScript{B}}$, and  $\Lambda^{\max}T \EuScript{M}$ extends to the larger space as  the determinant index bundle $\mathrm{Det}(\EuScript{F}) \to \EuScript{B}$, the natural line bundle with fibres  $\Lambda^{\max} \ker(F_x)\otimes \Lambda^{\max}(\coker (F_x))^*$. One then studies $w_1$ of this bundle. 

The Fredholm operators we wish to consider are Cauchy--Riemann operators over a Riemann surface with boundary, subject to a Lagrangian boundary condition.  Let $S$ be a compact Riemann surface with one boundary component. Let $E\to S$ be a (trivial) symplectic vector bundle of rank $2n$. We consider the space $\Lag(S,E)$ of Lagrangian subbundles $F\subset  E| \partial S$. By fixing a diffeomorphism $\partial S\cong S^1$ and a trivialisation of $E|\partial S$, we can identify $\Lag(S,E)$ with the space of unbased loops $\EuScript{L} \Gr_n$ in $\Gr_n$, the Grassmannian of Lagrangian subspaces of the standard symplectic vector space $(\C^n,\omega_{\C^n})$. Notice that the Maslov index decomposes $\Lag(S,E)$ (or equally $\EuScript{L}\Gr_n$) into 
connected components $\Lag_k(S,E)$ (or $\EuScript{L}_k\Gr_n$) indexed by $k\in \Z$.

By choosing a conjugate-linear bundle map $a \in V:=\Hom^{0,1}(TS, E)$, we get an $\R$-linear, Fredholm Cauchy--Riemann problem for $(E,F)$: that of solving $ (\dbar  + a) (u) =0 $, where $u \in L_1^p(E)$, for some fixed $p>2$,
and $u(z)\in F_z$ for $z\in \partial S$. Our object of study is the determinant index bundle $\Det(\EuScript{F})\to V\times \EuScript{L}\Gr_n$. 

The space $ V\times \EuScript{L}\Gr_n$ obviously deformation-retracts to the subspace $\{0\}\times \EuScript{L}\Gr_n$. This in turn deformation-retracts to $\U(n)/\orth(n)$, via a deformation retraction of $\Sp(2n,\R)$ to $\U(n)$. Hence $\Lag(S,E)$ is homotopy-equivalent to $\EuScript{L}(\U(n)/\orth(n))$. 
\begin{Rk}
It is worth noting that if one considers only the based loops, and passes to the limit $\U /\orth = \varinjlim_n\U(n)/\orth(n)$, there is a homotopy equivalence 
\[ \Omega (\U/\orth) \simeq  \Z \times B\orth  \]
with the classifying space of $KO$-theory. Bott \cite{Bot} constructed such an equivalence using Morse theory, as the final step in the periodic octagon starting with $\Z\times B\orth$; there is also an equivalence $\Omega^7(\Z\times B\orth)\simeq \mathrm{U}/\orth$. We shall not need Bott periodicity, but we shall be concerned with the functorial map 
\[ [X,\EuScript{L}(U/\orth)] \to  KO(X) \]
defined by stabilising the virtual index bundle, and with its restriction to $\Omega(\U/\orth)$.
\end{Rk}
Following \cite{Se6} (which itself extends \cite{Sil}), we now give the formula for $w_1(\Det(\EuScript{F}))\in H^1(\Lag(S,E);\Z/2)$.  Evaluation $e\colon S^1\times \EuScript{L}_k\Gr_n \to \Gr_n$ induces a map 
\[    e^*\colon H^*(\Gr_n;\Z/2) \to  H^*(S^1\times \EuScript{L}_k\Gr_n;\Z/2). \]
For $c\in H^i(\Gr_n;\Z/2)$, let $T(c)= \int_{S^1}e^*c\in H^{i-1}(\Gr_n;\Z/2)$, and let $U(c) = (e^*c)|( \{z\}\times \Gr_n)\in H^i(\Gr_n;\Z/2)$. It is easy to check that
\[ H^1(\EuScript{L}_k\Gr_n;\Z/2) = \Z/2\oplus \Z/2,\quad n\geq 3,  \]
with generators $T(w_2)$ and $U(\mu)$, $\mu \in H^1(\EuScript{L}_k\Gr_n;\Z/2)$ being the Maslov index reduced mod 2. The first Stiefel--Whitney class for the determinant index bundle over the component $\EuScript{L}_k\Gr_n$ is then 
\begin{equation}\label{w1}
	w_1(\Det(\EuScript{F})) = T(w_2) + (k+1)U(\mu). 
\end{equation}
Given a loop $\gamma\colon S^1\to \Lag(S,E)$, lying in the component where the Maslov index is $k$, one can take the torus of boundary values $\partial \gamma\colon S^1\times \partial S \to L$. What the formula says is that, 
if $k$ is odd, one has
\begin{equation}
\langle w_1(\gamma^* \Det(\EuScript{F})),[S^1]\rangle = 
		\langle (\partial \gamma)^*w_2(TL), [S^1\times \partial S]\rangle.
\end{equation} 
If $k$ is even one must add a correction term  $\langle w_1 l^*TL,[S^1]\rangle$, where $ l$ is the loop $(\partial\gamma) | S^1\times \{x\} \colon S^1 \to L$.
In particular, $ w_1(\gamma^* \mathrm{Det})$ vanishes when $L$ is spin.

Fukaya et al. \cite{FOOO} show that one gets an actual orientation for the moduli space $\EuScript{M}_S(M,L)$ of $J$--holomorphic maps $(S,\partial S)\to (M,L)$ by giving a `relative spin structure', that is, an oriented vector bundle $\xi \to M$, with $w_2(\xi|L) = w_2(TL)$, and a spin structure on $\xi|L \oplus TL$.  This is not quite enough for our purposes, because it precludes non-orientable Lagrangian boundary conditions. We can get around this by working not with spin but with \emph{pin} (again, cf. \cite{Se6}).

The Lie group $\pin^+(n)$ is one of the central extensions of $\orth(n)$  by $\Z/2$. Such central extensions of  topological groups are 
classified by $H^2(B\orth(n);\Z/2)$, which is the direct sum of two copies of $\Z/2$ generated by the universal characteristic classes $w_2$ and $w_1^2$. The central extension associated with $w_2$ (resp. $w_1^2+w_2$) is $\pin^+(n)$ (resp. $\pin^-(n)$).  A well-known algebraic construction assigns to any quadratic vector space $(V,Q)$ a group $\pin(V,Q)$ as the subgroup of the unit group of the Clifford algebra $\mathrm{Cl}(V)$ (in which the relation $v^2 = Q(v) 1$ holds) generated by $Q^{-1}\{\pm 1\} \subset V\subset \mathrm{Cl}(V)$. In these terms, $\pin^\pm (n)= \pin(\R^n, \pm \| \cdot \|^2)$.

A $\pin^+$--structure for a real $n$-plane bundle $\zeta\to Z$ is a homotopy class of homotopy-liftings $Z\to B \pin^+(n)$ of the classifying map $Z\to B \orth(n)$. 
The obstruction to such a lift is $w_2(\zeta)$; when this is zero, the $\pin^+$--structures form an affine space modelled on $H^1(Z;\Z/2)$. 
The inclusions $\orth(n)\to\orth(n+1)$ lift to the $\pin^+$--groups, and this gives the meaning of `stable'. 

\begin{Ex}\label{odd pin}
When $n$ is odd, one has an isomorphism
\[\orth(n)\stackrel{\cong}{\longrightarrow} \Z/2 \times \SO(n) , \quad A\mapsto (\det(A),\det(A)^{-1} \cdot A). \] 
The four central extensions of $\Z/2\times \SO(n)$ by $\Z/2$ are:
\begin{enumerate}
\item 
$\Z/2\times \Z/2\times \SO(n)$, 
\item 
$\Z/4 \times \SO(n)$, 
\item
$\Z/2\times \Spin(n)$, and 
\item $(\Z/4 \times \Spin(n))/ (\Z/2)$ (quotient by the `diagonal' involution). 
\end{enumerate}
Which of these is $\pin^+(n)$? Certainly one of the last two, since the covering $\pin^+(n)\to \orth(n)$ must restrict to the spin covering $\Spin(n)\to \SO(n)$. The algebraic fact \cite[Lemma 11.14]{Se6} which allows us to decide between them is that the two preimages of $-\id\in \orth(n)$ in $\pin^+(n)$ both have order 2 when $n\equiv 1 \mod 4$, and 4 when $n\equiv 3 \mod 4$. Hence $\pin^+(n)$ is $\Z/2\times \Spin(n)$ or $(\Z/4 \times \Spin(n))/ (\Z/2)$ according to whether $n$ is 1 or 3 mod 4. 
\end{Ex}

Now suppose that $(X^{2n+2},\pi,\Omega)$ is a LHF over a compact matched surface $S$, and $Q\subset \partial X$ a Lagrangian matching condition. Thus
$Q$ is a sub-bundle of $Y_+ \times_{\partial S_- } Y_-$, where $Y_-= \pi^{-1}(\partial S_-)$ and $Y_+$ is the pullback to $\partial S_-$ of $ \pi^{-1}(\partial S_-)$ by the matching $\tau\colon \partial S_-\to \partial S_+$. 
\begin{Defn}
A {\bf relative pin structure} for $(X,Q)$ is 
\begin{itemize}
\item
a stable oriented vector bundle $\xi \to X$, with $w_2(\xi|Q)=w_2(\Tv Q)$; and
\item
a stable $\pin^+$--structure on $ \xi'  \oplus \Tv Q $. Here $\xi'$ denotes the restriction to $Q$ of the bundle $\mathrm{pr}_1^* (\tau^* \xi|Y_+)  \oplus \mathrm{pr}_2^*(\xi|Y_-)$
over $Y_+ \times_{\partial S_-} Y_-$.
\end{itemize}
\end{Defn}

\begin{Lem}\label{canonical pin}
Suppose that the first projection embeds $Q$  in $Y_+$, and that the second projection is an $S^1$--bundle over $Y_-$.
There is then a canonical relative pin structure for $(X,Q)$.
\end{Lem}
\begin{pf}
The boundaries $\partial S_-$ and $\partial S_+$ lie on different components of $S$, since the respective fibres of $\pi$ have different dimensions, $2n$
and $2(n-1)$ respectively.  Let $\xi$ be the bundle $\Tv X$ over the components containing $\partial S_-$, and the trivial bundle over those containing $\partial S_+$. Thus
$\xi' \to Y_+ \times_{\partial S_-} Y_-$ is $\mathrm{pr}_2^*\Tv Y_-$. We exhibit a pin structure on $\xi' | Q \oplus \Tv Q$. 

$  \Tv Q$ contains a distinguished line-subbundle $\lambda$, the kernel of the projection to $\Tv Y_-$. Choose a Euclidean metric on $ \Tv Q$, 
and notice that $\lambda^\perp = \mathrm{pr}_2^* \Tv Y_-$. The splitting $ \Tv Q=\lambda \oplus \mathrm{pr}_2^* \Tv Y_-$ reduces the structure group 
of $ \xi'|Q \oplus \Tv Q$ to the subgroup 
\[  \orth(1)\times \mathrm{SO}(2n-2) \hookrightarrow \orth(4n-3),\quad (A,B)\mapsto \mathrm{diag}(A,B,B). \]
We claim that this homomorphism factors through $\pin^+(4n-3)$.  This will give a canonical relative pin structure.

The restriction to identity components, $\{1\} \times \SO(2n-2) \to \SO(4n-3)$, lifts to $\Spin(4n-3)$ because it kills $\pi_1 \SO(2n-2)$. Write this lift as 
$l$. Let $\sigma \in \pin^+(4n-3)$ be an element of order 2 which maps to $(-1,I_{2n-2},I_{2n-2})$; this exists by Example \ref{odd pin}. The map $(A,B)\mapsto \sigma(A) l(A)$ gives the claimed lift. 
\end{pf}

\begin{Prop}\label{pin}
Suppose $J$ is a compatible complex structure for $(E,\Omega)$ such that $D\pi\circ J = \ii \circ D\pi$. Consider the moduli space of $J$-holomorphic sections $\EuScript{M}_J(X,Q)$; assume it is transversely cut out. A relative pin structure for $(X,Q)$ induces an orientation for any component of $\EuScript{M}_J(X,Q)$ on which the Maslov index is odd. A relative pin structure together with an orientation for one fibre $Q_x$ induces an orientation for components of even Maslov index.
\end{Prop}
\begin{pf}
We assume, up until the last moment, that the base is a disc $D$.
 
Step (i). The regularity assumption means that the cokernels of the relevant linear Cauchy--Riemann operators are trivial, and hence that the moduli space is smooth. The tangent space to $\EuScript{M}_J(X,Q)$ at $u$ is then the index of a Fredholm operator associated with the Cauchy--Riemann problem $(u^*\Tv X, \partial u^* \Tv Q)$. Thus the top exterior power of the tangent space is the determinant line of the corresponding Fredholm operator. The usefulness of working over a disc is that $u^*\Tv X$ can be trivialised in a canonical way, up to homotopy. Hence the top exterior power of the tangent bundle to $\EuScript{M}_J(X,Q)$ is isomorphic to the pullback of $\Det(\EuScript{F})$ (the determinant bundle of the universal family of Cauchy--Riemann operators over $\EuScript{L}_k\Gr_n$) by a classifying map 
\[\tau \colon \EuScript{M}_J(X,Q) \to \EuScript{L}_k\Gr_n.\] 

Step (ii). On the other hand, $\Det(\EuScript{F})$ gives rise to a double covering $p\colon \tilde{\EuScript{L}}\to \EuScript{L}_k\Gr_n$ as the $S^0$--bundle sitting inside the line bundle. The pullback line bundle $p^*\Det(\EuScript{F})$ has a \emph{tautological} trivialisation. It will therefore suffice to show that a stable pin structure determines a lift of $\tau$ to a map $\tilde{\tau}\colon  \EuScript{M}_J(X,Q) \to \tilde{\EuScript{L}}$, so that the pullback by $\tilde{\tau}$ of the tautological trivialisation of $p^*\Det(\EuScript{F})\to \tilde{\EuScript{L}}$ gives the sought orientation. 

Step (iii). Now consider the pullback of the tautological $\R^n$--bundle $U\to \Gr_n$ (whose fibre at $\Lambda$ is $\Lambda$ itself) by the evaluation map $e \colon S^1\times \EuScript{L}\Gr_n \to \Gr_n$. There is a double covering of $\EuScript{L}\Gr_n$ whose fibre at $F$ is a stable pin structure in the bundle $(e^*U)|S^1\times \{F\}$ over $S^1$. This double covering is the $S^0$--bundle of a unique line bundle $\upsilon$. By construction, $w_1(\upsilon)=T(w_2)$.

We are given a stable pin structure on $\xi|Q\oplus \Tv Q$. This pulls back to a stable pin structure on $(\partial u)^*(\xi|Q \oplus \Tv Q) $, for any $u\in \EuScript{M}_S(X,Q)$. But $(\partial u)^* (\xi|Q)$ is canonically trivial, since it extends to an oriented vector bundle $u^*\xi$ over the disc. Hence we obtain a stable pin structure on $(\partial u)^* \Tv Q$. In other words,  our relative pin structure gives  rise to a section of the covering $\tau^*\upsilon$.

Step (iv). Our task is now to compare the information found in steps (ii) and (iii).
Suppose that $k$ is odd. According to Equation (\ref{w1}), one then has $w_1 \Det(\EuScript{F})=w_1(\upsilon)$; thus, $\Det(\EuScript{F})$ and $\upsilon$
are isomorphic line bundles. To specify an isomorphism, one has only to do so for some fibre. The particular choice will not really matter for our purposes; notionally, however, we follow the convention of \cite[Lemma 11.17]{Se6}. \emph{A fortiori}, one then has an isomorphism $\tau^*\Det(\EuScript{F})\cong\tau^*\upsilon$. As in (iii), the pin structure gives rise to a section $\sigma$ of the $S^0$--bundle in $\tau^*\upsilon$; composed with the isomorphism with $\tau^*\Det(\EuScript{F})$, one then gets a lift $\tilde{\tau}$ as demanded in step (ii).

When $k$ is even, one proceeds along similar lines, but because of Equation (\ref{w1}), the result is only an isomorphism $\Det(\EuScript{F}) \cong \upsilon \otimes \lambda$, where $\lambda$ is the pullback to $\EuScript{L}_n\Gr_n$ of the line bundle $\Lambda^n U \to \Gr_n$. The pin structure then trivialises $\tau^*\upsilon$, but  one also needs the given orientation of $Q_x$ to trivialise $\tau^* \lambda$.

Finally, in the general case where $S$ is not a disc, one proceeds by a trick of `pinching off the boundary', as explained in both \cite{FOOO,Se6}.
\end{pf}

\begin{Rk}
Suppose $k$ is even, that $(X,Q)$ admits a relative pin structure, and
that  the fibres $Q_z$ of $Q\to S^1$ are orientable but $Q$ itself is not.  Choosing a relative pin structure and an orientation for $Q_z$, we determine an orientation for $\EuScript{M}_S(X,Q)$. On the other hand, sliding $z$ around the circle, we eventually come back to the \emph{opposite} orientation of $Q_z$, and hence to the opposite orientation of the moduli space, which is absurd. Hence the moduli space is empty in this case---an observation which one can also prove in more elementary fashion.
\end{Rk}

\subsubsection{Invariants}
We can now define our prototypical Lagrangian matching invariants, using a familiar procedure. Suppose that $(E,\pi,\Omega,J_0,j_0)$ is a symplectic Morse--Bott fibration, $Q$ a Lagrangian matching condition, $h\in \pi_0\sect(E,Q)$. Let $I\subset \mathrm{int}(S)$ and $J\subset \partial S$ be finite sets. Suppose \emph{either} that the fibres, critical manifolds and Lagrangians are monotone, \emph{or} that they satisfy the negativity assumptions of the last proposition. 

We write $(E^\#,Q^\#)$ for $(E,Q)$ with the enhancements of a relative pin structure for $(E,Q)$ and, if $Q$ is orientable, an orientation.

The invariant takes the form of homomorphism
\begin{equation}\label{def Phi}
 \Phi_{E,Q;  I,J} (h) \colon \bigotimes_{s\in I}{H_*(E_{s};\Z) }\otimes \bigotimes_{s'\in J} {H_*(Q_{s'};\Z)} \to \Z.  
\end{equation}
Thinking of the products of cycles $Z$ and $Z'$ as defining monomials $[Z]$, $[Z']$ in the tensor products, we put
\[ \Phi_{E,Q; I, J } (h) ( [Z] \otimes [Z'])= 
\begin{cases}
 \#  \EuScript{M}_{J,j}(E,Q) \times_{(\ev ,\zeta,\times\zeta')}(Z\times Z' ), & D(h; Z, Z')=0;  \\
 0 & D(h; Z, Z')\neq 0.
\end{cases}\] 
Here $\#$ is the signed count of points in a compact, oriented zero-manifold. This is independent of choices  because the one-parametric moduli spaces are compact, oriented one-manifolds with boundary.

\subsubsection{Floer homology}

Here we review the definition of Floer homology for symplectic automorphisms. This is a well-established theory, and we shall be brief (Fredholm and gluing theory are as in \cite{Sch}; transversality and compactness as in \cite{HoS, Se2}).

\emph{Preliminaries:}
\begin{itemize} 
\item
Recall that we have defined a locally Hamiltonian fibration (LHF) to be a triple $(E,\pi,\Omega)$, where $E$ is a smooth manifold with boundary, $\pi$ a smooth proper submersion $\pi\colon E\to B$ to another manifold with boundary, mapping $\partial E$ to $\partial B$, and $\Omega\in \Omega^2(E)$ a closed two-form which is non-degenerate on the fibres. 
\item
For a commutative ring $R$, the universal Novikov ring $\Lambda_R $ is the ring of formal series $\sum_{c\in \R}{a(c) t^c}$, where $a\colon \R\to R$ has the property that $\supp(a)\cap (-\infty,C)$ is a finite set for each $C\in \R$. When $R$ is a field, so is $\Lambda_R$.
\end{itemize}
Floer homology theory assigns to each  LHF $(Y,\pi,\sigma)$ over a closed, oriented 1--manifold $Z$, a $\Z/2$-graded $\Lambda_R$-module
\[   HF_*(Y,\sigma)=HF_0(Y,\sigma)\oplus HF_1(Y,\sigma), \]
though we should give two caveats: (i) the fibres $Y_z$ of $\pi$ should be weakly monotone symplectic manifolds in the sense of \cite{HoS} (the definition was given in the discussion preceding Definition \ref{kneg}); (ii) the definition is simpler when $R$ has characteristic 2.

There is a splitting (of $\Z/2$--graded modules) into `topological sectors', i.e. components of the space of sections, $\sect(Y)$:
\[ HF_*(Y,\sigma) =  \bigoplus_{\gamma\in \pi_0 \sect (Y)}HF_*(Y,\sigma)_\gamma. \]
There is a continuation isomorphism $HF_*(Y,\sigma) \to  HF_*(Y,\sigma') $ when there is an $\alpha\in \Omega^1(Y)$ with 
\[ \sigma-\sigma' = d\alpha, \quad  \alpha(\Tv Y)=0.\]

$HF_*(Y,\sigma)$ is, formally, the Morse--Novikov homology of the action 1-form $\mathcal{A}_{Y,\sigma}\in \Omega^1(\sect(Y)) $:
\[\mathcal{A}_{Y,\sigma}(\gamma; \xi) = \int_{Z}{\sigma(\dot{\gamma}, \xi )},  \quad \xi \in C^\infty_Z(\gamma^* \Tv Y). \]
The set of zeros of $\mathcal{A}_{Y,\sigma}$ coincides with the set $\hor(Y,\sigma)$ of horizontal sections defined by the Hamiltonian connection:
\[ \hor(Y,\sigma)=\{  \nu\colon Z\to Y : \, \pi\circ\nu=\id_Z, \, \im D_z \nu \subset \Th X_{\nu(z)}   \}.  \]
One should again enlarge the space of sections $\sect(Y)$ to a Banach manifold, namely $\sect^2_1(Y)$, which has tangent spaces $L^2_1(\gamma^* \Tv Y)$. The transversality condition for a zero $\nu$ of $\mathcal{A}_{Y,\sigma}$ is surjectivity of the operator
\[\xi\mapsto \nabla_{\partial_t}\xi; \quad L_1^2(\nu^*\Tv Y)\to L^2(\nu^* \Tv Y),\]
where $\nabla$ is the intrinsic connection along $\nu$. After fixing basepoints $z_i\in Z$, one in each component, we can consider the linear holonomy maps $L_{\nu,i} \in \End(\Tv Y_{\nu(z_i)} )$. Transversality of $\nu$ is equivalent in turn to the invertibility of the linear maps $\id - L_{\nu,i}\in \End(\Tv_{\nu(z_i)} Y)$. This is the condition that $\nu(z)$ is a non-degenerate fixed point of the monodromy $\phi \in \aut(Y_z,\sigma|_{Y_z})$. It follows that when $Y$ is \emph{non-degenerate}, meaning that every $\nu\in\hor(Y,\sigma)$ is transverse, the set $\hor(Y,\sigma)$ is finite. 
In the next section, we shall see how to assign to each $\nu$ a free $R$--module of rank 1, denoted by $|o(\nu)|$. When $R$ has characteristic 2, $|o(\nu)|=R$.
Now, for a non-degenerate LHF $(Y,\pi,\sigma)$, one sets
\[ CF(Y,\sigma) = \bigoplus_{\nu\in \hor(Y,\sigma)}{ |o(\nu)|\otimes_R \Lambda_R},\]
%Denote the generators of this free $\Lambda_R$-module by $\langle \nu \rangle$. 
This is the module underlying Floer's (co)chain complex. The definition does not yet use the orientation of $Z$. 
For a general LHF, one perturbs it to a non-degenerate one and proceeds as before.

Each $\nu \in \hor(Y,\sigma)$ has a Lefschetz number $l_\nu \in \Z/2$, defined when $Z$ is connected to be $0$ if $\det(\id-L_\nu)>0$ and $1$ if $\det(\id-L_\nu)<0$, and in general by summing the Lefschetz numbers over components of $Z$. These give the $\Z/2$-grading of the complex, $ CF = CF_0\oplus CF_1,$
with $CF_i$ generated by $\{\nu: l_\nu=i\}$. It follows that when $Z=\coprod{Z_i}$ there are canonical $\Z/2$-graded isomorphisms $CF(Y)\cong \bigotimes_i{CF(Y|{Z_i})}$. 

A compatible vertical almost complex structure $J^{v} \in \EuScript{J}(\Tv Y,\sigma)$ extends uniquely to an almost complex structure $J_0$ on $T Y \oplus \varepsilon^1$ (where $\varepsilon^1$ is the trivial real line bundle over $Y$) such that $J(\widetilde{\partial_t})=1\in \R$ (here $\partial_t$ is the unit length vector field on $S^1= \R/\Z$ and $\widetilde{\partial_t}$ its horizontal lift).

Since $T(Y\times \R)= TY\oplus \varepsilon^1$, $J_0$ induces a translation-invariant complex structure on $T(Y \times \R)$. It is still denoted $J_0$, and called {\bf cylindrical}. If we regard $Y\times \R$ as a locally Hamiltonian fibration over the cylinder $S^1\times \R$, with two-form $\mathrm{pr}_1^* \sigma$, then $J_0$ preserves both vertical and horizontal subbundles, and the projection to $S^1\times \R = \C/\ii \Z$ is holomorphic. 

The differential on $CF_*(Y,\sigma)$ is defined via the moduli space $ \EuScript{M}_{J,j}(Y,\pi ,\sigma) $ of finite-action pseudo-holomorphic sections of 
\[\pi\times\id\colon Y\times \R \to Z\times \R,\] i.e. sections $u$ satisfying $ J \circ (D u) = (D u) \circ j$ and $\int_{Z\times\R}{u^*\sigma}  <\infty$. Under the assumption of non-degeneracy, any $u\in \EuScript{M}_{J,j}(Y\times\R,\pi\times \id,\mathrm{pr}_1^*\Omega)$ has the following asymptotic behaviour: 
\begin{itemize}
\item
As $s\to \pm \infty$, the loops $ u (s,\cdot)$ converges pointwise towards a horizontal section $\nu^\pm \in \hor(Y,\sigma)$. 
\item
The convergence is exponentially fast with respect to the Riemannian metric $g$ given by $\sigma(\cdot, J \cdot)$ on $\Tv Y$ and for which $\widetilde{\partial_t}$ is a unit-length vector field.
\end{itemize}
When $Y$ has weakly monotone fibres, there is a good transversality theory for gradient trajectories (see \cite{FHS, HoS}). This shows that a dense set of cylindrical almost complex structures are `regular' in the sense that (i) the deformation operators are cylindrical for all trajectories $u\in \EuScript{M}_{J,j}^{\leq 2}(Y,\pi,\sigma)$ (the superscript $\leq 2$ refers to the index); (ii) these trajectories do not hit any point of $Y$ which lies on a pseudo-holomorphic sphere $S\subset Y_t$ with $c_1(S)\leq 0$. 
For regular almost complex structures, the moduli spaces $\EuScript{M}_{J,C}^{\leq 2}(\nu_- | Y | \nu_+ )$ of trajectories of index $\leq 1$, asymptotic limits $\nu^\pm$, and action $\leq C$, have compactness (as well as regularity) properties: any sequence has a subsequence converging in the Gromov--Floer topology to a broken trajectory.

We shall show in the next section how to assign to an arbitrary $u\in \EuScript{M}^1_{J,j}(\nu_-|Y|\nu_+)/\R$ an isomorphism $i_u\colon |o(\nu_-)| \to |o(\nu_+)|$. 
To define the differential in the Floer complex, we formally sum these isomorphisms, weighting them by their action. Thus, for $x\in |o(\nu_-)|$, we set
\[ \partial x  =\sum_{\nu_+\in \hor(Y,\sigma)} { \sum_{u\in \EuScript{M}_{J,j}^1(\nu_-|Y|\nu_+)/\R}{ i_u(x) \otimes t^{\mathcal{A}(u) }}},  \]
and extend $\partial$ to a $\Lambda_R$--linear endomorphism of $CF_*(Y,\sigma)$.
There are chain homotopy equivalences (`continuation maps') between the Floer complexes for $(Y,\sigma;J)$ and $(Y,\sigma',J')$ when there exists $\alpha\in \Omega^1(Y)$ with $ \sigma-\sigma' = d\alpha$ and $\alpha(\Tv Y)=0$. The argument by which one proves that $\partial ^2 =0$ and that continuation maps are chain maps is famous enough that we will not even adumbrate it. 
However, we do still need to say what $|o(\nu)|$ and $i_u$ are, and we shall do that in the next section.

\subsection{Coherent orientations in fixed point Floer homology} 
Coherent orientations for Floer-theoretic moduli spaces were introduced by Floer and Hofer \cite{FH}. The essential technical ingredient is their (linear)  gluing construction for cylindrical Cauchy--Riemann operators, together with an index theorem which expresses the determinant of the glued operator as the tensor product of those of the factors. Besides \cite{FH}, we draw on Seidel's  treatment  of the Lagrangian case \cite{Se6}. 

In our sketch of Floer homology, we have so far avoided discussing the relevant linearised operators. These are \emph{cylindrical Cauchy--Riemann operators}. Let $S$ be a Riemann surface with cylindrical ends (and compact outside the ends). We compactify it to a Riemann surface $\overline{S}$ with one boundary component for each end: this is done by identifying an outgoing end $S^1\times (0,\infty)\subset S$ with $S^1\times (0,1)$ via the diffeomorphism $(0,\infty) \to (0,1)  $, $s\mapsto s(1+s^2)^{-1/2}$, then adding a circle $S^1\times \{1\}$; similarly for incoming ends. Take a hermitian vector bundle $(E,\| \cdot \|)$ over $\overline{S}$, and a unitary connection $\nabla$ in it. The Cauchy--Riemann operator $\dbar_\nabla$ is the  operator $C^\infty(S,E) \to \Omega^1(S,E)$ defined by
\[  \dbar_\nabla \xi  =  (\nabla \xi)^{0,1}.  \] 
The behaviour of $\nabla$ over the boundary of $\overline{S}$ is of importance. Namely, for each component $T$ of the boundary, one wants the linear map $\nabla_{\partial/\partial t}$, operating on sections of $E| T$, to have trivial kernel. In this case, $\dbar_\nabla$ is called a \emph{non-degenerate} Cauchy--Riemann operator, and it extends to a Fredholm operator $L^p_1(S,E)\to L^p(S, T^*S\otimes E)$.

Take $\nu\in \hor(Y,\sigma)$, and form the pullback $\nu^*\Tv Y\to S^1$ as a symplectic vector bundle. It carries a canonical symplectic connection $\nabla_{\nu,\sigma}$, the resulting of linearising the connection in $\Tv Y$ defined by $\sigma$, and the non-degeneracy assumption is that the monodromy $L=L_\nu$ of $\nabla_{\nu,\sigma}$ does not have 1 as an eigenvalue. 

Let $x=\nu([0])\in M$, and consider the space $\EuScript{P}_{L}$ of paths $\gamma\colon [0,1] \to \U(T_x M)$ with $\gamma(0)=1$ and $\gamma(0)=L$. The component-set $\pi_0(\EuScript{P}_{L}) $ is an affine copy of $\Z$; indeed, the homotopy classes of such paths $\gamma$ make up the fibre over $L$ of the universal covering $\widetilde{\U}(T_x M) \to \U(T_x M)$.  We denote by $L^\#[n]$ the result of acting on $L^\# \in \pi_0(\EuScript{P}_{m} ) $ by $n\in \Z$.

Associated with any $\gamma\in \EuScript{P}_{L} $ is a Cauchy--Riemann operator $\dbar_{\gamma}$ in the trivial vector bundle $T_x M\times \C \to \C$. The assumption that $L$ does not have 1 in its spectrum guarantees its non-degeneracy. It has a determinant line $\delta_\gamma = \det \dbar_{\gamma}$, and these lines make up a bundle $\delta \to \EuScript{P}$. 

Now, each component of $\EuScript{P}$ is simply connected, since its fundamental group is isomorphic to $\pi_2(\U(T_x M)) = 0$; hence $\delta$ is a trivial line bundle.  We can therefore associate with a lift  $L^\# \in \widetilde{\U}(T_x M)$ of $L$ a line $o(L^\#)$, unique up to canonical isomorphism, by putting $o(L^\#)= \delta_\gamma$ for any $\gamma$ representing the component $L^\#$.

The orientation group $|o (L^\#)| $ is the abelian group generated by the two orientations $\omega_1$ and $\omega_2$ of $o_\gamma$, modulo the relation that $\omega_1 + \omega_2 =0$. Thus a choice of isomorphism determines an isomorphism $|o (L^\#) |\cong \Z$. 

How does $o(L^\#[n])$ compare to $o(L^\#)$? The answer is pleasantly simple: there is a canonical homotopy-class of isomorphisms between them. 
This is proved using the Floer--Hofer gluing theorem: $o(L^\#[n])^*\otimes o(L^\#)$ is the determinant line of a Cauchy--Riemann operator over $S^2$, obtained by gluing. Such determinant lines are canonically oriented since they form a simply-connected space which contains determinant lines of complex-linear operators.  

The last observation shows that the orientation groups $|o (L^\#)[n]|$ are canonically isomorphic to a single group which we may denote by $|o(\nu)|$. Formally, $|o(\nu)|$ is the group of families $(a_n)_{n\in \Z}$  where $a_n \in |o(L^\#)[n]|$ and $a_n$ maps to $a_0$ under the canonical isomorphism $ |o(L^\#)[n]|\to |o(L^\#)|$.

We now have the promised orientation lines $|o(\nu)|$. The next step is to show how $u\in \EuScript{M}_J(\nu_-,\nu_+)$ determines an isomorphism $|o(\nu_-)|\to |o(\nu_+)|$. For this we invoke, once again, the linear gluing theorem of Floer--Hofer. The tangent space to $\EuScript{M}_J(\nu_-,\nu_+)$ is the kernel of a surjective Fredholm operator $D_u$---namely, a Cauchy--Riemann operator $\dbar_\nabla$ in the bundle $u^* \Tv Y$ over the cylinder, where $\nabla$ extends the canonical symplectic connection over the ends $\{ \pm \infty\}\times S^1$.  We cap off the two ends of the cylinder to make $S^2$. The bundle $\nu^* \Tv Y$, being trivial, extends to a hermitian vector bundle over $S^2$, and the connection $\nabla$ also extends. The resulting Cauchy--Riemann operator $L$ over $S^2$ is homotopic to a complex-linear operator, hence its real determinant line is canonically oriented. But it is also homotopic to an operator obtained by gluing $\dbar_\nabla$ to two cylindrical Cauchy--Riemann operators over discs, and this gives an isomorphism $\det L \cong \det(D_u)\otimes o(L_-^\#) \otimes o(L_+^\#)^* $, where $L_\pm$ is the monodromy of $\nabla_{\nu_\pm,\sigma}$, $L_-^\#$ is an arbitrary lift of $L_-$, and $L_+^\#$ the resulting lift of $L_+$ (obtained from it by path-lifting).  The upshot is that $u$ determines an isomorphism $\Lambda^{\mathrm{top}}T_u \EuScript{M}_{J}(\nu_-|Y|\nu_+) \cong o(L_-^\#)^* \otimes o(L_+^\#)$ (up to a positive factor).

Let us abbreviate $\EuScript{M}_{J}(\nu_-|Y|\nu_+) $ to $\EuScript{M}$, and
$\EuScript{M}_{J}(\nu_-|Y|\nu_+) /\R$ to $\EuScript{M}^*$. If $u$ is a non-constant trajectory, translation in the $\R$-direction determines a canonical tangent vector $\eta \in T_u \EuScript{M}$, and hence an isomorphism 
$\iota(\eta) \colon \Lambda^{\mathrm{top}}T_u \EuScript{M} \to  \Lambda^{\mathrm{top}}T_{[u]} \EuScript{M}^* $.
Thus there is a composite isomorphism 
\[ \Lambda^{\mathrm{top}}T_{[u]} \EuScript{M}^*\cong o(L_-^\#)^* \otimes o(L_+^\#). \] 
Now when $T_{[u]} \EuScript{M}^*$ is $0$-dimensional, its top exterior power is simply $\R$, which has a canonical generator $1$. Hence, plugging this generator into our isomorphism, we obtain a generator for $o(L_-^\#)^* \otimes o(L_+^\#)$. This, finally, induces an isomorphism
\[ i_u \colon |o(\nu_-)|\to |o(\nu_+)|. \]
It is not difficult to see that these isomorphisms are compatible with gluing, in the sense that when $u_1 \# u_2$ is a glued trajectory, $i_{u_1 \# u_2}$ is the composite of $i_{u_1}$ and $i_{u_2}$ (cf. the discussion in \cite[Section 12b]{Se6}). This, of course, is what is needed for the proof that $\partial^2=0$.

\subsection{A field theory}
As observed by Seidel in \cite{Se4}, the Floer homology theory just described has an extension to an open-closed topological field theory in $(1+1)$ dimensions, coupled to singular fibrations. This unifies it with the Lagrangian matching theory. An `object' in the field theory is an LHF over a disjoint union of oriented circles. `Cobordisms' between objects 
are symplectic Morse--Bott fibrations over punctured matched surfaces, equipped with Lagrangian boundary conditions. To be precise, a cobordism from $\bigcup_{i=1}^p{(Y_i,\pi_i,\sigma_i)}$ to $\bigcup_{j=1}^q(Y'_j,\pi'_j,\sigma'_j)$ is defined by the following data:
\begin{itemize}
\item
A compact surface $S$ with matched boundary, together with a finite set $F=\{ z_1,\dots, z_{p+q}\}$ of interior marked points, of which $p$ are labelled as `incoming' and $q$ as `outgoing'.
\item 
Small disjoint coordinate discs  $ \xi_i \colon (\Delta,0)\to (S,z_i) $ contained in $\mathrm{int}(S)$. 
\item
A symplectic Morse Bott fibration $(E,\pi,\Omega,J_0,j_0)$ over $S\setminus F$, equipped with a Lagrangian boundary condition $Q$.
\item
Diffeomorphisms $\phi_j$ fitting into commutative diagrams 
\[ \begin{CD} \xi_j^* E @>{\phi_j}>>  Y_j' \times [0,\infty) \\
@VVV @VVV \\
\Delta^*  @>{re^{\ii\theta}\mapsto (e^{\ii\theta}, -\log r)}>> S^1\times [0,\infty)
\end{CD} \]
at each of the outgoing punctures, where $\Phi_j^*\Omega = \mathrm{pr}_1^*\sigma_j'$. At the incoming punctures, one has rather 
\[ \begin{CD} \xi_i^* E @>{\psi_i}>>   Y_i \times (-\infty,0] \\
@VVV @VVV \\
\Delta^*  @>{re^{\ii\theta}\mapsto (e^{-\ii\theta}, \log r)}>> S^1\times (-\infty,0].
\end{CD} \]
\end{itemize}
There is then an obvious notion of composition of cobordisms.

\begin{Rk}In formulating the field theory, it is convenient to work over the Novikov ring $\Lambda_k$ of an arbitrary base field $k$, rather than over $\Lambda_\Z$. This is itself a field, so taking homology commutes with $\otimes$ and $\Hom$.\end{Rk} 
 
Assuming that $(E,Q)$ is \emph{either} a fibrewise monotone fibration with a fibrewise monotone boundary condition, \emph{or} a fibrewise $3$-negative fibration with a fibrewise $2$-negative boundary condition, there are now relative invariants
\begin{equation}\label{rel invariants} 
\Phi_{E,Q} \in \Hom_{\Lambda_{k}} \left(\bigotimes_{i=1}^p {HF_*(Y_i, \sigma_i)} ,  \bigotimes_{i=1}^q {HF_*(Y_j', \sigma_j')}\right), \end{equation}
subject to a gluing law (see \cite{Se2,Se4}; the gluing theory is done in \cite{Sch}).\footnote{Note that this is a \emph{topological} field theory; the invariants are defined by choosing a complex structure on the base, and are then independent of it. We do not consider the more subtle question of constructing a \emph{conformal} field theory.}

Another general property of Floer homology is the Poincar\'e duality isomorphism:
\begin{equation}\label{PD}
HF_*( -Y,\sigma) \cong HF_*(Y,\sigma)^*,\end{equation} 
where $(-Y,\sigma)$ refers to the fibration obtained by switching the orientation of the base circle. A change in the incoming/outgoing label of a point $\zeta_i$ dualises the corresponding group.

\subsection{Quantum cap product}
If one chooses further finite sets $I\subset S$ and  $J\subset \partial S_+$, disjoint from one another and from the set of punctures $F$, the relative invariant $\Phi_{E,Q}\in \Hom(HF_*(Y,\sigma),HF_*(Y',\sigma'))$ generalises  to a homomorphism
\begin{equation} \label{rel Phi}
\Phi_{E,Q;I,J} \colon \bigotimes_{s\in I}{ H_* (E_{s};k) } \otimes \bigotimes_{s'\in J} {H_*(Q_{s'};k) } 
 \to \Hom_*(HF_*(Y,\sigma), H F_*(Y',\sigma')).\end{equation}
The homomorphism is the map on homology associated with a chain-level map $C\Phi$. The latter is defined using moduli spaces of index $0$, finite action pseudo-holomorphic sections with boundary on $Q$, hitting chosen cycles in the marked fibres $E_{s_a}$ and $Q_{s'_b}$. 
When $S$ is a twice-punctured sphere, $\Phi_{E;\{s\},\emptyset}\colon H_*(X_s) \to \End_* HF_*(X)$ is, by definition, the quantum cap product on Floer homology. 

For simplicity, we state the next lemma in a version which concerns a single marked point (the reader will have no difficulty in generalising it to allow more points).
\begin{Lem}
Let $\Gamma \colon H_*(E_{s_0})\to H_*(E_{s_1})$ be the isomorphism obtained by parallel transport of cycles over a path $\{s_t\}_{t\in[0,1]}$ in $S$.
Then  
\[\Gamma\circ \Phi_{E,Q; \{ s_0\}, \emptyset}(c) =\Phi_{E,Q; \{ s_1\}, \emptyset}(\Gamma (c)) .\] 
Hence $\Phi_{E,Q;I,J}(c)$ depends only upon the image of $c$ in $H_*(E)$. 
When $s \in \partial S_+$, one has 
\begin{equation}\label{moving between components}
\Phi_{E,Q,\{s\},\emptyset}(c)=\Phi_{E,Q,\emptyset,\{s\}}(i_*\circ \mathrm{pr}_1^!(c) ),\end{equation} 
where $i\colon Q_s \hookrightarrow E_s\times E_{\tau(s)}$ is the inclusion. One also has  \[\Phi_{E,Q,\{\tau(s)\},\emptyset}(c)=\Phi_{E,Q,\emptyset,\{s\}}(i_*\circ \mathrm{pr}_2^!(c)).\]
\end{Lem}
\begin{pf}
The first assertion is proved by considering the ends of 1-dimensional moduli spaces associated of sections passing through the smooth singular chain obtained by parallel transporting a cycle from one fibre to another. 

For the assertions about cycles in boundary fibres, one should consider a cycle $b$ in the product $E_{s} \times E_{\tau(s)}$. The fibre product of $b$ with a moduli space of sections $M$---taken via the evaluation map of $M$ into $E_{s} \times E_{\tau(s)}$---is identified with the fibre product of $b \cap Q_s$, taken via the evaluation map of $M$ into $Q$. Applying this to cycles $b$ of the form $c \times E_{\tau(s)}$ or $E_{s}\times c$ gives the two formulae.
\end{pf}

\subsubsection{A reduction to homology} 
Take two symplectic manifolds $(M,\omega)$ and $(\bar{M},\bar{\omega})$ with a Lagrangian submanifold $L\subset (M,-\omega)\times(\bar{M},\bar{\omega})$. Let $Y= S^1\times M$, $\sigma$ the pullback of $\omega$. Let $S= \Delta \cup \Delta$ be the union of two discs (with their boundaries labelled as $+$ then $-$) and equip it with the obvious matching $\tau$. Take the union of trivial fibrations, $E= (M\times \Delta) \cup (\bar{M}\times \Delta)\to S$. Puncture $S$ at interior points $\zeta_+$ in the first disc and $\zeta_-$ in the second. One obtains, by restriction, an LHF $E^*\to S\setminus \{ \zeta_+,\zeta_-\}$ with Lagrangian matching condition $Q$.

On the other hand, one has canonical (Piunikhin--Salamon--Schwarz, or `PSS') isomorphisms\footnote{See \cite{PSS}, or alternatively, \cite[Chapter 12]{MS}. The formal argument has to be backed up by transversality and compactness arguments (unproblematic in the monotone or 3--negative cases), by Floer's gluing theorem for  pseudo-holomorphic trajectories, and by a Gromov--Witten style gluing theorem for holomorphic spheres \cite[Chapter 10]{MS}.} 
\begin{align*}
 & \Pi_M \colon H_*(M;\Lambda_k)\to HF_*(M\times S^1, \mathrm{pr}_1^*\omega), \\
 & \Pi_{\bar{M}} \colon H_*(\bar{M};\Lambda_k)\to HF_*(\bar{M}\times S^1, \mathrm{pr}_1^*\omega).
 \end{align*} 
Note that these are valid with our orientation convention (the reason is almost the same as that which permitted the construction of the maps $i_u$ in our construction of Floer homology: the PSS map and its inverse are defined by means of holomorphic planes asymptotic to horizontal sections and subject to incidence conditions;  the moduli space of such planes, when it is $0$-dimensional, is `canonically oriented relative to the cylindrical end').
\begin{Prop}
Suppose that either (i) $M$ and $\bar{M}$ are monotone and $L$ a monotone Lagrangian, or (ii) $M$ and $\bar{M}$ are 3-negative and $L$ 2-negative. Then
one has
\[  (\Pi_{\bar{M}})^{-1}\circ \Phi_{E^*,Q}\circ \Pi_M (a) =   \pm \mathrm{pr}_{2*}([L] \cap \mathrm{pr}_1^! \cdot ) + \mathrm{h.o.t.} \] 
The letters h.o.t., abbreviating `higher order terms', designate a map which carries the subgroup $H_*(M;\Lambda^{\geq 0})\subset H_*(M;\Lambda_k)$ into $H_*(\bar{M};\Lambda^{>0})$. Here $\Lambda^{>0}\subset \Lambda^{\geq 0 }\subset \Lambda_k$ are  the subgroups consisting of series $\sum{a(r)t^r}
$ with $r>0 $ or $\geq 0$.
\end{Prop}
In words: the leading order part of the map on Floer homology (with respect to the filtration by the action functional, recorded by the Novikov coefficients)  is the classical map between homology groups associated with the correspondence, up to a sign which we do not check.
\begin{pf}
The relative invariant $\Phi_{E,Q; \{\zeta_+,\zeta_-\}, \{s\}}$ for the `completed' fibration $E$ may be formulated as a map 
\[ \Phi_{E,Q; \{\zeta_+,\zeta_-\}, \{s\}} \colon H_*(\bar{M}; \Lambda_k)\otimes H^*(L;\Lambda_k) \to H^*(\bar{M}; \Lambda_k). \]
We can decompose this map as $\Phi_{E,Q;\{\zeta_+,\zeta_-\}, \{s\}} = \sum_{\lambda\in \R}{\Phi_\lambda t^\lambda}$, where $\Phi_\lambda$ takes values in $H_*(\bar{M};k)$.
We have $\Phi_\lambda = 0$ for $\lambda >0 $. In fact, in the negative (rather than monotone) case, $\Phi_\lambda=0$ for $\lambda \neq 0$ for dimension reasons.

The vertical isomorphisms in the diagram are obtained from the relative invariants of trivial fibrations over the once-punctured two-sphere, and commutativity of the diagram follows from the gluing property for the relative invariants. It therefore suffices to prove that the leading term $\Phi_0$ of $\Phi(E,Q; \{\zeta_+,\zeta_-\}, \{s\})$ is $a\otimes c\mapsto  \mathrm{pr}_{2*}(c\cap \mathrm{pr}_1^! a)$.

Trivialising the fibrations, $\sect(E,Q)$ is identified with the space of maps
\[ (\Delta,\partial \Delta) \to (M\times \bar{M}, L) .   \]
We use a complex structure of the form $-J \oplus \bar{J}$. The relative invariant is then defined using two smooth singular cycles, $Z_\pm$, lying over $\zeta_\pm$, and and a third cycle $Z'$ in $L$. We then consider discs $\delta (\Delta,\partial \Delta)\to (M\times \bar{M}, L)$ such that $\delta(0) $ lies in the image of $Z_+\times Z_-$, and $\delta(1)$ lies in the image of $L$. The leading term $\Phi_0$ is computed using discs $\delta$ of index zero. Because of our monotonicity/negativity assumptions, index zero implies area zero, hence these discs are constant. Moreover, regularity for constant discs is the same thing as transversality for the cycles $Z_+\times Z_-$ and $Z'$. This reduces us to classical intersection theory, and so gives the result. \end{pf}

\begin{Rk}
The proposition says that the relative invariant gives a homomorphism $H^*(L)\otimes QH^*(M)\to QH^*(\bar{M})$, where $QH^*$ is quantum cohomology, deforming the homomorphism on classical cohomology groups induced by the correspondence $L$. These homomorphisms are little explored, even when $\bar{M}=\{\mathrm{pt.}\}$. 
\end{Rk}

\subsubsection{A speculative interlude}
Lagrangian matching invariants should have computational consequences in Floer homology. The following conjecture (which is the subject of work in progress) explains how this might work. 

Suppose that $(M,\omega)$ and $(\bar{M},\bar{\omega})$ are symplectic, $\rho\colon V\to\bar{M}$ an $S^k$--bundle with structure group $\SO(k+1)$, and $i\colon V\hookrightarrow M$ 
an embedding such that $i^*\omega=\rho^*\bar{\omega}$. Suppose that $\phi\in \aut(M,\omega)$ satisfies $\phi(V)=V$ and that $\phi|V$ covers $\bar{\phi} \in \aut(\bar{M},\bar{\omega})$. 
The cohomology of $\bar{M}$ acts on $HF_*(\bar{\phi})$ by quantum cap product; in particular, the Euler class of $V$ defines an endomorphism $e= e(V)\cap \cdot \colon HF_*(\bar{\phi})\to HF_*(\bar{\phi})[-k-1]$.
We form its mapping cone---a chain complex $\cone(e)$. 

There is a natural cobordism, in the sense of our field theory, from $\torus(\phi)$ to $\torus(\bar{\phi})$, and hence (assuming $\widehat{V}$ monotone in $-M\times \bar{M}$, say)
a Lagrangian matching invariant 
\[ \Phi \colon  HF_*(\bar{\phi}) \to HF_*(\phi). \]
It turns out that, when $\widehat{V}$ has high enough minimal Maslov index (at least $k+2$), the composite $ \Phi\circ  e$ is  nullhomotopic, in an essentially canonical way. Using the nullhomotopy, one builds a map
$ a  \colon H_*\cone(e)  \to HF_*(\phi)$.
On the other hand, there is the fibred Dehn twist $\tau_{V}\in \aut(M,\omega)$ studied in I.2. By considering a symplectic Morse--Bott fibration over an annulus, with one critical value, one gets a map $b\colon HF_*(\phi)\to HF_*( \phi\circ \tau_V)$. 
 \begin{Conj} \label{LES}
Under conditions which render $a$ and $b$ well-defined, they fit into a long exact sequence
\begin{equation}\label{qiso}
\xymatrix{
 H_*\cone(e) \ar^{a}[r]  & HF_*(\phi) \ar^{b}[r]  & HF_*(\phi\circ \tau_V) \ar@/^2pc/^{[-1]}  [ll]
}\end{equation} 
\end{Conj}
In the case of our Lagrangian matching conditions for relative Hilbert schemes, this can be translated into Seiberg--Witten theory, where it appears to be a somewhat unorthodox 
variant of the surgery triangle in monopole Floer homology \cite{KMOS} in which all the three--manifolds involved are fibred. Exactness of this triangle can in fact be deduced from that of the usual one together with a connected sum formula.

\section{Lagrangian matching invariants for broken fibrations}\label{invariants}
The previous section was devoid of examples. This section is devoted to a single class of examples: the relative Hilbert schemes associated with broken fibrations, and the Lagrangian matching conditions of Theorem \ref{Theorem B}.
\subsection{Monotonicity for symmetric products and their vanishing cycles}\label{mono}
\subsubsection{The non-separating case}\label{nonsep}
We consider $\sym^n(\Sigma)$, where $\Sigma$ is connected of genus $g$, and $n\geq 2$.  
\begin{enumerate}
\item
The Hurewicz map $\pi_2(\sym^n(\Sigma))\to H_2(\sym^n(\Sigma))$ has rank one.  Here is a different proof to those I have seen in the literature  \cite{BT, DS}: we deduce the result from Hopf's classic theorem that, in general, $\coker(\pi_2(X)\to H_2(X)) \cong H_2(\pi_1(X)).$\footnote{An isomorphism is induced by the map $X\to B\pi_1(X)$ classifying the universal cover. When $X=\sym^n(\Sigma)$ that map can be taken to be the Abel--Jacobi map to the Jacobian torus.} Since $\pi_1(\sym^n(\Sigma))=H^1(\Sigma)$, and the second integral homology of a free abelian group is its exterior square, we have $H_2(\pi_1) \cong \Lambda^2 H^1(\Sigma)$. But $H_2(\sym^n(\Sigma))\cong\Z\oplus \Lambda^2 H^1(\Sigma)$. This tells us $\im(\pi_2\to H_2)$ must be isomorphic to $\Z$, and that its generator must be a primitive class (i.e., not a non-unit multiple of another integer class).
\item
One can easily identify a spherical homology class $h\in H_2(\sym^n(\Sigma))$:   take any pencil of linearly equivalent effective divisors. One has $\langle \eta, h\rangle=1$ (which shows that $h$ is primitive) and $\langle \theta, h \rangle =0$. (The notation was explained in the introduction to this paper.) Since $c_1(T\sym^n(\Sigma))=(n+1-g)\eta-\theta$, the Chern number of $h$ is $n+1-g$. Thus $\sym^n(\Sigma)$ is monotone if $n\geq g$, regardless of the choice of symplectic form. 
\item
Let $\bar{\Sigma}$ be the result of surgery on a non-separating loop $L\subset \Sigma$.  Consider $\sym^n(\Sigma)\times\sym^{n-1}(\bar{\Sigma})$ with a symplectic form $-\omega\oplus\bar{\omega}$  that it is a sum of K\"ahler forms, with a sign reversal on the first factor, in cohomology classes $[\omega] = \eta_{\Sigma}+ \lambda \theta_{\Sigma}$ and $[\bar{\omega}] = \eta_{\bar{\Sigma}}+ \lambda \theta_{\bar{\Sigma}}$. This symplectic manifold is also monotone. The two generators for $\im(\pi_2\to H_2)$ (chosen with appropriate signs) each have Chern number $n+1-g$ and area $ 1$.
\item
Now consider the Lagrangian vanishing cycle
\[ \widehat{V}_L \subset \big( \sym^n(\Sigma)\times\sym^{n-1}(\bar{\Sigma}) , -\omega\oplus\bar{\omega}\big)  \]
associated with the loop $L\subset \Sigma$ (Part I: Theorem A). It was shown in Part I that every disc in $\pi_2(\sym^n(\Sigma)\times\sym^{n-1}(\bar{\Sigma}),\widehat{V}_L)$ lifts to a sphere in the product (Lemma 3.21). Hence the positive area discs $D$ for which $|\mu_{\widehat{V}_L}(D)|$ is least all have
\[\mu_{\widehat{V}_L}(D)= 2(n+1-g)= 2((n-1)+1-(g-1)), \]
and area 1. Thus $\widehat{V}_L$ is monotone as a Lagrangian when $n\geq g$.
\item
In the same situation, one checks that $\widehat{V}_L$ is 2--negative when $n\leq (2g(\Sigma)-1)/4$. This condition does not imply 3--negativity of $\sym^n(\Sigma)\times \sym^{n-1}(\bar{\Sigma})$, but it does imply 3--negativity of the two factors.
\end{enumerate}

\subsubsection{The separating case}\label{sep}
The case where $\Sigma$ is connected but $\bar{\Sigma}$ disconnected, with components $\Sigma_1$ and $\Sigma_2$, is more delicate.  One has
\[\sym^{n-1}(\bar{\Sigma}) = \bigcup_{i=0}^{n-1}{\sym^i(\Sigma_1)\times \sym^{n-1-i}(\Sigma_2)}.\] 
For brevity, we write $M$ for $\sym^n(\Sigma)$, $N$ for $\sym^{n-1}(\bar{\Sigma})$ and $N_i $ for the component $\sym^i(\Sigma_1)\times \sym^{n-1-i}(\Sigma_2)$ of $N$.  We put a symplectic form on $M\times N$ which, as before, is of shape $\omega\oplus -\bar{\omega}$ for K\"ahler forms $\omega$ on $M$ and $-\bar{\omega}$ on $N$. Since we want to find Lagrangian correspondences, we should suppose that $\omega$ represents a class $\eta_{\Sigma}+ \lambda \theta_{\Sigma}$ for some $\lambda>0$, and that  $\bar{\omega}$ 
represents $e + \lambda \theta_{\bar{\Sigma}}$, where $e$ restricts to the component $N_i$ as 
\[  { \mu_i \eta_{\Sigma_1}+ (1-\mu_i)   \eta_{\Sigma_2}}  \]
for some $\mu_i \in (0,1)$. When this is true, Theorem A from Part I gives us a correspondence $\widehat{V}_{L,k}\subset M\times N_i$ for each $i\in\{0,\dots,n-1\}$.

The factors $\sym^i(\Sigma_1)$ and $\sym^{n-1-i}(\Sigma_2)$ are all 3--negative when 
\[n-1 \leq \min_{j=1,2} {(g(\Sigma_j)-1)/2};\] 
under this hypothesis, we will have adequate control so far as bubbling of holomorphic \emph{spheres} is concerned.

The awkward point is that it is no longer true in general that every disc with boundary in  $\widehat{V}_L$ lifts to a sphere: when $\Sigma_1$ and $\Sigma_2$ have positive genus, there is a cyclic group's worth of non-spherical discs (Part I: Lemma 3.21) for each component of $\widehat{V}_L$. We have to determine the Maslov indices for the generators for these cyclic groups. 
 \begin{Lem}\label{Maslov comp}
Assume that $g(\Sigma_1)>0$, $g(\Sigma_2)>0$ and $n>1$. Then the group $\pi_2(M\times N_k, \widehat{V}_{L, k} )$ contains primitive, non-spherical classes $u_k$ such that 
\[\mu_{ \widehat{V}_{L,k}}(u_k)= k+1- 2g(\Sigma_1).\]  
The class $u_k$ generates $\pi_2(M\times N_k, \widehat{V}_{L, k})$ modulo the image of $\pi_2(M\times N_k)$.
\end{Lem}

Before we can prove this, we need to obtain a slightly more precise picture of $\widehat{V}_{L,k}$. We will achieve this by considering the model considered in Part I, Section 3.2.1, in which $\Sigma$, $\Sigma_1$ and $\Sigma_2$ are all 2-spheres. By interpreting the $n$th symmetric product of $\PS^1$ as $\PS^n$, we obtained an explicit equation for $\widehat{V}_{L,k}$ (the explicit model was denoted by $\EuScript{L}_{n,k}$).

\begin{Lem}
One of the $S^1$-fibres of $\EuScript{L}_{n,k}\subset \sym^n(\PS^1)$ consists of $n$-tuples of the form $[e^{i\theta}, e^{i\theta} \eta,\dots,e^{i\theta} \eta^{k-1};  \infty,\dots, \infty]$, where $\eta=e^{2\pi i/k}$. Indeed, under the quotient map
$\EuScript{L}_{n,k}\to \sym^{k-1}(\PS^1)\times \sym^{n-k}(\PS^1)$, this circle maps to the point
\[( (k-1)  \infty , (n-k) 0  ).\]
\end{Lem}

\begin{pf}
Our notation is from Part I, Section 3.2.1. The $S^1$-action $a_{n,k}$ preserves this circle, so it suffices to take $\theta=0$. Our identification $\sym^n(\PS^1)=\PS^n$ sends an $n$-tuple $\mathbf{x}=[x_1,\dots,x_k; \infty,\dots,\infty]$, where all the $x_j$ lie in $\C\subset\PS^1$, to the point $(\sigma_k(\mathbf{x}): \dots :\sigma_1(\mathbf{x}) :1 : 0 : \dots : 0 )$ ($n-k$ zeros). Here $\sigma_j$ stands for the $j$th elementary symmetric polynomial, and our convention is such that $\sigma_1(\mathbf{x}) = - \sum{x_i}$ and $\sigma_k(\mathbf{x})= (-1)^k\prod{x_i}$. The $n$-tuple $\mathbf{x}=[1,\eta,\dots,\eta^{k-1}; \infty,\dots,\infty]$ is mapped to $(1: 0:\dots : 0 :1 : 0 : \dots : 0 :)$, i.e. to the point $(z_0:\dots:z_n)$ where $z_0=z_k=1$ and $z_j=0$ otherwise. To see this, note that $\sigma_k(1,\eta,\dots,\eta^{k-1})=1$ (clear) but that $\sigma_j(1,\eta,\dots,\eta^{k-1})=0$ for $j<k$, since by homogeneity of $\sigma_j$,
\[ \eta^j \sigma_j(1,\eta,\dots,\eta^{k-1}) = \sigma_j(\eta ,\eta^2,\dots,\eta^k)=\sigma_j(1,\eta,\dots,\eta^{k-1}).\] 
This point $(z_0:\dots:z_n)$ certainly satisfies the defining equation for $\EuScript{L}_{n,k}$, which is $\sum_{j=0}^{k-1}{|z_j|^2=\sum_{j=k}^n}{|z_j|^2}$. The reduction map (dividing out the $S^1$-action) carries it to $( (1:\dots:0 : 1), (1:0: \dots : 0)\in \PS^{k-1}\times\PS^{n-k}$, which corresponds to the divisor $( (k-1)(1:0), (n-k) (0:1)) \in \sym^{k-1}(\PS^1)\times \sym^{n-k}(\PS^{n-k})$.
\end{pf}

\begin{pf}[Proof of Lemma \ref{Maslov comp}]
By exchanging the roles of $\Sigma_1$ and $\Sigma_2$, and of $k$ and $n-1-k$, we see that it is sufficient to assume that $k>1$.

Compactify each of the two components of $\Sigma \setminus \gamma $ to a Riemann surface-with-boundary $\tilde{\Sigma}_i$, $i=1,2$. Thus there is a map $\tilde{\Sigma}_i\to \Sigma_i$ collapsing the boundary to a point. Consider a holomorphic branched covering $b \colon \tilde{\Sigma}_1\to D$ over the closed unit disc, of degree $k$, restricting to an unramified cyclic covering $\partial \tilde{\Sigma} \to \partial D$ (also of degree $k$). (We are free to adjust the complex structure on $\Sigma$ to make sure that the covering exists.) Fixing a point $z\in \Sigma_2$, we obtain a holomorphic map
\[  b_z \colon D \to \sym^n(\Sigma),\quad x\mapsto b^{-1}(x) + (n-k)z. \]
Now, the coisotropic hypersurface $V_L\subset \sym^n(\Sigma)$ depends on the choice of a K\"ahler form on a relative Hilbert scheme, so it does not make sense to ask whether the boundary of $b_z$ lies on $V_L$. However, if we allow ourselves to move $V_L$ by totally real isotopies by the method of Part I, Section 3.6, adjusting the (non-closed) `good two-form' (Part I Section 2.3.2) defining the vanishing cycle, we may assume that this is indeed so. Furthermore we may arrange that the boundary curve of this holomorphic curve  circumnavigates the fibre of the projection $V_L \to N_k$ over a point $(P,Q)$. Indeed, by the last lemma, this behaviour is exactly what happens in the genus-zero model, which may be transplanted into higher-genus situations as in Part I, Section 3.6. (Because of the homotopical nature of Maslov indices, it is perfectly acceptable to sacrifice closedness of the non-degenerate two-forms.)

We can then promote $b_z$ to a holomorphic map 
\[ B_z \colon  D\to M\times N_k, \quad x\mapsto (b_x(z), P,Q).\]  
with boundary on $\widehat{V}_L$. We can see that $[B_z]$ is primitive in $\pi_2(M\times N_k, \widehat{V}_{L,k})$ simply by observing that its boundary is primitive in $\pi_1(\widehat{V}_L)$ (this requires $g(\Sigma_1)>0$).

To calculate the Maslov index of $B_z$, we shall exhibit a trivialisation $(Y_1,\dots, Y_{2n-1})$ of $(B_z| \partial D)^* T\widehat{V}_{L,i}$. Via a complex structure $J$, we then obtain a trivialisation 
$(Y_1,\dots, Y_{2n-1};  \\ J Y_1,\dots, JY_{2n-1}) $ for $(B_z| \partial D)^* T (M\times N_i)$, and hence an orientation $o$ for this last bundle. We then extend this orientation to a section of $\Lambda^{\mathrm{top} } B_z^* T (M\times N_k)$, and count its zeros with signs; the Maslov index is this count.

A warm-up exercise, left to the reader, is to derive the Chern number $n+1-g$ for the generating sphere in $\pi_2(\sym^n(\Sigma))$ by thinking of this sphere as a degree $n$ branched cover of $S^2$ and applying Riemann--Hurwitz. This exercise is useful in getting the local contributions to the Maslov index correct.

Start by trivialising $T(\partial D)$ via an anticlockwise-pointing vector field $Z$. 
Let $a_z=b_z|\partial D$ and $A_z=B_z|\partial D$.
Then $a_{z*}Z$ is a non-vanishing section of $a_z^* T V_L$. The required trivialisation for $A_z^* T \widehat{V}_L$ is obtained by putting $Y_i  =  X_i$ for $i=1,\dots, 2n-2$, where  $(X_1,\dots, X_{2n-2})$ is a positively oriented trivialisation of $T_{(P,Q)}(N_k)$, and $Y_{2n-1}=a_{z*}Z$. Now, $Z$ extends to a vector field $\tilde{Z}$ on $D$ with a single, non-degenerate zero $s$ of sign $+1$, which we may assume does not coincide with any of the critical values of the covering $b$. We may extend the fields $Y_i$ to fields $\tilde{Y}_i$  over $D$ by putting $\tilde{Y}_i= X_i$ for $i<2n-1$ and $\tilde{Y}_{2n-1}= b_{z*}\tilde{Z}$. 

The field $\tilde{Y}_{2n-1}$ has a $k$-fold zero at $s$, which contributes $2k$ to the Maslov index. 
The other contributions come from the critical points of $b$; each adds $-1$ to the Maslov index, and there are $k-1 + 2g(\Sigma_1)$ of them, by Riemann--Hurwitz. This gives us a Maslov index of
\[  \mu_{\widehat{V}_L}(B_z) =   2k - (k- 1+ 2g(\Sigma_1)) = k + 1 - 2g(\Sigma_1). \]
Putting $u_k=[B_z]$, we obtain the Maslov index claimed. 
\end{pf}

Now, we wish $\widehat{V}_{L,k}$ to be $2$--negative. We assume that $n \leq \frac{1}{2} \min(g(\Sigma_1),g(\Sigma_2)) $. It is not clear to the author whether $2$--negativity can be achieved simultaneously for all $k\in \{0 ,\dots, n-1\}$; rather, we will adjust the available parameters so that it is $2$--negative for any particular $k$. We first choose the coefficients $\mu_k$ so as to ensure that the product manifolds $M \times N_k$ are all negatively monotone, i.e., that the areas of the two generators for $\pi_2$ are proportional to their respective Chern numbers (with a negative constant of proportionality, $\lambda)$. We have to arrange that the area of the disc $u_k$ is in the ratio $2\lambda$ to its Maslov index. Now, $\widehat{V}_{L,k}$ depends on the Hamiltonian isotopy class of $L\subset \Sigma$. The areas of the disks that it bounds are sensitive to this Hamiltonian isotopy class: in fact, under a \emph{Lagrangian} isotopy of $L$, the area of the resulting $u_k$ will change proportionally to the flux of the isotopy.  Hence by moving $L$ we can ensure that the area of $u_k$ is $2\lambda$ times its Maslov index.

Once we have arranged this negative monotonicity, we can apply Lemma \ref{Maslov comp} to see that $2$--negativity holds provided that $n \leq \frac{2}{3} \min(g(\Sigma_1),g(\Sigma_2))$, which is certainly true when $n \leq \frac{1}{2} \min(g(\Sigma_1),g(\Sigma_2)) $.

Our conclusion is expressed by the following lemma.
\begin{Lem} \label{Maslov conclusion}
When $n \leq \frac{1}{2} \min(g(\Sigma_1),g(\Sigma_2)  $, all the symplectic manifolds $M$ and $N_k$ are 3--negative. For each $k\in \{0,\dots, n-1\}$, after an isotopy of $L$ inside $\Sigma$, the corresponding Lagrangian  $\widehat{V}_{L,k}\subset M\times N_k$ is 2--negative.
\end{Lem}

\subsection{Defining the invariants for broken fibrations on closed four-manifolds}\label{definition}
We are, at last, ready to give the definition of the Lagrangian matching invariants for broken fibrations on closed four-manifolds.

The invariants are indexed by $\spinc$-structures. For now we explain in schematic form the connection between moduli spaces of sections and $\spinc$-structures. The details will be given afterwards. 

Let $(X,\pi)$ be a broken fibration over a closed surface ($X$ connected) and let $Z\subset X^\crit$ be the 1--dimensional part of the set of critical points. 

(i) The `Taubes map'  $\tau_X$ is a bijection 
\[ \spinc(X) \to \delta^{-1}([Z])\subset H_2(X,Z),  \]
where $\delta\colon H_2(X,Z)\to H_1(Z)$ is the boundary homomorphism, and $Z$ is oriented by a vector field $v$ such that $i \pi_*(v)$ points into the side of $\pi(Z)\subset S$ on which the Euler characteristic is higher (here $i$ is the complex structure on $S$). The map $\tau_X$ arises from the canonical $\spinc$-structure $\mathfrak{s}_{\mathrm{can}}$ on the almost complex manifold $X\setminus Z$. It is characterised by the relation
\begin{equation}\label{taubes map}  
\tau_X(\mathfrak{s})=\beta\quad \text{if}\quad \mathfrak{s}|(X\setminus Z)=\mathrm{PD}(\beta)\cdot \mathfrak{s}_{\mathrm{can}}.   
\end{equation}
In the context of near-symplectic manifolds, the map $\tau_X$ plays an important role in Taubes' programme. In  \cite{Ta2}  it is proved that $\im(\tau)\subset \delta^{-1}([Z])$; it is then easy to see that $\tau_X$ is bijective.

(ii) Construct $X^{[\nu]}$ as in (\ref{Xnu}), where $\nu\colon S^{\mathrm{reg}}\to \Z_{\geq 0}$ satisfies $2\nu(s)+\chi(X_s)=2d$, with $d$ constant; and its Lagrangian matching condition $\EuScript{Q}$. We shall see presently how to  build a  map 
\[ \alpha_\nu \colon \pi_0 \sect(X^{[\nu]}, \EuScript{Q})\to \delta^{-1}([Z])\subset H_2(X,Z;\Z).\] 
If $\beta\in \im (\alpha_\nu) $ then one has $\langle \beta, X_s \rangle = \nu(s) $. Also, $ \langle \beta, [F] \rangle \geq 0$ for every homology class $[F]\in H_2(X\setminus Z;\Z)$ represented by a component of a regular or nodal fibre.

Define the map
\begin{equation}   \label{sigma}
\sigma_\nu = \tau_X^{-1}\circ \alpha_\nu \colon  \pi_0 \sect(X^{[\nu]}, \EuScript{Q})\to \spinc(X). \end{equation}
We will prove that $\sigma_\nu$ is injective, and, when $\inf(\nu)\geq 2$, describe its image.

Now take $\beta=\tau_X^{-1}(\mathfrak{s})\in H_2(X,Z;\Z)$. Let us recall a definition from I.1:
\begin{Defn}
Suppose that the fibres of $\pi$ are all connected. The $\spinc$--structure $\mathfrak{s}$ is then {\bf admissible} if (i) for every homology class $[F]\in H_2(X\setminus Z;\Z)$ represented by $F$, a regular fibre or an irreducible component of a nodal fibre, one has $ \langle c_1(\mathfrak{s} ), F \rangle \geq \chi(F)$ (equivalently, $ \langle \beta, [F] \rangle \geq 0$); and (ii) one of the following two conditions holds:  for each regular fibre $X_s$ one has \emph{either}
\begin{enumerate}
\item 
$  \langle  c_1(\mathfrak{s}) ,[X_s]\rangle  \leq \chi(X_s)/2 $
(equivalently, $ \langle \beta,[X_s] \rangle  \leq  -  \chi(X_s) / 4$); \emph{or}
\item
$ \langle c_1(\mathfrak{s}) ,  [X_s]  \rangle > 0 $ (equivalently, $ \langle \beta, X_s \rangle \geq   g(X_s)$). 
\end{enumerate} 
\end{Defn}
According to the results of Section \ref{nonsep}, condition (2) implies that the fibres $\sym^{\nu(s)}(X_s)$ are monotone, and that $\EuScript{Q}$ is a fibrewise-monotone matching condition. Condition (1)  implies that the fibres $\sym^{\nu(s)}(X_s)$ are 3--negative, and that the fibres of the Lagrangian matching condition are 2--negative. These are the conditions under which we are able to construct invariants.

\begin{Defn} In a broken fibration some of whose fibres have two connected components (but none more than two), we say that $\mathfrak{s}$ is admissible when 
\[  \chi(F) \leq  \langle  c_1(\mathfrak{s}) ,[F]\rangle  \leq \chi(F)/2 \]
for each fibre-\emph{component} $F$.
\end{Defn}

By the results of Section \ref{sep} (more particularly, Lemma \ref{Maslov conclusion}) admissibility guarantees 3--negativity of the relevant symmetric products, and---after suitable adjustments are made---2--negativity of the Lagrangian correspondences.

Let $I_\pi \subset  H_1(X;\Z)$ be the subgroup of classes supported on a fibre of $\pi$. Let $I^*_\pi = \Hom(I_\pi,\Z)$. Write 
\begin{equation}  
\mathbb{A}(X,\pi) =\Z[U] \otimes_{\Z} \Lambda^*I_\pi^*,   \end{equation} 
and make it a graded ring by declaring $U$ to have degree 2. An element $l$ of $\mathbb{A}(X,\pi)$ is a finite sum $\sum{U^n\otimes l_n}$. It can be thought of equally as a homomorphism $\Z[U]\otimes_{\Z} \Lambda^*I_\pi \to \Z$ (sending $U^m\otimes \lambda$ to $l_m(\lambda)$) and we write $l(m)$ to denote its value on  $m\in \Z[U]\otimes_{\Z} \Lambda^*I_\pi$. 

When $I_\pi=H_1(X;\Z)$, we write $\mathbb{A}(X)$ instead of $\mathbb{A}(X,\pi)$.  The reader can verify that this is always true when the base is $S^2$.

\begin{Rk}
At this point I should confess to a technical inconsistency 
with the announcements made in the introduction to Part I, 
arising because of a careless oversight. (This refers to the published version; the ArXiv version has been corrected.) In Part I, Section 1, $\mathbb{A}(X,\pi)$
was defined using the the dual group of $K_\pi = \ker (\pi_*\colon H_1(X)\to H_1(S))$ in place of $I_\pi^*$.  In general, $K_\pi$ will be larger than $I_\pi$, 
and I have not established how to work with the larger group. I apologise for 
this change. 
\end{Rk}

\begin{Defn}
Fix an orientation for each orientable attaching surface $Q_i$ for the broken fibration $(X,\pi)$. The {\bf Lagrangian matching invariant} 
\[\EuScript{L}_{(X,\pi)}(\mathfrak{s}) \in \mathbb{A}(X,\pi)\]
is a homogeneous element of degree $d(\mathfrak{s})$, where 
\begin{equation}\label{formal dimension}
d(\mathfrak{s}) = \frac{1}{4}\left (c_1(\mathfrak{s})^2  -2\chi(X)- 3\sigma(X)\right), \end{equation} 
characterised as follows. Let $m = U^a\otimes l_1 \wedge\dots \wedge l_b$ be a monomial of degree $d(\mathfrak{s})$.

Represent $l_i$ by a loop $\gamma_i \subset X_{s_i}$ on a regular fibre. The points $s_i$ should all be distinct. 

We associate with each $\gamma_i$ a smooth codimension-1 cycle in a symmetric product: 
\[ \delta_{\gamma_i}= \gamma_i + \sym^{\nu(s)-1}(X_s)  \subset \sym^{\nu(s)}(X_s).\] 
Choose $ a $ points $x_j \in X$, lying in distinct regular fibres of $\pi$, disjoint from the the fibres $X_{s_i}$. With each of these we associate the codimension-2 cycle 
\[ \delta_{x_j} = x_j +\sym^{\nu(s)-1}(X_{\pi(x_j)}) \subset \sym^{\nu(s)}(X_{\pi(s)}).\]
The symplectic Morse--Bott fibration $X^{[\nu]} \to S'$, and its Lagrangian matching condition $\EuScript{Q}$, give rise to a homomorphism 
\[\Phi =  \Phi_{ X^{[\nu]} ,
\EuScript{Q} ; \{ s_i \} \cup \{ \pi(x_j) \} , \emptyset}\colon \bigotimes_{i=1}^{a+b}{H^*(\sym^{\nu(s_i) }(X_{s_i}) ;\Z) }\to \Z, \]
as in (\ref{def Phi}). Set 
 \[ \EuScript{L}_{(X,\pi)}(\mathfrak{s}) \big(U^a\otimes [\gamma_1] \wedge\dots \wedge [\gamma_b]\big) 
   = \Phi([\delta_{x_1}] \otimes \dots \otimes [\delta_{x_a }]  \otimes  [\delta_{\gamma_1}] \otimes \dots \otimes [\delta_{\gamma_b}])   \] 
(the cohomology classes here have degrees 2 and 1, and are Poincar\'e dual to the cycles $\delta_{x_j}$, $\delta_{\gamma_i}$).

The definition involves an orientation for the moduli space of 
pseudo-holomorphic sections. This is defined via the relative pin structure for $\EuScript{Q}$ constructed by the general mechanism of Lemma 
\ref{canonical pin}. 

We also need to orient the components $\EuScript{Q}_i$ of $\EuScript{Q}$ for which an orientation exists. This is equivalent to orienting the attaching tori $Q_i \subset X$ (the tori which degenerate to critical circles): indeed, $\EuScript{Q}_i$ is an $S^1$-family of circle-bundles over symplectic manifolds. One orients these circle-bundles by orienting a fibre then orienting the base.  The base has its symplectic orientation, and one of the fibres may be taken to be a copy $\gamma \times (n-1)x$ of the attaching circle $\gamma \subset \Sigma$. Finally, one orients $\EuScript{Q}_i$ as 
\begin{align*}
 &(\text{canonical orientation of }S^1)\\
  \otimes & (\text{orientation of fibre of } \widehat{V}_\gamma \to \sym^{n-1}(\bar{\Sigma}))\\
   \otimes & (\text{symplectic orientation of }\sym^{n-1}(\bar{\Sigma})).
  \end{align*} 
\end{Defn}
That the degree of $\EuScript{L}_{(X,\pi)}(\mathfrak{s})$ is $d(\mathfrak{s})$ is a non-trivial matter: it is a direct consequence of our index formula, Theorem \ref{index formula}. Other than that, the well-definedness of $\EuScript{L}_{(X,\pi)} (\mathfrak{s})$ is  an easy consequence of the general theory we have set up. The following points are worth remarking:
\begin{itemize}
 \item
We must check independence from the choices of cycles $\delta_{x_j}$ and $\delta_{\gamma_i}$. One can certainly move the point $x_j$ inside $X'$: a path $x_j(t)$ gives rise to a homotopy of cycles $\sym^{\nu(x_j)-1}(X_{\pi(x_j(t))})\subset \sym^{n(x_j)}(X_{\pi(x_j(t))})$, and the parametrised moduli space is a compact one-manifold with boundary. It is also legitimate to move $x_j$ between different components of $X'$, because of (\ref{rel Phi}).

Similarly, one can replace $\gamma_j$ by a homotopic loop in a nearby fibre, or by a homologous point loop in the same fibre, and one can also pass between different components of $X'$.
\item 
We must show that $\Phi$ is symmetric in the $\delta_{x_j}$--factors and antisymmetric in the $\delta_{\gamma_i}$--factors. This holds because a pair of  cycles in distinct fibres can be replaced by the intersection of two cycles in a single fibre; thus symmetry and antisymmetry hold because the relevant cycles have, respectively, even and odd codimension.
\end{itemize}

\subsection{Sections of $\EuScript{Q}$ versus $\mathrm{Spin}^c$-structures}
We now explain the map $\alpha_\nu$ involved in the definition of $\sigma_\nu$ given in equation (\ref{sigma}), starting with a prototypical case.

We use the same notation as in Part I, Section 4. So, $\pi\colon Y\to S^1$ and $\bar{\pi}\colon \bar{Y}\to S^1$ are fibred three--manifolds, and $(X_{\mathrm{br}},\pi_{\mathrm{br}} )$ an elementary broken fibration over an annulus $A=\{ z\in \C : 1/2\leq |z|\leq2\} $ realising a cobordism between them: 
\[Y = \pi_{\mathrm{br}}^{-1}(\{|z|=1/2\}),\quad \bar{Y}= \pi_{\mathrm{br}}^{-1}(\{|z|=2\}).\] 
The fibre $\Sigma = \pi_{\mathrm{br}}^{-1}(1/2)$ is connected of genus $g$, whilst $\bar{\Sigma}:=\pi_{\mathrm{br}}^{-1}(2)$ is either connected of genus $g-1$, or else disconnected with components of genera $g_1$ and $g-g_1$.

We write $(Y^{[n]},\pi^{[n]})$ for the relative symmetric product $\sym^n_{S^1}(Y)\to S^1$, and  $(\bar{Y}^{[n-1]},\bar{\pi}^{[n-1]}$) for  $\sym^{n-1}_{S^1}(\bar{Y})\to S^1$. As usual, after choosing complex structures
on $\Tv Y$ and $\Tv \bar{Y}$ these become differentiable families of complex manifolds.

The first step is to look more carefully at the Taubes map $\tau_{X_{\mathrm{br}}}$. Let $X=X_{\mathrm{br}}$ and let $Z\subset X$ be the circle of critical points. Then we have a commutative diagram
\[\xymatrix{
 \spinc(X)\ar[r]^{\mathrm{restr.}}\ar[d]_{\tau_X} & \spinc(\bar{Y})\ar[d]^{\tau_{\bar{Y}}}_{\cong}\\ 
  H_2(X, \partial X\cup Z;\Z) \ar^{\delta}[r] &H_1(\bar{Y};\Z)\\
}\]
where
\begin{itemize}
\item The map $\delta$ is the obvious boundary homomorphism. A straightforward computation verifies that it is surjective. 
\item $\tau_X$ is the Taubes map, defined by the condition (\ref{taubes map}). In similar fashion, $\tau_{\bar{Y}}$ is defined by the relation  
\[  \mathfrak{t} =\mathrm{PD}(\tau_{\bar{Y}}(\mathfrak{t}))\cdot \mathfrak{t}_{\mathrm{can}},   \]
where $\mathfrak{t}_{\mathrm{can}}\in \spinc(\bar{Y})$ is the canonical $\spinc$-structure associated with the vertical two-plane field on $\bar{Y}$ (it is the restriction to $\bar{Y}$ of $\mathfrak{s}_{\mathrm{can}})$. 
\end{itemize}
For each integer $d$, let 
\begin{equation} \label{spinc d}
\spinc(\bar{Y})_d =\{ \mathfrak{t}\in \spinc(\bar{Y}):  \langle c_1(\mathfrak{t}) , [\bar{\Sigma}]\rangle = 2d  \},  \end{equation} 
and let $\spinc(X)_d$ be the preimage of $\spinc(\bar{Y})_{d}$ under the restriction map. 

What is the image of $\spinc(X)_d$ (resp. $\spinc(\bar{Y})_d$) under $\tau(X)$ (resp. $\tau_{\bar{Y}}$)? Write
\[ 2d =\chi(\Sigma)+ 2n = \chi(\bar{\Sigma})+ 2 (n-1).\] 
Then $\tau_X(\spinc(X)_d) $ is the affine subgroup 
\[ H_2^{(n)} = \{ \beta \in  H_2(X,Z\cup \partial X): \partial \beta = [Z] \; \text{mod} \; H_1(\partial X;\Z) \text{ and } \beta \cdot [\Sigma]=n \} . \] 
This follows immediately from the fact that $\mathfrak{s}_{\mathrm{can}}$ has degree $\chi(\bar{\Sigma})$ over $\bar{\Sigma}$. Similarly, $\tau_{\bar{Y}}(\spinc(\bar{Y})_d)$ is the affine subgroup 
\begin{equation}\label{H1 n}
H_1(\bar{Y};\Z)_n  = \{ \gamma \in H_1(\bar{Y};\Z): \gamma \cdot [\bar{\Sigma}]  =n \}.  \end{equation} 
If $\bar{\Sigma}$ is disconnected, with components $\bar{\Sigma}_1$ and $\bar{\Sigma}_2$, it become necessary to refine the set $\spinc(X)_d$. Write
\[ \spinc(\bar{Y})_d=\bigcup_{d_1+d_2 = d} {\spinc(\bar{Y})_{d_1,d_2} },\]
where $\mathfrak{t}\in \spinc(\bar{Y})_{d_1,d_2 }$ if $\langle c_1(\mathfrak{t}) , [\bar{\Sigma}_i]\rangle = 2d_i$. 
Let  $\spinc(X)_{d_1,d_2}$ be the preimage of $\spinc(\bar{Y})_{d_1,d_2}$ under restriction. Then $\tau_X(\spinc(X)_{d_1,d_2}) = H_2^{(n_1,n_2)}$, where $n_i = d_i-\chi(\Sigma)/2$, and $H_2^{(n_1,n_2)}$ is the obvious refinement of $H_2^{(n)}$.

Now we bring in the space of sections of the Lagrangian matching condition $\EuScript{Q}$.
\begin{Prop}\label{def alpha} 
Assume $\bar{\Sigma}$ is connected. Then there is a unique map 
\[\alpha\colon\pi_0\sect(\EuScript{Q})\to H_2^{(n)}\] 
which makes the diagram
\[\xymatrix{
   \pi_0 \sect(\EuScript{Q}) \ar@{-->}[r]^{\alpha} \ar[d]   & H_2^{(n)}\ar@{_{(}->}[d]^{\partial}\\ 
    \pi_0 \sect(Y^{[n]}) \times \pi_0 \sect(\bar{Y}^{[n-1]})\ar[r]   	&  H_1(Y;\Z)_n \oplus H_1(\bar{Y};\Z)_{n-1}\\ 
  }\]
commute. Here the lower horizontal arrow is the natural `cycle map' sending homotopy classes of sections of symmetric product bundles to the homology classes of the cycles swept out by the points in their support. The map $\alpha$ is surjective, and if $n-1\geq 2$ also injective. Hence the composite
\[\sigma := \tau_X^{-1}\circ \alpha\colon \pi_0 \sect(\EuScript{Q}) \to \spinc(X)_d\] 
is surjective, and bijective if $n-1\geq 2$.
\end{Prop}
\begin{pf}
The elements of $H_2^{(n)}$ are represented by surfaces with boundary, of form $C_1\cup\dots\cup C_{r-1}\cup C'$, where the $C_i$ are cylinders which are sections of the broken fibration $X\to A$, and $C'$ is a cylinder with $\partial C' = Z \cup \gamma$ for some circle $\gamma \subset Y$. Recall the attaching surface $Q\subset Y$ (torus of Klein bottle) which collapses to $Z$. The set $\pi_0\sect(\EuScript{Q})$ maps injectively to $H_1(Q;\Z)$, and in this way it may be considered as an affine abelian group, isomorphic to  $\Z$ or $\Z/2$ according to whether or not $Q$ is orientable. There is an `exact sequence'
\[  0\to \pi_0 \sect(Q) \to H_2^{(n)}\to H_1(\bar{Y};\Z)_{n-1}\to 0 \]
(we use quotation marks because these are only affine abelian groups).

Notice that $H_2(X,\partial X\cup Z;\Z)$ injects into $H_1(\partial X;\Z)$ (and even into $H_1(Y;\Z)$). Thus $H_2^{(n)}$ defines a correspondence between $H_1(Y;\Z)_n$ and $H_1(Y;\Z)_{n-1}$. Because of this, $\alpha$---if it exists---is fully determined by the commutative square in the statement of the theorem. We must show that the image of any section of $\EuScript{Q}$ in $H_1(Y;\Z)_n\oplus  H_1(\bar{Y};\Z)_{n-1}$ lies in $\partial (H_2^{(n)})$. 

To do so, we study $\pi_0\sect (\bar{Y}^{[n-1]})$, and relate it to $ H_1(\bar{Y};\Z)_{n-1}$. In general, given a self-diffeomorphism $\psi\in \diff(M)$ of a connected manifold $M$, once one chooses a reference section of the mapping torus $\torus(\psi)$, the set of homotopy classes of sections, $\pi_0 \sect\torus(\psi)$, is described by a natural commutative diagram
\[\xymatrix{ 
\pi_0\sect \torus(\psi) \ar[r]^{\cong}\ar[d] & \pi_1(M)/\sim \ar[d] \\
H_1(\torus(\psi);\Z)&  H_1(M)/ \im(1-\psi_*) \ar@{^{(}->}[l]   
} \] 
where $\sim$ is the relation of `$\psi$-twisted conjugacy' in $\pi_1(M)$:
\begin{equation}
\label{twisted conj}
 \gamma \sim  \beta\cdot \gamma \cdot (\psi\circ\beta)^{-1}
 \end{equation} 
(concatenation of paths is read from left to right). This applies to $\bar{Y}$, the mapping torus of $\bar{\phi}\in \diff^+(\bar{\Sigma})$. In that case, the vertical map on the right (the abelianisation map) is always surjective, but when the genus of the fibre is $>1$ it is not injective. 

When we consider the symmetric products $\bar{Y}^{[n-1]}=\sym^{n-1}_{S^1}(\bar{Y})$, we can slightly modify the targets of the vertical arrows in the commutative diagram:
\[\xymatrix{ 
\pi_0\sect (\bar{Y}^{[n-1]}) \ar[r]^(.4){\cong}\ar[d] & \pi_1(\sym^{r-1}(\bar{\Sigma}))/\sim \ar[d] \\
H_1(\bar{Y}) &  H_1(M)/ \im(1-\bar{\phi}_*). \ar@{^{(}->}[l]   
} \] 
Here we are  considering sections of $\bar{Y}^{[n-1]}$ as $(n-1)$-fold sections of $\bar{Y}$. If $n-1\geq 2$, $ \pi_1(\sym^{n-1}(\bar{\Sigma}))$ is abelian (isomorphic to $H_1(\bar{\Sigma};\Z)$), and the vertical map on the right is a bijection. 

We conclude that $\pi_0\sect(\bar{Y}^{[n-1]})$ maps surjectively to $H_1(\bar{Y};\Z)_{n-1}$, and also injectively if $n-1\geq 2$.

For any section of $\gamma \in \sect (\bar{Y}^{[n-1]} )$,  $\gamma^*\EuScript{Q}$ is an $S^1$--bundle over $S^1$, and thus $\gamma$ lifts to a section of $\EuScript{Q}\to S^1$. If the surface $Q$ is a torus, the pullback bundle  $\gamma^*\EuScript{Q}$ is always orientable, and the available lifts (up to homotopy) are parametrised by $\Z$. If $Q$ is a Klein bottle, it is never orientable, and the lifts are parametrised by $\Z/2$ (cf. Part I, section 4). In either case, the lifts are in bijection with $\pi_0 \sect(Q)$. Now $\bar{Y}$ contains a  distinguished braid $B$ (which degenerates inside $X$ to the critical circle $Z$). The section $\gamma$ is homotopic to a section $\gamma'$ whose image in $\bar{Y}$ is well away from $B$. We noted in Part I (Remark 4.2) that, after moving $\EuScript{Q}$ by an isotopy, and restricting to an open subset $U\subset \bar{Y}^{[n-1]}$ of points $[x_1,\dots,x_{n-1}]$ where none of the $x_i$ is close to $B$, we can suppose that the $S^1$-bundle $\EuScript{Q}\to\bar{Y}^{[n-1]}$ is identified with $\bar{Y}^{[n-1]}\times_{S_1}Q$. That is to say, to lift $\gamma'$ to a section $\widetilde{\gamma}'$ of $\EuScript{Q}$, one just has to add a section of $Q$. But that means that $\alpha(\widetilde{\gamma}')\in H_2^{(n)}$ exists. 

Hence the map $\alpha$ exists, and we obtain a commutative diagram 
\[\begin{CD}
  \pi_0 \sect(Q) @>>> \pi_0\sect (\EuScript{Q}) @>>> \pi_0 \sect (\bar{Y}^{[n-1]} ) \\
 @| @V{\alpha}VV @VVV\\
 \pi_0 \sect(Q) @>>> H_2^{(d)}@>>> H_1(\bar{Y};\Z)_{n-1}
 \end{CD} \] 
where the upper row is a `short exact sequence' of sets, the lower row a short exact sequence of affine abelian groups. We have already established the surjectivity (and when $n\geq 3$, injectivity) of the right-hand vertical arrow. This gives the result.
\end{pf}

If $\bar{\Sigma}$ is disconnected, one should beware that $\sect(\bar{Y}^{[n-1]})$ (and hence $\sect(\EuScript{Q})$) may be empty, because $\bar{\phi}$ permutes the components of $\bar{\Sigma}$. However, we can still set up a commutative diagram
\[\xymatrix{
   \pi_0 \sect(\EuScript{Q}) \ar@{-->}[r]^{\alpha} \ar[d]   & \bigcup_{0\leq d_1\leq d}H_2^{(d_1,d-d_1)}\ar@{_{(}->}[d]^{\partial}\\ 
    \pi_0 \sect(Y^{[n]}) \times \pi_0 \sect(\bar{Y}^{[n-1]})\ar[r]   	&  H_1(Y;\Z)_n \oplus H_1(\bar{Y};\Z)_{n-1}\\ 
  }\]
which uniquely characterises $\alpha$. This again leads to composite maps
\[\sigma = \tau_{X}^{-1} \circ \alpha\colon \pi_0 \sect(\EuScript{Q}) \to \bigcup_{0\leq d_1\leq d}\spinc(X)_{(d_1,d_2)}. \]
If $n-1\geq 2$, $\alpha$ (and hence $\sigma$) is injective. In general, however, it is not surjective.

\subsubsection{Sections of relative Hilbert schemes versus $\mathrm{Spin}^c$-structures} 
Now take any broken fibration $(X,\pi)$ over a surface $S$. Consider the boundary map 
\[\delta\colon H_2(X,\partial X \cup Z;\Z) \to H_1(Z;\Z). \] 
Let $\beta \in \delta^{-1}([Z])$. The intersection number of $\beta$ with a regular fibre $X_s$ is not constant, and nor is the Euler characteristic of the fibre, but one has
\begin{Lem}
The number $2 \beta \cdot X_s + \chi(X_s)$ is independent of $s\in S\setminus S^\crit$.
\end{Lem}
\begin{pf}
The Lefschetz dual $D(\beta)$ of $\beta$ lies in $H^2(X\setminus Z;\Z)$. The submanifold $X\setminus Z$ supports an almost complex structure compatible with $\omega$ and preserving the fibres of $\pi$, and whilst $c_1(X\setminus Z)$ does not extend to $X$, $c_1(X\setminus Z)+2 D(\beta)$ does (cf. Taubes \cite{Ta2}). The result follows.
\end{pf}
Note that one can define the Taubes map 
\[\tau\colon  \spinc(X) \to  H_2(X,\partial X\cup Z;\Z)  \] 
precisely as in (\ref{taubes map}). As before, let $\nu \colon S\setminus S^\crit \to \Z_{\geq 0}$ be a function satisfying  $2\nu(s)  + \chi(X_s) = d$, and form the relative Hilbert scheme $X^{[\nu]}$ and its matching condition $\EuScript{Q}$. Then any section $u\in \sect(X^{[\nu]},\EuScript{Q})$ defines a relative homology class 
\[\alpha(u) \in \delta^{-1}([Z])\subset H_2(X,\partial X\cup Z;\Z). \]
Indeed, over $S'=S\setminus \mathrm{nd}(\pi(Z))$, $u$ tautologically defines a cycle $u'$ representing a class in $H_2(X , \pi^{-1}(\partial  S);\Z)$. This extends to a class in $H_2(X,Z;\Z)$ as in Prop. \ref{def alpha}. 

The composite of $\tau^{-1}$ and $\alpha$ defines the map
\[ \sigma_\nu = \tau^{-1}\circ \alpha \colon \pi_0\sect(X^{[\nu]} , \EuScript{Q} ) \to \spinc(X)_d  \]
introduced in (\ref{sigma}).
 
\begin{Prop}
Assume $\inf(\nu)\geq 2$. Then the map $\sigma_\nu$ is injective; if all fibres are connected, it is also surjective.
\end{Prop} 
\begin{pf}
We claim that, if $S=S_1 \cup _T S_2$ is a splitting along a circle of regular values $T$, and $X= X_1\cup_Y X_2$ is the corresponding splitting, then injectivity of $\sigma_\nu$ for $X$ is implied by injectivity for $X_1$ and $X_2$. Indeed, suppose sections $s_1$, $s_2 \in \sect (X^{\nu},\EuScript{Q})$ satisfy $\sigma_\nu(s_1)=\sigma_\nu(s_2)$. Then, by injectivity of $\sigma_\nu$ for the $X_i$, $s_1$ is homotopic to a new section $s_1'$ which differs from $s_2$ only over a small neighbourhood of $T$. The obstruction to extending the homotopy over $T$ is a class $\delta\in \pi_1(Y^{[n]},\gamma)$, where $n=\nu(s)$, $s \in T$, and $\gamma = s_2|T$. But there is a natural injection $i \colon \pi_1(Y^{[n]},\gamma) \to H_2(Y;\Z)=H^1(Y;\Z)$, and $i(\delta)$ is the lift of $\sigma_\nu(s_1)-\sigma_\nu(s_2)$ to $H^1(Y;\Z)$. Hence $\delta =0$.

Similarly, surjectivity of $\sigma_{\nu}$ for $X_1$ and $X_2$ implies surjectivity for $X$, providing that $\inf \nu \geq 2$. The point is that $\pi_0\sect(Y^{[n]})$ maps surjectively to $H_1(Y;\Z)_n $ if $n\geq 2$ ($\pi_0\sect(Y^{[n]})$ is an affine space modelled on the abelianisation of $\pi_1(\sym^n(\Sigma))$, where $\Sigma$ is the fibre of $Y\to S^1$).

With this patching procedure to hand, it suffices to check the result in three cases: (i) $S$ is a connected surface with non-empty boundary and $X^\crit =\emptyset$; (ii) $S$ is a disc and $X\to S$ is an elementary Lefschetz fibration; (iii) $S$ is an annulus and $X\to S$ an elementary broken fibration. Case (i) is easy once one observes that $S$ retracts to a bouquet of circles; (ii) is also straightforward and left to the interested reader. Proposition \ref{def alpha}, and the remarks following its proof, deal with (iii). The result follows.
\end{pf}

\subsection{Computing the index} The map $\sigma_\nu \colon \pi_0 \sect(\EuScript{Q})\to \spinc(X)_d$ gives a tidy way of partitioning the sections of $\EuScript{Q}$, but its importance goes beyond book-keeping: we can use it to understand the index problem for pseudo-holomorphic sections with boundary on $\EuScript{Q}$. Whilst setting up invariants from such sections depends critically on compactness and transversality theorems, the index computation can be done in a `soft' topological context. Underlying this and virtually all such  computations is the invariance of the Fredholm index under compact perturbations. 

The proof of the index formula is rather complicated. As a warm-up, and so as to 
make clear which aspects are tied specifically to our matching conditions, we first consider the index problem for the relative Hilbert scheme $X^{[n]} = \hilb^n_S(X)$ of a \emph{Lefschetz} fibration $\pi\colon X\to S$ over a closed surface. In the case $S=S^2$, this was done by Smith \cite[Proposition 4.3]{Smi} using an analytic form of Grothendieck--Riemann--Roch. It can be done in more elementary fashion as follows.

\begin{Lem}\label{Lefschetz index}
Let $u\in \sect(X^{[n]})$, and let $A\in H_2(X;\Z)$ be its underlying homology class. Then the index of $u$ (which, according to Riemann--Roch, is equal to $2 \langle c_1(\Tv X^{[n]}),[u]\rangle +n\chi(S)$) is 
\[ \ind(u) = A\cdot A + \langle c_1(TX), A\rangle.\]
\end{Lem}

The proof is preceded by an algebraic digression. Suppose that $X$ is actually a complex surface (not necessarily compact) and $\pi\colon X\to S$ a holomorphic submersion (not necessarily proper). To give a holomorphic  section $u\colon S\to X^{[n]}$ of the relative symmetric product $\sym^n_S(X)$ is to give a commutative triangle
\[ \xymatrix{
 T_u \ar[r] ^{\tilde{u}}\ar[d]^{t_u} & X  \ar[dl]^{\pi}  \\
  S
} \]
where $T_u$ is a Riemann surface, $t_u$ a proper holomorphic map of degree $d$, and $\tilde{u}$ a holomorphic mapping into $X$. One constructs $T_u$, as a branched cover of $S$, as the space of germs of sections of $\pi$ mapping each point $s$ into $\supp(u(s))$.

\begin{Lem} For any holomorphic section $u\colon S\to X$, there is an isomorphism of holomorphic vector bundles over $S$,
\[  u^*\Tv \sym^n_S(X)  \cong t_{u*} \tilde{u}^*(\Tv X). \]
\end{Lem}
\begin{pf}
This is essentially the infinitesimal version of the correspondence just described. The standard algebro-geometric description of $T \sym^n(X_s)$ at an effective divisor $D$ is as the vector space of `Laurent tails' $H^0(X_s;\mathcal{O}_{X_s,D}(D))$ (here $\mathcal{O}_{X_s,D}(D)$ is the cokernel of the natural inclusion of sheaves $\mathcal{O}_{X_s}\to \mathcal{O}_{X_s}(D)$). Thus the fibre of $u^*\Tv \sym^n_S(X)$ at $x$ is  $H^0(X_s;\mathcal{O}_{X_s,D}(D))$, where $D$ is $u(x)$, counted with multiplicities. On the other hand, the fibre of $t_{u*} \tilde{u}^*(\Tv X)$ at $x$ is the space of $D$--jets of holomorphic maps into $X_s$ at $D$, which is also naturally identified with $H^0(\mathcal{O}_{X_s,D}(D))$.
\end{pf}

\begin{pf}[Proof of Lemma \ref{Lefschetz index}]
Because $u$ is a smooth section, it does not map to any critical points of $\pi^{[n]}$. The index only depends on $u$ only through its homotopy class; hence, by making small perturbations to $u$, we may arrange that $u$ has generic properties. Namely, we can assume that it is transverse to each stratum of the diagonal in $X^{[n]}$, and hence only intersects the diagonal in its top stratum (i.e., there are no triple self-intersections). Moreover, we may suppose that, for any such intersection point $u(s)$, there is a small disc $D_s\subset S$ centred at $s$ over which $u$ is either holomorphic or antiholomorphic, according to the sign of the intersection: having identified the self-intersection points $u(s)$, we adjust the almost complex structure on $X$ so that it becomes integrable in a neighbourhood of $X_s$. 

In this situation, we can again construct a `parametrisation' of $u$ as a commutative diagram
\[  \xymatrix{
 T_u \ar[r] ^{\tilde{u}}\ar[d]^{t_u} & X  \ar[dl]^{\pi}  \\
  S  . 
  }\]
Here $t_u\colon T_u\to S$ is a degree $n$ map of smooth surfaces whose fibre over $x$ is $\supp( u(x))$, \emph{except} when $u(x)$ hits the diagonal. In the latter case,  $u$ is locally (anti)holomorphic, and $T_u$ can be constructed in the way discussed in the preamble to this proof. Thus $t_u$ has simple ramification. We still have an isomorphism of complex vector bundles $u^*\Tv \sym^n_S(X)  \cong t_{u*} \tilde{u}^*(\Tv X)$: the previous lemma gives such an isomorphism near the branch values, and there is then an obvious (and canonical) extension to the rest of $S$.

The image of $ \tilde{u}$ is an embedded surface $U$ in $X$; the branch points of $t_u$ (in $T_u$) map to tangencies with the fibres of $X$. Of these tangencies, some ($N_+$, say), have coincident orientations with the fibre of $X$, and some (say $N_-$) have opposite orientation. The two local models for $t_u$ are respectively $z\mapsto z^2$ and $z\mapsto \bar{z}^2$. In this context, the Riemann--Hurwitz formula reads
\[ \chi (T_u)  = n \chi(S) - (N_+ + N_-). \]
On the other hand, for a complex line bundle $L\to T_u$,
\begin{equation}\label{pushforward}
 \langle c_1( t_{u*}L), [S] \rangle  =   \langle c_1( L), [T_u] \rangle + \frac{1}{2}(N_+ - N_-). \end{equation}
We now calculate
\begin{align*}
\ind(u) &= 2 \langle c_1(u^*\Tv X^{[n]}),[S]\rangle +n\chi(S) 
							&& \text{(Riemann-Roch)} \\
&= 2  \langle \tilde{u} ^*c_1(\Tv X) , [T_u] \rangle  +  N_+ - N_- + n \chi(S)  
							&& \text{(by equation (\ref{pushforward}))}\\
& =  2 \langle \tilde{u} ^*c_1(T X)  - \tilde{u}^* \pi^*c_1(TS) , [T_u] \rangle +   N_+ - N_-  + n\chi(S) 				
					&& \text{(since } TX=\Tv X \oplus \pi^*TS	\text{)}\\
& =  2 \langle \tilde{u} ^*c_1(T X) , [T_u] \rangle -  n \chi(S) +  N_+ - N_-  
							&&\text{(since} \deg(t_u)=n\text{)}) \\
& =  2 \langle \tilde{u} ^*c_1(T X) , [T_u] \rangle -  \chi(T_u)  - 2 N_-
							&& \text{(Riemann--Hurwitz).}
\end{align*}
We can now make the link with the homology class $A$, since
\[\langle \tilde{u} ^*c_1(T X) , [T_u] \rangle =  \langle c_1(TX), A\rangle.\] 
If $N_-=0$, we can homotope the almost complex structure on $X$ so that it preserves $TU$, whereupon the adjunction formula gives 
\[ - \chi(T_u) = A\cdot A - \langle c_1(TX), A\rangle .\] 
If $N_- >0$ then we may take $U$ to be almost complex except at the $N_-$ orientation-reversing vertical tangencies. There is a relation 
\[- \chi(T_u) + k N_-  = A\cdot A - \langle c_1(TX), A\rangle\] 
for a \emph{universal} integer $k$; indeed, $k$ is the first Chern class of the trivial $\C^2$-bundle over the disc, $\Delta$, relative to a certain trivialisation over its boundary. Hence
\[ \ind(u) = \langle c_1(TX),A\rangle + A \cdot A - (k+ 2)N_-. \]
The cheapest way to deduce that $k=-2$ is to observe that (there is some example in which) one can deform an initial $u_0\in \sect(X^{[n]})$ which has $N_-=0$ to a homotopic section $u_1$ which is generic and has $N_->0$. We have $\ind(u_0)=\ind(u_1)$, but the only way this can occur is if $k=-2$. The result now follows.
\end{pf}
The $\spinc$-structure $\mathfrak{s}= \tau_X^{-1}(A)$ corresponding to $A$ is $L_A \otimes \mathfrak{s}_{\mathrm{can}}$, where $L_A$ is the line bundle with a section with transverse zero-locus $A$; thus $c_1(\mathfrak{s}) = c_1(TX) + 2 \mathrm{PD}(A)$. Since  $c_1(TX)^2 = 2e(X)+ 3\sigma(X)$, one has $A^2 + \langle c_1(TX),A\rangle =  d(\mathfrak{s})$, where $d(\mathfrak{s})$ is as in equation (\ref{formal dimension}).

We now turn to broken fibrations, and to one of the results previewed in Part I of this pair of papers.
\begin{TheoremD}\label{index formula}
Suppose that $u\in \sect(X^{[\nu]},\EuScript{Q})$ represents the $\spinc$-structure $\mathfrak{s}$, i.e. that $\sigma_\nu(u)=\mathfrak{s}$. Then the index of $u$ is
\[ \ind(u) = d(\mathfrak{s}). \]
\end{TheoremD}

A part of the proof is to verify the result by hand for each of two broken fibrations
(and at least one section $u$ in each case). The attaching surface in one of these fibrations is a torus; in the other, a pair of Klein bottles. We do this now.
\begin{Ex}\label{torus ex}
This example is a broken fibration $\pi\colon X\to S^2$. The total space is $X=X_+\cup X_0\cup X_- $ where $X_+$ is a trivial $T^2$-bundle over a north-polar disc $D_+$; $X_-$ is a trivial $S^2$-bundle over the south-polar disc $D_-$; and these are joined in the trivial, `untwisted', way by a broken fibration $X_0$ over the equatorial annulus. The vanishing surface $Q\subset \partial X_+$ is a torus, and there is a section $u$ of $(D_+,\partial D_+)\to (X_+,Q)$ of Maslov index 0. Thus $\ind(u) =1$. In fact, for a standard almost complex structure, the moduli space of pseudo-holomorphic sections is $S^1$. 

It is easy to see that $e(X)=2$ and $b_1(X)=1$, so $b_2(X)=2$  and, since there is a section, $\sigma(X)=0$ (in fact, as shown in the last section of \cite{ADK}, $X\cong (S^1\times S^3)\# (S^2\times S^2)$). The section $u$ can be completed to a section $v$ over $S^2$, which we may take to be pseudo-holomorphic (adjusting the almost complex structure on $X\setminus Z$). We have $[v]\cdot [v] = 0$, so by the adjunction formula, $\langle c_1(X\setminus Z) , [ v] \rangle = 2 $. 

The $\spinc$-structure $\mathfrak{s}$ associated with $u$ has $c_1(\mathfrak{s}) = 2 [\mathrm{fibre}]+ 2 [v]$. Note here that, on $X\setminus Z$, $c_1(\mathfrak{s})=c_1(X\setminus Z)+ 2\mathrm{PD} [u]$, and that one can calculate the intersection of $c_1(X\setminus Z)+ 2\mathrm{PD} [u] $ with a surface in $X\setminus Z$ as a sum of two intersection numbers. Thus $c_1(\mathfrak{s})^2 = 8$, and $d(\mathfrak{s}) = 2 - 1 = 1$.
\end{Ex}
\begin{Ex}
In this example of a broken fibration $X\to S^2$ (a homotopy--$S^2\times S^2$), the vanishing surface $Q$ is the union of two Klein bottles. We have $X=X_+\cup X_0\cup X_- $, where $X_+\to D_+$ and $X_-\to D_-$ are now trivial $ S^2$--bundles. The equatorial part $X_0$ is a broken fibration with two circles of critical points, mapping to parallel circles (the `tropics'). The fibre over the equator is $S^2\cup S^2$, and the monodromy around the equator interchanges the two components.  The vanishing Klein bottle $Q_- \subset S^2\times \partial D_- $ is $\bigcup_{z\in S^1}{l_z\times \{z\}}$, where $l_z$ is the circle $z^{1/2}\R\cup \{\infty\}\subset \C\cup \{\infty\} = S^2$. The fibration is symmetric under reflection in the equator in the base $S^2$; in particular, the other Klein bottle $Q_+$ is the reflection of the first.

We have $\pi_1(X)=\{1\}$ (as one verifies with a little care) and $e(X)=4$, so $b_2(X)=2$. There exists a two-fold section, so the fibres are non-trivial in homology and the signature is zero. 

There is a section $u_-$ over $D_-$ with boundary on $Q_-$, given by $u_-(z) = 0 \in \C \subset S^2$. This has Maslov index $+1$. There is a similar section $u_+$ over $D_+$. Their union, $u=u_+\cup u_-$, has index $\ind(u) = 2(1+1) = 4$. 

To calculate $c_1(\mathfrak{s})^2$ for the associated $\spinc$--structure $\mathfrak{s}$, we look more closely at $H_2(X;\Z)$. Observe that there is a two-fold section $s$, given on $D_-$ by $s(z) = \{ 0,\infty \}$, and similarly on $D_+$. In the equatorial annulus, $s(z)$ is given by one point on each of the two components of the fibre $X_z$.
Moreover, we can arrange firstly that  $S:=\im(s)$ is pseudo-holomorphic for some almost complex structure compatible with the broken fibration, and secondly that $S$ is an embedded sphere with $[S]\cdot [S] = 4 $. By adjunction, $\langle c_1(X\setminus Z), [S] \rangle = 6$. Notice that $S$ is homologous inside $X\setminus Z $ to another surface $S'$ transverse to $u$ and meeting $\im(u)$ with (multiplicity $+1$) at $0\in D_+$ and at $0\in D_-$. Thus $\langle c_1(\mathfrak{s}), [S] \rangle = 6 + 2.2 = 10$. We also have $\langle c_1(\mathfrak{s}), F \rangle = 2 $, where $F$ is an equatorial fibre. Hence $c_1(\mathfrak{s})= 2 [S ] + 2 F $. Then $c_1(\mathfrak{s})^2 = 24 $, and $d(\mathfrak{s})=4 = \ind(u)$.
\end{Ex}

\begin{pf}[Proof of Theorem D.]

We will show that $\ind(u)-d(\mathfrak{s})$ is some universal multiple of the number of tori, plus a universal multiple of the number of Klein bottles. The foregoing computations, showing that $\ind(u)-d(\mathfrak{s})=0$ in each of the two examples above, then complete the proof of the theorem.

The argument is a little less intricate when $S=S^2$ and there is just one critical circle $Z$, so we will deal with that case then indicate how it is generalised.
So, $S' = S_+\amalg S_-$, the union of two discs, and there is an orientation-reversing diffeomorphism $\widetilde{\tau} \colon S_+\to S_-$.

One can think of $u$ as a section $(u_+,u_-)$ of $X_+^{[n]}\times_{S_+}\widetilde{\tau}^* X_-^{[-1]}$. The Riemann--Roch theorem for surfaces with boundary gives the index of $u$ as
\[ \ind(u) =  \mu_{\EuScript{Q}}(u) + n+ (n-1),  \]
where $\mu_{\EuScript{Q}}(u)$ is the Maslov index for $u$ relative to $\EuScript{Q}$, when $u$ is considered as a section of $X_+^{[n]}\times \widetilde{\tau}^* X_-^{[n-1]}$.
We compute this index in a number of steps.

{\bf Step 1.} 
Consider the matching condition $\EuScript{Q}$ lying over $\partial S_+$. At a homotopical level, we can make a sequence of simplifications. Let $Y=\pi^{-1}(\partial S_+)$ and $\bar{Y}=\pi^{-1}(\tau(\partial S_+))$.

\begin{enumerate}
\item
$\bar{Y}$ contains a distinguished braid, $B$ a two-fold section of $Y\to S^1$, associated with the broken fibration: see the discussion at the beginning of Section 4.1 in Part I. Using Remark 4.2 from Part I, we can isotope $\EuScript{Q}$ through totally real sub-bundles
so that, over an open set $U\subset \bar{Y}^{[n-1]}$ of points $[x_1,\dots,x_{n-1}]$ with $x_i \neq \mathrm{nd}(B)$, we have $\EuScript{Q} \cong \bar{Y}^{[n-1]} \times_{S^1} Q$ as an $S^1$-bundle over $\bar{Y}^{[n-1]}$. Furthermore, the embedding $\EuScript{Q} \to Y$ sends $([x_1,\dots x_{n-1}];q)$ to $[x_1,\dots, x_{n-1},q]\in Y^{[n]}$.  Such a totally real isotopy leaves Fredholm indices unaltered.
\item
Suppose $\Gamma = (\gamma,\bar{\gamma})\in \sect (\EuScript{Q} $). We can perturb $\bar{\gamma}$ to another section of $\bar{Y}^{n-1}$ which avoids $U$, and lift the homotopy to $\EuScript{Q} $, so obtaining a homotopic section of $\EuScript{Q}$. 
\item After a further (generic) perturbation, $\gamma$ does not intersect the  diagonal in $Y^{[n]}$, and $\bar{\gamma}$ does not intersect the diagonal in $\bar{Y}^{[n-1]}$.
\end{enumerate} 
When $\Gamma$ satisfies point (3), we can canonically associate with it an $n$-fold covering $t_\gamma\colon T_\gamma \to S^1$; an $(n-1)$-fold covering $t_{\bar{\gamma}} \colon T_{\bar{\gamma}}\to S^1$; and maps
\[ \tilde{\gamma} \colon T_\gamma \to Y ,\quad \tilde{\bar{\gamma}}\colon T_{\bar{\gamma}}\to \bar{Y} \]
such that $\gamma = \tilde{\gamma}\circ t_\gamma^{-1}$ and $\bar{\gamma} = \tilde{\bar{\gamma}} \circ t_{\bar{\gamma}}^{-1}$. We then have 
\[ \gamma^* \Tv Y^{[n]} =  (t_\gamma)_* \tilde{\gamma}^{*}\Tv Y,\quad  \bar{\gamma}^* \Tv \bar{Y}^{[n-1]} =  (t_{\bar{\gamma}})_* \tilde{\bar{\gamma}}^{*}\Tv \bar{Y}. \]
By point 2, $T_\gamma$ is isomorphic to a disjoint union $S^1\cup T_{\bar{\gamma}}\to S^1$. Moreover, $\gamma^* \Tv Y^{[n]} \cong q^*\Tv Y\oplus \bar{\gamma}^* \Tv \bar{Y}^{[n-1]} $, where $q\colon S^1\to Q_i$ is the restriction of $\gamma$ to $S^1\subset T_\gamma$. To put this more plainly: $\gamma$ consists of a section $q$ of the vanishing surface $Q_i \subset Y$, together with an $(n-1)$-fold section disjoint from $Q_i$ which matches with the $(n-1)$-fold section $\bar{\gamma}$ of $\bar{Y}$.

With a generic $u\in \sect(X^{[\nu]},\EuScript{Q})$ we can associate, as in the proof of Lemma \ref{Lefschetz index}, a surface with boundary $T_u$, a branched covering $t_u\colon T_u\to S' $, and a smooth map $\tilde{u}\colon T_u\to X'$, such that $\tilde{u}\circ t_u^{-1}=u$ and 
\[u^*\Tv X^{[\nu]} =(t_u)_* \tilde{u}^*\Tv X'. \]
$\Tv \EuScript{Q}$ defines a (linearised) boundary condition for deformations of $u$, and this translates into a boundary condition for deformations of $(T_u, \widetilde{u})$. We understand the boundary condition using point (1): $T_u$ has a distinguished component $C$ over $\partial S_+$, on which the map $\tilde{u}$ is $q$. On this boundary component, the boundary condition is $q^* \Tv Q \subset q^* \Tv Y$. The remaining components of $T_u$ over $\partial S_+$, say $\gamma_1,\dots,\gamma_k$, are matched with the components of $T_u$ over $\partial S_-$. On each pair $\gamma_i$, the boundary condition is the diagonal Lagrangian matching condition $\mathrm{diag}_{\gamma_i^* \Tv \bar{Y}} \subset \gamma_i^* \Tv \bar{Y} \oplus \gamma_i^* \Tv \bar{Y} $. 

{\bf Step 2.} We can go further by joining up the matching ends of $T_u$, eliminating the diagonal boundary condition. We have $S=S'\cup N$, where $N$ is an annulus. We can find a covering $T_{tunnel}\to N$ which joins up with $T_u$ to give a surface with boundary $T_{joined}=T \cup T_{tunnel}$ with a map $t_{joined}\colon T_{joined}\to S$. Moreover, there is a map $v= u_{joined}\colon T_{joined}\to X$, lifting $t_{joined} $ and extending $\tilde{u}\colon T_{u}\to X'$. The boundary $\partial T_{joined}$ maps diffeomorphically to $\partial S_+$ under $t_{joined}$, while $ u_{joined}(\partial T_{joined})\subset Q$.
Now, one has
\[  \mu_{\EuScript{Q}}(u ) =  \mu_Q( v) + (N_+ -N_-),  \]
where $N_+$ and $N_-$ are numbers of branch points of $T_{joined}\to S$, as in the proof of Lemma \ref{Lefschetz index}. (This is a small modification of (\ref{pushforward}).)

{\bf Step 3.} The next step is to cap off $T_{joined}$---not in $X$, but in a manifold $X_{Z,f}$ obtained by surgery on $X$ along $Z$. That is, we choose a framing $f$ for $Z$, excise a small neighbourhood $Z\times D^3$ of $Z$ from $X$, and glue in $D^2 \times S^2$ in its place. Then $\partial T_{joined}$ bounds a standard embedded disc $D$ in $X_Z$. 

Notice that the framings can be chosen according to a standard recipe. If $Q$ is a torus, we may parametrise a neighbourhood of $Z_i\cong S^1$ as $S^1\times D^3$ in such a way that $Q= S^1 \times \{ (x,y,0) : x^2+y^2 =1 \} \subset \partial (S^1\times D^3)$ and such that the curve $\im(q)\subset Q$ is $S^1\times \{  (1,0,0)\}$. If $Q$ is a Klein bottle, we can choose fix some other model in which $Q$ and $q$ appear in a standard way in $\partial (S^1\times D^3)$.

The almost complex structure $J$ on $X\setminus Z$ then extends to an almost complex structure $J'$ on $X_{Z,f}^* := X_{Z,f} \setminus \{p\}$, where $p$ is a point in the added $D^2\times S^2$. There is no obstruction to  choosing $J'$ so that the disc $D$ is pseudo-holomorphic. 

Now $D$ has complex normal bundle, and tangents to $Q$ form a real subbundle over $\partial D$. There is therefore a Maslov index $\mu_Q(D)$, which is actually zero when $Q$ is a torus. We can pin down our universal model in the Klein bottle case by stipulating that the Maslov index (which can be any odd integer) should be $1$. 

Let $W=T_{joined} \cup \bigcup_i{D_i}$ be the capped surface, and $w\colon W\to X_{Z,f}$ the natural map extending $v$. The complex line bundles $v^* \Tv X$ and $N_{D/X_{Z,f}}$ are isomorphic over the common boundary of $T_u$ and $D$, so join together to give a line bundle $L\to W$. We then have
\begin{equation}\label{joined index}  \mu_Q(u_{joined}) + \mu_Q( D) = 2 \langle c_1(L) , [W] \rangle .   \end{equation}
Notice that that there is a degree $n$ map $p\colon W\to S$ extending  $t_u\colon T_u \to S'$. We now compute as in Lemma \ref{Lefschetz index}: 
\begin{align*}
\ind(u) & =  \mu_{\EuScript{Q}} (u) + (2n-1) 		
					&& \text{(Riemann--Roch)}\\
&= \mu_{Q}(u_{joined}) +  (2n-1) + N_+ - N_-  
					&& \text{(by Step 2)}\\
 & =  2\langle c_1(L) , W \rangle  - \mu_Q(D) + (2n-1) +   N_+ - N_- 
 					&& \text{(by equation (\ref{joined index}))}\\
 & =  2 \langle w^*c_1(T X^*_{Z,f})-w^* p^*c_1(TS) ,  [W] \rangle - \mu_Q(D)  + (2n-1)+  N_+ - N_- 
 					&& \text{(splitting } w^*TX^*_{Z,f}\text{)} \\
& =  2 \langle w^*c_1(T X^*_{Z,f}) ,  [W] \rangle- \mu_Q(D)   -2n - 1 +  N_+ - N_-  						&& (\deg(p)=n) \\
& =  2 \langle w^*c_1(T X^*_{Z,f}) ,  [W] \rangle - \mu_Q(D)  - \chi(W)  -  2 N_-  - 1					&& \text{(Riemann--Hurwitz)}\\
& =  \langle w^*c_1(T X^*_{Z,f}) ,  [W] \rangle +  w_* [W]\cdot w_*[W]  - \mu_Q(D)   - 1.		
					&&\text{(adjunction)}
\end{align*}
The last line uses the variant of the adjunction formula discussed in the proof of Lemma \ref{Lefschetz index}.

Thus $\ind(u)$ differs from $ [W]^2 + \langle c_1(TX^*_{Z,f}) ,W\rangle$ by $\mu_Q(D)   + 1$, which is $1$ or $2$, depending on the topology of $Q$. But $c_1(TX^*_{Z,f})^2$ differs from $2\chi(X_{Z,f})+3\sigma(X_{Z,f})$ just by contributions (`Hopf invariants') from the singularities of $p$ \cite{HH}. These depend only on the topology of $Q$. It follows that
\[ \ind(u)  - \frac{1}{4}[c_1(\mathfrak{s}_W)^2- 2 e(X_{Z,f}) - 3\sigma(X_{Z,f}) ]  \]
depends only on the topology of $Q$, where $\mathfrak{s}_W $ is the $J'$-canonical $\spinc$-structure on $X^*_{Z,f}$, twisted by $\mathrm{PD}(W)$. Hence the same is true of $\ind(u) -d(\mathfrak{s})$. This establishes our claim. 

We now drop the assumptions that $S=S^2$ and that $Z$ is a single circle. By the index formula of Lemma \ref{matched index}, we now have
\[\ind(u) = \mu_{\EuScript{Q}}(u)+ \sum_{S_i\in \pi_0(S')} {\nu(S_i)\chi(S_i)}\]
where $\mu_{\EuScript{Q}}(u)$ is a generalised Maslov index. One can form $T_{joined}$ and $v$ as before, and $\mu_{\EuScript{Q}}(u) = \mu_Q(v)+ (N_+ +N_-) $. The computation now proceeds  as before. It shows, as claimed, that $\ind(u)-d(\mathfrak{s})$ is a universal multiple of the  number of torus components of $Q$, plus a universal multiple of the number of Klein bottle components.
\end{pf}

\subsection{The relative theory: Floer homology}
We now begin to lay out the relative version of the theory---the version for broken fibrations over surfaces with boundary. We begin by setting out the special features of Floer homology for relative symmetric products of fibred three-manifolds. These groups first appear, in connection with Seiberg--Witten theory, in the work of D Salamon \cite{Sa2}. They are studied in the author's thesis \cite{Per} and also by Usher \cite{Ush}, who noticed how to achieve the monotonicity condition explained below.

\subsubsection{Topological sectors} Let $(Y,\pi)$ be a surface bundle over $S^1$, $\Sigma=\pi^{-1}([0])$ its fibre, $j\in \EuScript{J}(\Tv Y)$ a vertical complex structure, and $Y^{[n]}=\sym^n_{S^1}(Y)$ the relative symmetric product. Let $H_1(Y;\Z)_n\subset H_1(Y;\Z)$ be the affine subgroup of classes $\gamma$ with $\gamma\cdot [\Sigma] = n $. As discussed in the last section, there is a `tautological'  map
\[  a \colon \pi_0 \sect(Y^{[n]}) \to H_1(Y;\Z)_n . \]
Assuming $\Sigma$ connected, this map is surjective and, if $n>1$, also injective (see the proof of Proposition \ref{def alpha}).

The map $\tau\colon H_1(Y;\Z)\to \spinc(Y)$, sending $\gamma$ to $\mathrm{PD}(\gamma)\cdot \mathfrak{t}_{\mathrm{can}}$, maps $H_1(Y;\Z)$ bijectively to the subset $\spinc(Y)_{n+\chi(\Sigma)/2}$ of $\spinc$-structures $\mathfrak{t}$ with $\langle c_1(\mathfrak{t}),[ \Sigma] \rangle = 2n+\chi(\Sigma)$. We put 
\[HF_*(Y , \mathfrak{t}) = \bigoplus_{\gamma \in (\tau\circ a)^{-1}(\mathfrak{t}) } HF_*(Y^{[n]}, \gamma).\]
The group $HF_*(Y, \mathfrak{t}) $ should really be notated as $HF_*(Y^{[n]},\pi^{[n]}, \mathfrak{t};\Omega)$, since it depends upon a choice of closed two-form $\Omega$. The matter of choosing $\Omega$ takes us back to the homomorphism $\kappa_{n,\lambda}$ constructed at the beginning of this paper. The monodromy $m$ of $Y\to S^1$ is only defined modulo isotopies, but by choosing a suitable two-form on $Y$ we can find a lift $\tilde{m}$ to the symplectic mapping class group of $\Sigma$. We can then apply $\kappa_{n,\lambda}$ to it. Thus we choose $\Omega$ to be a fibrewise-K\"ahler form whose monodromy represents $\kappa_{n,\lambda}(m)$.

The adjustable parameters in this procedure are the choice of lift $\tilde{m}$ of the monodromy $m$, and the value of $\lambda$. There is no obligation to make the same choices for different $\mathfrak{t}$, but we can pick out a canonical, $\mathfrak{t}$--dependent choice by insisting that a monotonicity condition holds. This is the topic of the following paragraphs.

\subsubsection{Periods and monotonicity}
Consider the mapping torus $\torus(\phi)$ of a self-diffeomorphism $\phi\in \diff^+(M)$ of a connected manifold $M.$ Assume that $\phi(x)=x$; we then use $x$ as a basepoint of $M$, and $c_x$ (the `constant section at $x$') as a basepoint for the twisted free loopspace $\sect(\torus(\phi))$.
\begin{Lem}
The map $e\colon \sect(\torus(\phi))\to M$ given by evaluation at $[0]\in S^1$ gives rise
to an exact sequence
\[  \pi_2(M)/\im(1-\phi_*) \to \pi_1(\sect(\torus(\phi)), c_x) \stackrel{e_*}{\to} \pi_1(M)^\phi \to 1 \]
where $\pi_1(M)^\phi\subset \pi_1(M)$ is the group of invariants.
\end{Lem}
A very similar assertion appears in M. Po\'zniak's thesis \cite{Poz}. 
\begin{pf}
The evaluation map $e$ is a Serre fibration---just like the evaluation map $\EuScript{L}M\to M$ on the free loopspace---and the sequence can be obtained from the associated exact sequence of homotopy groups. However, it can also be derived directly. 

Think of loops in $\sect(\torus(\phi))$ as maps $\Gamma\colon S^1\times [0,1]\to M$ satisfying $\Gamma(t,1)=\phi\circ \Gamma(t,0)$. Given a based loop $\gamma\colon S^1\to M$ representing a $\phi$--invariant class $[\gamma]$, a choice of homotopy $\phi_*\gamma \simeq \gamma$ gives rise to a loop in $\pi_1(\sect(\torus(\phi)), c_x )$, hence $e_*$ maps onto $\pi_1(M)^\phi$.  One obtains a map $d\colon \pi_2(M)\to \pi_1(\sect(\torus(\phi)), \ell )$ by modifying the constant loop at $l$ by a map $(D^2,\partial D^2)\to (M, \{x\})$, where $D^2$ is embedded as a small disc inside $S^1\times [0,1]$. Elements of $\pi_2(M)$ shape $s - \phi_*s$ then lie in the kernel of $d$, and $\im(d)\subset \ker(e_*)$. Moreover, if  $[\Gamma] \in \ker (e_*)$ then, after a homotopy, $\Gamma$ is constant along $S^1 \times \{0\}$, and so after a further homotopy differs from $l$ by a map supported in  in a disc contained interior of $S^1\times [0,1]$. Hence $[\Gamma]\in \im(d)$.  If $[r] \in \pi_2(M)$ lies in $\ker(d)$ then one choose a  nullhomotopy $\{\Gamma_t\}_{t\in [0,1]}$ of $\Gamma_0$ (the image of $r$). Then $e \circ \Gamma_t$ defines a 2-sphere $s$ in $M$, and one has $r\simeq s-\phi_* s$.
\end{pf}
\begin{Prop}
If $n>2$, the natural map 
\[  p \colon  \pi_1 \sect(Y^{[n]},\gamma) \to H_2(Y;\Z) = H^1(Y;\Z),\] 
is injective, with cokernel isomorphic to $\Z/n$.
\end{Prop}
\begin{pf}
When $n>2$ one has $\pi_2(\sym^n(\Sigma))=\Z$: this is true when $n \gg 0$ because $\sym^n(\Sigma)$ is then a projective space bundle over a torus. It then follows for all $n\geq 3$ by descending induction, using the Lefschetz hyperplane theorem. The action of $\phi_*$ on $\pi_2$ is trivial. After homotoping the monodromy of $Y^{[n]}$  so that it has a fixed point, the previous lemma then gives the upper row in the following commutative diagram:
\[\begin{CD}
 0  @>>> \Z @>>>  \pi_1 \sect(Y^{[n]},\gamma) @>>> H_1(\sym^n(\Sigma))^\phi @>>>0 \\
@. 		@VVV  @VVV						@VVV \\
0 @>>> H_2(\Sigma) @>>> H_2(Y) @>>> H_1(\Sigma)^\phi @ >>> 0
\end{CD} \]
The lower row is the short exact sequence describing homology for a mapping torus; the maps in it are the ones induced by inclusion of a fibre and by intersection with a fibre. The vertical maps are the `tautological' ones. The one on the right is an isomorphism;
that on the left is injective, with cokernel $\Z/n$. By the snake lemma, $\pi_1 \sect(Y^{[n]},\gamma) \to H_2(Y;\Z)$ is also injective  with cokernel $\Z/n$.
\end{pf}

The Chern class of $\Tv Y^{[n]}$ is given, in the notation of the introduction, by 
\begin{equation}\label{c1} c_1(\Tv Y^{[n]}) = (c_1(\Tv Y )^{[1]}+ 1^{[2]} )/2.  \end{equation} 
The two sides agree on the fibres $\sym^n(\Sigma)$ since the restrictions of $c_1(\Tv Y)^{[1]} $ and $1^{[2]}$ are respectively $(2-2g)\eta_\Sigma$ and $2n\eta_\Sigma -2\theta_\Sigma$. One way to prove this formula is to use Grothendieck--Riemann--Roch (the calculation is very similar to that of the index in \cite{Smi}). As pointed out to the author by Michael Usher (private communication), it has a more concrete expression: 
\begin{equation}
\langle c_1(\Tv Y^{[n]}), [\Gamma]\rangle = \langle  c_1(\Tv Y ) + 2 \mathrm{PD}(\gamma), p(\Gamma) \rangle / 2  = \langle c_1(\mathfrak{t}_\gamma) , p( \Gamma) \rangle/2, 
\end{equation}
where $\Gamma$ is a loop in $\sect(Y)$ based at $\gamma$ and $\mathfrak{t}_\gamma=\tau(\gamma)$. This can also be verified directly, using the method of Lemma \ref{Lefschetz index}. 

%We record this in a lemma.
%\begin{Lem}
%$\langle c_1(\Tv Y^{[n]}), [\Gamma] \rangle = \frac{1}{2} \langle c_1(\mathfrak{t}_\gamma),p(\Gamma) \rangle$.
%\end{Lem}
%Hence, for any closed, fibrewise symplectic form $\Omega$ on $Y^{[n]}$, one can define Floer homology with coefficients in the power series ring
%\[  \Lambda_Y := \Z[[ H^1(Y;\Z) ]]   \]
%instead of the universal Novikov ring (provided that $\sym^n(\Sigma)$ is weakly monotone, i.e. that  $n\geq g-1$ or $n\leq (g-1)/2$). This works even when $n=1$, because the periods of the action functional then factor through $p$.

Now, Floer homology for an LHF $(T\to S^1,\Omega)$ is at its simplest when the homomorphism
\[ [\Omega]\colon  \pi_1(\sect(T),\gamma)\to \Z  \]
is a multiple of the index homomorphism defined by $ c_1(\Tv T)$. This is the `monotonicity' condition, and has been exploited by a number of authors (most relevantly, by Seidel \cite{Se4} in connection with the mapping class group). In our set-up, monotonicity can be achieved provided that $n\neq - \chi(\Sigma)/2$.
\begin{Lem}\label{monotonicity}
Fix a $\spinc$-structure $\mathfrak{t}= \mathfrak{t}_\gamma \in \spinc(Y)_{2n+\chi(\Sigma)} $, where $ 2n+\chi(\Sigma)\neq 0$. 
Put
\[ W_\lambda =  (1+\lambda n) \, (w_\lambda)^{[1]} - \frac{\lambda}{2} 1^{[2]} \in H^2(Y^{[n]};\R), \]
where 
\[ w_\lambda=  \frac{ \lambda\mathrm{PD}(\gamma) + ( \chi(\Sigma) + 2 n)^{-1} c_1(\mathfrak{t}_\gamma) } {1+\lambda n} ,\quad \lambda > -\frac{1}{n}. \]
Then the following monotonicity relation holds:  $W_\lambda -   \frac{2}{\chi(\Sigma)+2n} c_1(\Tv Y^{[n]})$ vanishes on $\pi_1(\sect(Y^{[n]} ,\gamma)$.
\end{Lem} 
\begin{pf}
Take $w\in H^2(Y;\R)$ with $\langle w ,[\Sigma]\rangle =1$.
On $Y^{[n]}$ we consider the class
 \[ W= (1+\lambda n) \, w^{[1]} - \frac{\lambda}{2} 1^{[2]},\quad \lambda>0  \]
(cf. equation (\ref{linear comb})). One has 
\[ \langle W , [\Gamma]\rangle = \langle (1+\lambda n) \,w - \lambda  \mathrm{PD}(\gamma), p(\Gamma) \rangle .\] 
We have learned that
\[ c_1(\Tv Y^{[n]}, [\Gamma])  = \frac{1}{2} \langle c_1(\mathfrak{t}_\gamma) , p( \Gamma) \rangle. 
\]
If we put $w=w_\lambda= (1+\lambda n)^{-1}( \lambda\mathrm{PD}(\gamma) + ( \chi(\Sigma) + 2 n)^{-1} c_1(\mathfrak{t}_\gamma)$ then we will have, $W=W_\lambda$, $\langle w ,[\Sigma]\rangle =1$, and
\[   \langle W _\lambda - \frac{2}{\chi(\Sigma)+2n} c_1(\Tv Y^{[n]}), [\Gamma]\rangle = 0. \]
 \end{pf}
We can express this lemma in terms of monodromy. Fix an area-form $\alpha$ on $\Sigma$ of total area 1, and let $m\in \pi_0\diff^+(\Sigma)$ be the monodromy of $Y\to S^1$. Lift it to $\tilde{m} \in \aut(\Sigma, \alpha)/\ham(\Sigma,\alpha)$ by specifying that the class $[\alpha_{\tilde{m}}]$ of the closed two-form $\alpha_{\tilde{m}}$ induced by $\alpha$ on $\torus(\tilde{m})$ should equal $w$, as in the lemma. The monodromy of the form $\Omega$ constructed in the lemma is then $\kappa_{n,\lambda}(\tilde{m})$. 

In this paper we have usually insisted that $\lambda$ be strictly positive, so it is worth noting that for the last lemma that would be an unnecessary stricture.

We now go back to the Floer homology group $HF_*(Y,\mathfrak{t})$ which we incompletely defined above. We can complete the definition by  specifying that
the closed, fibrewise-symplectic two form on $Y^{[n]}$ should be drawn from the convex set of closed, \emph{fibrewise-K\"ahler} two-forms, representing the class $W_\lambda$ for some $\lambda \geq 0$.  

The monotonicity property established in the lemma tells us that we can define Floer homology over $\Z$, not $\Lambda_{\Z}$. Indeed, two index-zero trajectories with the same asymptotic limits must also have the same area (the argument is given in \cite{Se4}, for example). 

Thus, adopting these $\Z$ coefficients, \emph{$HF_*(Y,\mathfrak{t})$ is a finitely generated abelian group.} The convexity of the set of allowed two-forms (and the contractibility of the space of allowed vertical almost complex structures) implies that $HF_*(Y,\mathfrak{t})$ is a well-defined group, up to canonical isomorphism.

We have not included $\lambda$ in the notation, and in fact $HF_*(Y,\mathfrak{t})$ is independent of $\lambda \geq 0$, up to isomorphism. Indeed, it is a general principle (proved by Y-J Lee) that fixed-point Floer homology $HF_*(T,\Omega_0)$ does not change under a deformation $\{\Omega_t\}_{t\in [0,1]}$  provided that (i) the groups are well-defined for all $t$, over the same coefficient ring, and (ii) the periods do not change either. Condition (ii) means, more precisely, that the energy (or action) homomorphism $[\Omega_t]\colon K\to \R$ on the subgroup $\ker(c_1(\Tv T))\subset \pi_1 (\sect(T),\gamma)$ (where $\gamma$ is a reference section, acting as  a basepoint), should be constant along the deformation. When this is satisfied, the subgroup $HF_*(T,\Omega_t)_\gamma$, corresponding to the component of $\sect(T)$ that contains $\gamma$, is independent of $t$. The proof uses Lee's delicate bifurcation analysis \cite{LeeTors}; see also the discussions in \cite{Lee, Ush}.

We can make the following conclusion.
\begin{Prop}
Given $Y$ and $\mathfrak{t}=  \mathfrak{t}_\gamma \in \spinc(Y)_{2n+\chi(\Sigma)} $, where either $n\geq g$ or $n\leq (g-1)/2$, there is a Floer homology group $HF_*(Y, \mathfrak{t}) $ which is a finitely generated abelian group, and is an invariant of $(Y,\pi, \mathfrak{t})$.
\end{Prop}

For $Y = S^1\times\Sigma$, with its trivial fibration over $S^1$, the PSS isomorphism relating Floer homology to ordinary homology gives  
\[ HF_*(Y, \mathfrak{t}_n ) \cong H_*(\sym^n(\Sigma);\Z ) , \]
as $\Z/2$--graded abelian groups, where $\mathfrak{t}_n$ is characterised by $c_1(\mathfrak{t}_n) = (2n + \chi(\Sigma)) \mathrm{pr}_2^* [o_{\Sigma}]$. 

One has $HF_*(Y, \mathfrak{t}')=0$ when $c_1(\mathfrak{t}')$ is not the pullback of a class on $\Sigma$.

\subsubsection{Extending closed two-forms defined on the boundary}
There are relative invariants which take the form of homomorphisms with values in the groups $HF_*(Y,\mathfrak{t})$, defined as in (\ref{rel Phi}). 

Suppose $(X,\pi)$ is a broken fibration over a surface-with-boundary $S$. Label the components of $\partial S$ as incoming or outgoing, and decompose  $Y=\partial X$ as $Y_{\mathrm{in}}\cup Y_{\mathrm{out}}$ accordingly.

Recall the procedure discussed in Section \ref{Matching conditions} as a preamble to Theorem \ref{Theorem B}. As in that discussion, we let $S'$ be the complement of a neighbourhood of $\pi(Z)\subset S$ (where $Z$ is the one-dimensional part of the critical set), and $X' = X|S'$.  
Let $\mathfrak{s} \in \spinc(X)$ be a $\spinc$--structure such that $\langle c_1(\mathfrak{s}), [\mathrm{fibre}] \rangle =  2d$.  There is then a locally constant function $\nu\colon S'\to \Z_{\geq 0}$ such that $ \nu(s) + \chi(X_s)/2 = d $
for all regular values $s$, and we can form the relative Hilbert scheme $X^{[\nu]}\to S'$.

Recall that the Floer homology $HF_*(Y_{\mathrm{out}},\mathfrak{t}_{\gamma_{\mathrm{out}}} )$ is defined using a two-form $\Omega_{\mathrm{out}} $ on $Y_{\mathrm{out}}^{[\nu]}$ with 
\[  [\Omega_{\mathrm{out}}] = 
 (1+\lambda \nu) \, (w_\lambda)^{[1]} - \frac{\lambda}{2} 1^{[2]}  \] 
where
\[ w_\lambda=  \frac{ 
\lambda\mathrm{PD}(\gamma_{\mathrm{in}}) 
+ (2d)^{-1} c_1(\mathfrak{t}_{\gamma_{\mathrm{in} }} ) 
} 
{1+\lambda \nu} ,\quad \lambda \geq -\frac{1}{\nu}  . \]
The same goes for the Floer homology associated with the incoming boundary. 
We can then extend $[\Omega_{\mathrm{out}}] $ (and its incoming version) to $X^{[\nu]}$ as the class
\begin{equation}\label{extension}   
(1+\lambda \nu) \, (W_\lambda)^{[1]} - \frac{\lambda}{2} 1^{[2]} \in H^2(X^{[\nu]};\R).    \end{equation}
The notations ${}^{[1]}$ and ${}^{[2]}$ make sense on the relative Hilbert scheme, as distinct from the relative symmetric product: the definition is as at the beginning of this paper, except that the universal divisor is replaced by the universal sheaf. We refer to Part I, Section 3 for further details. 

By Lemma \ref{admissible existence}, this class is represented by an admissible two-form $\Omega$, in the sense of Definition \ref{admissible}, provided that $\lambda>0$. We can moreover arrange that $\Omega$ extends both $\Omega_{\mathrm{out}}$ and its incoming version $\Omega_{\mathrm{in}}$. 

Finally, theorem \ref{Theorem B}  applies to give us a Lagrangian matching condition $\EuScript{Q}$ for $(X^{[\nu]},\pi^{[\nu]},\Omega)$ associated with $\mathfrak{s}$. Notice that, according to this prescription, we choose a different two form $\Omega$, and hence are obliged to use a different $\EuScript{Q}$, for each $\spinc$--structure (which, however, we decline to notate). This is the price of working with the canonical, finitely generated Floer homology groups.

\subsubsection{Relative Lagrangian matching invariants}
Using the last paragraph, we conclude that one gets functorial relative invariants
\begin{equation} 
\EuScript{L}_{(X,\pi)}(\mathfrak{s})\colon HF_*(Y_{\mathrm{in}},\mathfrak{t}_{\mathrm{in}})  \to HF_*(Y_{\mathrm{out}},\mathfrak{t}_{\mathrm{out}}), 
\end{equation}
for any admissible $\spinc$--structure (cf. Part I, Section 1). 
These are defined as an instance of (\ref{rel Phi}), with no marked points, save for a caveat we shall come to momentarily. The fibration is $X^{[\nu]}\to S'$, equipped with a form $\Omega$ representing the class defined by formula (\ref{extension}), as discussed in the last paragraph. The Lagrangian matching condition $\EuScript{Q}$ then arises as in Theorem \ref{Theorem B}.   

The caveat is that \emph{we only count sections  representing the $\spinc$--structure $\mathfrak{s}$}. To be precise, consider the space $\sect(X^{[\nu]},\EuScript{Q}; x_{\mathrm{in}}, x_{\mathrm{out}})$ of sections of the cylindrical-end completion of $(X^{[\nu]},\EuScript{Q})$, asymptotic to fixed sections $x_{\mathrm{in}}$ and $x_{\mathrm{out}}$ over the ends.  There is a natural map
$\sigma_\nu \colon \pi_0 \sect(X^{[\nu]},\EuScript{Q}; x_{\mathrm{in}}, x_{\mathrm{out}}) \to \spinc(X)$, defined in just the same way as in (\ref{sigma}). 

Thus, to define $\EuScript{L}_{(X,\pi)}(\mathfrak{s})$ one counts finite-action, index 0 pseudo-holomorphic sections $v$ of the cylindrical-ends completion of
$(X^{[\nu]},\EuScript{Q})$ such that $\sigma_\nu(v)=\mathfrak{s}$. 

Since we are not using Novikov coefficients, it is important to check that the count is finite. When $\inf(\nu)\geq 2$, $\sigma_\nu(v)$ determines the free homotopy class of $v$ over $X^{[\nu]}$. By `free' we mean that we allow homotopies through sections which do not end on $x_{\mathrm{in}}$ and $x_{\mathrm{out}}$. This claim is a straightforward variant of Proposition \ref{def alpha}. Thus any two sections ($v$ and  $v'$, say) on which $\sigma_\nu$ agree differ by an element $(u_1,u_2)$ of $\pi_1(\sect(Y_{\mathrm{in}}, x_{\mathrm{in}}))\times \pi_1(\sect(Y_{\mathrm{out}}, x_{\mathrm{out}}))$, encoding the difference between free and constrained homotopy. We are only interested in rigid (i.e., index zero) pseudo-holomorphic sections. For these, the indices of $u_1$ and $u_2$ must sum to zero. The monotonicity property 
(Lemma \ref{monotonicity}) then implies that the energies of $v$ and $v'$ are equal. Gromov--Floer compactness then implies that the count is finite.

\subsubsection{Quantum module structure}
By applying (\ref{rel Phi}) to cylinders $Y\times [0,1]\to S^1\times [0,1]$, with one marked point in the base, one finds that $HF_*(Y,\mathfrak{t})$ is a module over $\Z[U]$. Here $U$ acts by quantum cap product with a codimension-two cycle:
\[ U\cdot c = \delta_x\cap c. \]
It decreases degree by $2$ (the degree comes into play when one imposes relative gradings on $HF_*(Y,\mathfrak{t})$). Elements $l\in H_1(\Sigma;\Z)$ also act (decreasing degree by $1$):
\[ l\cdot c =\delta_\gamma \cap c ,\]
where $\gamma $ is a loop in a fibre, representing $l$. Hence $HF_*(Y,\mathfrak{t})$ becomes a module over the ring
\[ \Z[U] \otimes \Lambda^*H_1(Y;\Z). \] 

At this point it is easy to check that the groups $HF_*(Y,\mathfrak{t})$ and the elements  $\EuScript{L}_{(X,\pi)}(\mathfrak{s})$ give a field theory having the familiar properties set out in the introduction to Part I. 

The only point there which we have not yet addressed is the somewhat unorthodox grading 
properties of the groups; we deal with that next.

\section{Geometric gradings for Floer homology groups}\label{Gradings}
Our treatment of the gradings of the Floer homology groups, and of the degrees of the cobordism-maps, is a variation of the conventional method in symplectic Floer theory. It was motivated by treatment of gradings and degrees in Kronheimer--Mrowka's monopole Floer homology theory (see \cite{KM, KMOS}), and it highlights the similarity to that theory.

If $Y$ is an oriented three--manifold, the set $J(Y)$ of homotopy classes of oriented two-plane fields has a convenient algebro-topological description which we now summarise.  
$\spinc$--structures provide a convenient formulation (which seems to originate in \cite{KM}; see also the introduction to \cite{KMOS}), but the substance is in Pontrjagin's homotopy-classification of maps $Y\to S^2$ \cite{Pon}.  

An oriented two-plane field effects a reduction of structure group of $TY \oplus \varepsilon^1$ ($\varepsilon^1$ the trivial real line bundle) from $\SO(3)$ to $\U(1)\subset \U(2)$, and hence determines a $\spinc$--structure. We write $J(Y, \mathfrak{t})$ for the set of homotopy classes of oriented two-plane fields underlying $\mathfrak{t}\in \spinc(Y)$. One obtains a transitive $\Z$-action on $J(T,\mathfrak{t})$ by modifying two-plane fields via automorphisms of $TY$ supported in a ball $B^3\subset Y$: $n\in \Z$ sends $[\xi]$ to $[\alpha_n^* \xi]$, where $\alpha_n\colon (B^3,\partial B^3) \to (\SO(3), \{1 \})$ is a smooth map of degree $2n$, acting on $TY$ via a trivialisation of $TB^3$.
When $c_1(\mathfrak{t}_\gamma)$ is torsion, the $\Z$--action is free; 
otherwise, its stabiliser is the divisibility $ \Div(c_1(\mathfrak{t}))$ of $c_1(\mathfrak{t})$ in $H^2(Y)/\mathrm{torsion}$. We declare that $\Div(c)=0$ when $c$ is a torsion class, so that  $J(T,\mathfrak{t})$ is always identified, up to a shift, with $\Z/ \Div(c_1(\mathfrak{t}))$. 

Now suppose that $X$ is a cobordism from $Y_0$ to $Y_1$, $j_i\in J(Y_i,\mathfrak{t}_i)$ for $i=0,1$, and  $\mathfrak{s}\in \spinc(X)$ a $\spinc$--structure which restricts to $\mathfrak{t}_i \in \spinc(Y_i)$. Write 
\[  j_0 \stackrel{\mathfrak{s}}{\sim} j_1 \] 
if there is an almost complex structure $I$ on $X$ which induces the $\spinc$-structure $\mathfrak{s}$ and preserves representatives for $j_0$ and $j_1$. Given $j_0$ and $\mathfrak{s}$, there is a unique $j_1$ such that  $j_0 \stackrel{\mathfrak{s}}{\sim} j_1$, so this relation defines a map of $\Z$-sets $J(X,\mathfrak{s})\colon J(Y_0,\mathfrak{t}_0)\to J(Y_1,\mathfrak{t}_1)$. 

Kronheimer--Mrowka's monopole Floer homology group $HM_*(Y,\mathfrak{t})$ is graded by $J(Y,\mathfrak{t})$.\footnote{Monopole Floer homology has three basic versions, signified by various decorations, as well as perturbed and local-coefficient variants. The grading properties are common to all of these; we do not, however, consider the completed groups $HM_\bullet(Y,\mathfrak{t})$.} Thus $HM_*(Y,\mathfrak{t})$ is the direct sum of subgroups $HM_j(Y,\mathfrak{t})$ indexed by $j\in J(Y,\mathfrak{t})$; any two homogeneous elements have a relative degree in $\Z/\Div(c_1(\mathfrak{t}))$. The cobordism map for $(X,\mathfrak{s})$ sends $HM_{j_0}$ to $HM_{J(X,\mathfrak{s})(j_0)}$. The following theorem asserts that the same is true in the symplectic Floer homology of symmetric products of fibred three-manifolds.

\begin{Thm}\label{gradings}
Let $(Y_0,\pi_0 , J_0) $ be a bundle of Riemann surfaces over an oriented 1--manifold. Let $\mathfrak{t}_0 \in \spinc(Y)$ be an admissible $\spinc$-structure. Then $HF_*(Y_0, \mathfrak{t}_0)$ carries a functorial grading by $J(Y,\mathfrak{t}_0)$. Given a broken fibration $(X,\pi)$ realising a cobordism from $Y_0$ to $Y_1$, and an admissible $\spinc$-structure $\mathfrak{s}$ with $\mathfrak{s}|Y_i = \mathfrak{t}_i$, the cobordism map
\[ \EuScript{L}_{(X, \pi) } (\mathfrak{s}) \colon HF_*(Y_0, \mathfrak{t}_0 ) \to HF_*(Y_1 ,  \mathfrak{t}_1 )  \]
sends $HF_{j_0}$ to $HF_{j_1}$, where $j_1= J(X, \mathfrak{s})(j_0)$.
\end{Thm}

To prove this result, we attach to each generator $\gamma$ for the Floer complex a member $\theta$ of the $\Z$-set of reductions of structure group of the vector bundle $\gamma^* \Tv Y^{[n]}$ from $\Sp(2n,\R)$ to $\orth(n)$. These `graded sections' $(\gamma,\theta)$ are related to the two-plane fields involved in the monopole grading through a variant of the Pontrjagin--Thom construction.

\begin{Defn}
Let $(T^{2n+1},\pi,\sigma)$ be a locally Hamiltonian fibration (LHF) over $S^1$.
\begin{itemize}
\item 
A {\bf grading} for a section $\gamma\in \sect(T)$  is a reduction $\theta$ of the structure group of the vector bundle $\gamma^* (\Tv T)$ from $\Sp(2n,\R)$ to the diagonally-embedded subgroup $\orth(n)$. Thus $\theta$ is the isomorphism class of a pair consisting of a principal $\orth(n)$--bundle $P\to S^1$ and an isomorphism $P\times_{\orth(n)}\C^n \cong \gamma^* \Tv T$. The pair $(\gamma,\theta)$ is called a {\bf graded section}. 
\item
A grading for $\gamma$ determines a Lagrangian subbundle $\Lambda_\theta\subset \gamma^*\Tv T$, up to homotopy: $\Lambda_\theta = P \times_{\orth(n)}\R^n \subset P\times_{\orth(n)}\C^n =  \gamma^*\Tv T$. We say that $(\gamma,\theta)$ is {\bf even} if $\Lambda_\theta$ is an orientable vector bundle; otherwise it is {\bf odd}. 
\end{itemize}
\end{Defn}
The space $\sect_{gr}(T)$ of graded sections is a covering space of $\sect(T)$.  The fibre is $\pi_1(\Sp(2n,\R)/\orth(n))=\Z$. Given a graded section $(\gamma_0,\theta_0)$ and a path $\Gamma=\{\gamma_t\}_{t\in[0,1]}$ in $\sect(T) $ starting at $\gamma_0$, path lifting gives a grading $\theta_1$ for $\gamma_1$. For any other grading $\theta_1'$, we write $\delta(\theta_0,\theta_1')$ for the difference $\theta_1'-\theta_1$. If $\theta_0$ and $\theta_1'$ are both even, one has
\[  \delta(\theta_0,\theta_1') = 2 c_1(\Gamma;\theta_0,\theta_1'), \]  
where $c_1(\Gamma;\theta_0,\theta_1')$ is the relative Chern number---the Euler number of $\Gamma^*\Tv T$ relative to a trivialisation over the boundary which lifts $\theta_0$ and $\theta_1'$.

The connection between graded sections and index problems in Floer theory begins to emerge when one looks at the set $\hat{J\,}(T,[\gamma])$ of homotopy classes of graded sections of $T$ representing the same class as $\gamma$ in $\pi_0 \sect(T)$. Recall that the indeterminacy in the relative grading of two homotopic sections $\gamma$, $\gamma'$ which are horizontal (and so represent generators for the Floer complex) is $2N_\gamma \Z$, where
\begin{equation}\label{N gamma}  N_\gamma \Z= \langle c_1( \Tv T),\pi_1 (\sect(T),\gamma) \rangle. \end{equation}
But gradings $\theta_0$, $\theta_1$ for the same section $\gamma$ are homotopic if and only if
\[ \theta_1 - \theta_0 \equiv 0 \mod 2 N_\gamma. \]
Thus the action of $\Z$ upon $\hat{J\,}(Y,[\gamma])$ is transitive with stabiliser $2N_{\gamma}$. 

\subsection{Lifting to $\hat{J\,}(T,[\gamma])$} 
We now explain how to assign a grading to each generator for the Floer complex. 

A generator is a section $\gamma$ of $T\to S^1$ which is horizontal with respect to the connection determined by $\sigma$, i.e., $\sigma(\dot{\gamma}(t), v)=0$ for all $t\in S^1$ and all $v\in \Tv_{\gamma(t)}T $. Because $\gamma$ is horizontal, there is a canonical symplectic connection on $\gamma^*\Tv T$. One says that $\gamma$ is non-degenerate if the monodromy of this connection has no fixed points besides zero. We wish to assign a grading to a non-degenerate horizontal section.

What a non-degenerate horizontal section gives us is:
\begin{itemize}
\item
a symplectic vector bundle $\xi \to S^1 $; and
\item
a symplectic connection $\nabla$ on $\xi$ such that the monodromy $L_\nabla$ does not have 1 as an eigenvalue.
\end{itemize}

$\nabla$ determines an $\orth(n)$-reduction of $\xi$ as follows: Choose a symplectic vector bundle $\eta\to \bar{D}$ and an isomorphism 
$\iota\colon \eta|\partial \bar{D}\cong \xi$. By trivialising $\eta$ one obtains, \emph{a fortiori}, a reduction $\theta_\eta$ of $\xi$ to the orthogonal group. However, there is also a Conley--Zehnder index $ \CZ(\eta,\nabla)\in \Z $.\footnote{This number is derived from the usual Conley--Zehnder index for paths in the symplectic group as follows. One chooses an extension of $\iota^*\nabla$ to a symplectic connection on $\xi$ and a trivialisation of $\eta$ along a radial arc, whereupon one obtains a path in the symplectic group from the identity to the monodromy of $\nabla$. This has a CZ index.} 
\emph{We assign to $\gamma$ the unique grading $\theta$ such that} 
\begin{equation} 
\delta(\theta_\eta, \theta) =  \CZ(\eta,\nabla). \end{equation}
Notice that the parity of $\theta$ coincides with the usual parity (or mod 2 Conley--Zehnder index) of $(\xi,\nabla)$, derived from the sign of $\det(\id-L_\nabla)$. 
\begin{Prop}
The grading $\theta$ does not depend on the extension $\eta$. 
\end{Prop}
\begin{pf}
Suppose $\eta_0$ and $\eta_1$ are two such extensions, inducing gradings $\theta_0$, $\theta_1$. Glue the two copies of the disc together along their common boundary, reversing the orientation on the first copy; likewise, glue the two vector bundles together to obtain $\eta\to S^2$. We can evaluate the Chern number $2c_1(\eta)$ in two ways. First, by the additivity properties of the Conley--Zehnder index,
\[ 2c_1(\eta)  =  \CZ(\eta_1,\nabla) - \CZ(\eta_0,\nabla). \]
Second, let $f_0$, $f_1$ be symplectic trivialisations of $\xi$ which lift the gradings $\theta_0$, $\theta_1$ (if both are even) or $\theta_0+1$, $\theta_1+1$ (if both are odd). Then 
\begin{align*}
 2c_1(\eta) &= 2 c_1(\eta_1;f_1)- 2c_1(\eta_0; f_1) \\
 & =  2 c_1(\eta_1;f_1)- 2c_1(\eta_0; f_0) +  (\theta_0- \theta_1) \\
 &=  ( \theta_1 - \theta_{\eta_1}) - (\theta_0-\theta_{\eta_0})  +  (\theta_0- \theta_1)  \\
 & = \CZ(\eta_1, \nabla)-\CZ(\eta_0,\nabla) + (\theta_0- \theta_1)  .  
 \end{align*}
Hence $\theta_1-\theta_0 =0$.
 
\end{pf}

We have therefore established a grading map $\gamma\mapsto \hat{\gamma}=(\gamma,\theta)$. It is evidently compatible with the Floer-theoretic index: Suppose one has an LHF $(E,\pi,\Omega)$ over the cylinder $ S^1\times[0,1]$, realising a cobordism between non-degenerate locally Hamiltonian fibre bundles $(Y_0,\pi_0,\sigma_0)$ and $(Y_1,\pi_1,\sigma_1)$. Let $\Gamma$ be a smooth section such that $\partial \Gamma = \gamma_1-\gamma_0$ is horizontal. Then one has
\[ \delta( \hat{\gamma}_0,\hat{\gamma}_1) =  \CZ(\Gamma).  \] 

\subsubsection{Relation to two-plane fields}
We now specialise to surface-bundles $Y\to S^1$ and their relative symmetric products $Y^{[n]}=\sym^n_{S^1}(Y)$. A section $\gamma\in \sect(Y^{[n]}) $ determines a homology class $[\gamma]\in H_1(Y;\Z)$, and therefore a $\spinc$-structure $\mathfrak{t}_\gamma= \mathrm{PD}(\gamma)\cdot \mathfrak{t}_{\mathrm{can}}$. 

The object of this section is to compare the $\Z$-sets $\hat{J\,}(Y,[\gamma])$ and $J(Y,\mathfrak{t}_\gamma)$. A first task is to compare the stabilisers $2N_\gamma\Z$ and $\Div(c_1(\mathfrak{t}_\gamma))\Z$.

\begin{Lem}
$\Div(c_1(\mathfrak{t}_{\gamma}))$ divides $2N_\gamma$, which in turn divides $2(n+1-g)$. 
\end{Lem}

\begin{pf}
For a class $\Gamma \in H_2(Y^{[n]} ;\Z )$ represented by a loop $\{\gamma_t\}_{t\in S^1}$ of sections $\gamma_t\in \sect(Y^{[n]})$ one has, by Equation \ref{c1},
\[  \langle c_1(\Tv Y^{[n]}), \Gamma\rangle = \frac{1}{2} \langle c_1(\Tv Y) + 2h, p(\Gamma)  \rangle \]
where $p(\Gamma) \in H_2(Y;\Z)$ is the class tautologically determined by $\Gamma$, and $h\in H^2(Y;\Z)= H_1(Y;\Z)$ that determined by $\gamma_0$. Thus $ \langle c_1(\Tv Y^{[n]}), \Gamma\rangle = \langle  c_1(\mathfrak{t}_\gamma) , p(\Gamma) \rangle/2$, and $2N_\gamma$ is a multiple of $\Div(\mathfrak{t}_\gamma)$.

That $N_{\gamma}$ divides $n+1-g$ is due to the fact that $c_{\min}(\sym^n(\Sigma)) = |n+1-g|$. \end{pf}

We define a map 
\[ t \colon  \hat{J\,}(Y,[\gamma]) \to J(Y, \mathfrak{t}_{\gamma} )  \]
in three stages: 
\begin{enumerate}
\item
To define $t(\gamma,\theta)$, we first perturb $\gamma$ so that it does not meet the diagonal.
It then defines an embedded 1--manifold in $Y$, and the grading is represented by a line sub-bundle $l_\theta$ of its normal bundle. 
\item
The $[-\epsilon,\epsilon]$-subbundle of $l_\theta$ exponentiates to a surface in $Y$. The boundary of this surface is a 1-submanifold $\gamma'$ \emph{with a natural (outward-pointing) framing}. We call $\gamma'$ the {\bf framed double} of $(\gamma,\theta)$.
\item
We can twist $\Tv Y$ along $\gamma'$ (as we can for any framed 1-submanifold transverse to $\Tv Y$) using the inverse of the Pontrjagin--Thom construction. The framing determines a tubular neighbourhood $\gamma'\times D^2\hookrightarrow Y$. Projection onto $D^2$, followed by the collapsing map $D^2\to D^2/\partial D^2= S^2$ gives a map of the neighbourhood to $S^2$ which extends to a map $\alpha\colon Y\to S^2$. In the presence of a trivialisation of $TY$, any two plane-fields $\xi$, $\xi'$ differ by a map $\alpha \colon Y\to \SO(3)/\SO(2)=S^2$. Conversely, given a map $\alpha \colon Y\to S^2$, one can twist the plane-field $\xi$ to a new plane-field $\xi'=\alpha\cdot \xi$. Applying this twisting to the neighbourhood of $\gamma'$, with the trivialisation determined by the framing, and the map $\alpha\colon Y\to S^2$, we obtain a plane-field  $\alpha\cdot \Tv Y$. This field is $t(\gamma,\theta)$.
\end{enumerate}

The difference in Euler classes, $e(t(\gamma,\theta))-e(\Tv Y) $, is Poincar\'e dual to $[\gamma']=2[\gamma]$, hence $t(\gamma,\theta)$ does represent $\mathfrak{t}_\gamma$.

\begin{Lem}
$t$ is a map of $\Z$-sets.
\end{Lem}
\begin{pf}
It is evident that there is a universal constant $c$ such that $t(\gamma,\theta+1)=t(\gamma,\theta)+c$. The Pontrjagin--Thom map \emph{is} a map of $\Z$-sets: when one adds a twist to the framing of a submanifold, the resulting two-plane field changes by $+1$. But adding two to the grading $\theta$ results in adding two twists to the framing of $\gamma'$ (one twist for each component, if $\theta$ is even), hence $c=1$. 
\end{pf}

\subsubsection{Cobordism maps}

(i) Suppose $(X,\pi,\Omega)$ is a Lefschetz fibration over a surface with boundary, realising a cobordism between $(Y_0,\pi_0,\sigma_0)$ and $(Y_1,\pi_1,\sigma_1)$. 

Fix sections $\alpha_i \in \sect(Y_i)$; $r\in \mathbb{N}$; and a component $\beta$ of the space of sections $\Gamma$ of $X^{[n]}$ with $\partial\Gamma = \alpha_1-\alpha_0$. A grading of $\alpha_0$ extends to a unique trivialisation of $\Gamma^* \Tv X^{[n]}$, up to homotopy, and this defines a map 
\[  \hat{J\,}(X,\beta) \colon  \hat{J\,}(Y_0, \alpha_0) \to \hat{J\,}(Y_1, \alpha_1). \] 
This map encodes the degree of the cobordism map $HF(X,\beta)$ on Floer homology, in the sense that $\EuScript{L}_{(X,\pi)}(\mathfrak{s}_\beta)$ maps $HF_{\theta}(Y_0,\mathfrak{t}_{\alpha_0})$ into  $HF_{\hat{J\,}_{X,\beta}(\theta)}(Y_1,\mathfrak{t}_{\alpha_1})$. 

\begin{Lem}
The diagram
\[\begin{CD}
 \hat{J\,}(Y_0,[\gamma_0]) @>{\hat{J\,}_{X,\beta}}>> \hat{J\,}(Y_1,[\gamma_0]) \\ 
@V{t}VV		@VV{t}V\\
J(Y_0,\mathfrak{t}_{\gamma_0}) @>{J_{X,\mathfrak{s}_\beta}}>> J(Y_1,\mathfrak{t}_{\gamma_1})\\
\end{CD}  \]
commutes.  
\end{Lem}
\begin{pf}
It suffices to construct an almost complex structure, representing $\mathfrak{s}_\beta$ and preserving the two-plane fields $t(\gamma_i,\theta_i)$, where $\theta_0$ and $\theta_1$ extend to a grading $\{ \theta_t\}_{t\in[0,1]}$ for a homotopy $\Gamma= \{\gamma_t\}_{t\in[0,1]}$. 

We may assume that $\gamma_0$ and $\gamma_1$ are disjoint from the diagonals in $Y_i^{[r]}$, and that $\Gamma$ intersects the strata of the diagonal $\Delta \subset \hilb^n_S(X)$ transversely. Then $\Gamma$ hits only the top stratum of $\Delta$, and the intersections are isolated points. The subset $\Gamma' \subset X$ determined by $\Gamma$ is an embedded surface, and the intersections with $\Delta$ appear as points  $\{p_1,\dots, p_k\}$ where $\Gamma'$ is tangent to the fibres of $X$.

Suppose, for the time being, that $\Gamma$ is disjoint from the diagonal. We may then take $\Gamma'$ to be pseudo-holomorphic, with respect to an almost complex structure $I$ on $X$ preserving the vertical distribution $\ker(D\pi)$. The idea is to define an new almost complex structure by `twisting along the double of $\Gamma'$'. The procedure is similar to the one used to define the map $t$. One obtains a new almost complex structure $\alpha\cdot I := \alpha I \alpha^{-1}$ from $I$ from any map $\alpha\colon X \to \SO(TX)/\U(TX,I)$. But along $\Gamma'$ we can choose a unitary trivialisation of $TX$ (by trivialising $T\Gamma'$ and $\Tv X|\Gamma'$). Then, near $\Gamma'$, $\SO(TX)/\U(TX,I)$ is identified with $\SO(4)/\U(2)=S^2$. In a tubular neighbourhood $D^2\times \Gamma'$ of $\Gamma'$, we take $\alpha$ to be the composite $( D^2\times \Gamma' , \partial D^2 \times \Gamma' ) \to (D^2,\partial D^2) \to (S^2,*) \to (\SO(4)/\U(2), [\id])$. This defines a twisted almost complex structure $I'=\alpha\cdot I$, compatible with the twisted tangent distribution on $\partial X$.  

It remains to deal with the vertical tangencies. The failure of the square to commute is given $N  \Gamma \cdot \Delta $, where $N$ is some constant. It follows that, for Lefschetz fibrations on closed four-manifolds, one has $\ind(\Gamma) - d(\mathfrak{s}_\Gamma) = N  \Gamma \cdot \Delta $. Hence $N=0$ by Lemma \ref{Lefschetz index} (the Lefschetz fibration case of Theorem D). 
\end{pf}

The proof of Theorem \ref{gradings} is now almost complete; all that remains is to check functoriality for broken fibrations. If $(X,\pi)$ is a broken fibration, realising a cobordism from $Y_0$ to $Y_1$, and $\beta\in \pi_0  \sect(X^{[\nu]},\EuScript{Q})$, we define $\hat{J\,}(X,\beta)$ by sending $(\gamma_0,\theta_0)$ to the $(\gamma_1,\theta_1)$, where $\gamma_1-\gamma_0 =\partial \Gamma$ for a representative $\Gamma$ of $\beta$, and there is a trivialisation of $\Gamma^* \Tv X^{[\nu]}$ compatible with $\gamma_0$, $\gamma_1$ and $\Tv \EuScript{Q}$. Then $\hat{J\,}(X,\beta)$ defines the degree of the map $\Phi(X,\beta)$ on Floer homology, and we must compare $t\circ \hat{J\,}(X,\beta)$ with $J_{X,\mathfrak{s}_\beta} \circ t$. Again, the fact (Theorem D) that $\ind(\beta)$ equals $d(\mathfrak{s}_\beta)$ when $X$ is closed implies that the discrepancy must be zero.

\section{Calculations} \label{calculations}
We now compute the Lagrangian matching invariants in some simple cases, and show that they coincide with Seiberg--Witten invariants. The calculations presented here are limited to situations where (i) the monodromies around either side of each circle of critical values are trivial, or (ii) one considers $\sym^n$ only for $n\leq 1$.

\subsection{First examples}
\subsubsection{Trivial bundle over a disc}
Let $\Sigma\times \Delta\to \Delta$ be the trivial bundle over the unit disc. Let $\mathfrak{s}_{d}$ be the $\spinc$-structure with $\langle c_1(\mathfrak{s}) , [\Sigma] \rangle= d=\chi(\Sigma)+2n$, and $\mathfrak{t}_d$ its restriction to $\Sigma\times S^1$. When $\mathfrak{s}_{d}$ is admissible, the relative invariant is an element 
\[ \EuScript{L}(\mathfrak{s}_{d}) \in  HF_*(\Sigma\times S^1,\mathfrak{t}_d),    \]
the `fundamental class'. 
The ring $\Z[U]\otimes \Lambda^* H_1(\Sigma;\Z)$ operates on $HF_*(\Sigma\times S^1,\mathfrak{t}_d)$. Indeed, $l \in H_1(\Sigma;\Z)$ acts by $l\cdot x = \mu(l) \cap x$, where $\mu\colon H_1(\Sigma;\Z)\to H^1(\sym^n(\Sigma);\Z)$ is the usual $\mu$-isomorphism (see, e.g., Part I, Section 3.6), and $\cap$ is quantum cap product. The degree two element $U$ acts by $U\cdot x = \eta \cap x$ ($\eta$ as in the introduction to this paper). By the PSS isomorphism, the elements $m\cdot \EuScript{L}(\mathfrak{s}_{d})$ give a basis for $HF_*(\Sigma\times S^1)$ as $m$ runs over monomials
\[ m = U^i \otimes l_1\wedge \dots \wedge l_{n-i},\quad 0\leq i \leq n. \]
So we have an isomorphism
\[  \Phi\colon \Z[U]\otimes \Lambda^* H_1(\Sigma;\Z) \to HF_*(\Sigma\times S^1,\mathfrak{t}_d),\quad \eta^i \otimes l_1\wedge \dots \wedge l_{n-i} \mapsto U^i\cdot (l_1\cdots l_{n-1} \cdot \EuScript{L}). \]
If we reverse the orientation of the disc, we obtain a map $\EuScript{L}^\vee \colon HF_*(\Sigma\times S^1,\mathfrak{t}_d)\to \Z$. This sends $U^n\cdot \EuScript{L}$ to $1$, and all the other basis elements $m\cdot \EuScript{L}$ to $0$.

We can describe the action of $\Z[U]\otimes\Lambda^*H_1(\Sigma;\Z)$ on $HF_*(\Sigma\times S^1,\mathfrak{t}_d)$ by quoting Bertram and Thaddeus' formulae \cite{BT} for the quantum cup product on $\sym^n(\Sigma)$. We have 
\[  l \cdot  \Phi(c) = \Phi(  l \cdot  c) ,\quad l\in H_1(\Sigma;\Z) \]
(the quantum product with classes of odd degree is undeformed, as always). We also have
\[ U \cdot \Phi (c)  =  \Phi( \eta \cdot c  ),\quad n\leq (g-1)/2  \]
for dimension reasons. However,
\[ U \cdot \Phi (\eta^i)  =  \Phi ( \eta^{i+1} + \theta_{g-n+i} - \theta_{g-n} \eta^{i}) ,\quad n\geq g>0, \]
where $\theta_m = \theta^m/m! $ for $m\geq 0$, and $\theta_m=0$ for $m<0 $ (recall that $\theta\in \Lambda^2 H_1(\Sigma)= \Lambda^2H^1(\Sigma)^*$ is the cup-product form). 

When $g=0$, one has $U\cdot \eta^i = \eta^{i+1}+ \eta^{i-n}$ (negative powers are read as zero). In particular,
\[ U^{n+1}  \cdot  = \id. \] 

\subsubsection{The trivial $S^2$--bundle over $S^2$} 
This example---the fibre bundle $\mathrm{pr}_2 \colon S^2 \times S^2\to S^2$---is included just to show how to use the field theory before it gets mixed up with other ingredients. We use the basis $([S^2\times \{\mathrm{pt.}\}], [\{\mathrm{pt.}\}\times S^2 ])$ for homology. Consider the homology class $(m,n)$, $n\geq 1$ and its associated $\spinc$-structure $\mathfrak{s}(m,n)$. It has $c_1(\mathfrak{s})= \mathrm{PD} (2+2m,2+2n)$ and $d(\mathfrak{s}(m,n))= 2[mn +m+ n]$.

We compute $\EuScript{L}_{(S^2\times S^2, \mathrm{pr}_2)} (\mathfrak{s}(m,n))$ by splitting the two-sphere into the southern hemisphere $D_-$ and northern hemisphere $D_+$. We have
\[ \EuScript{L}_{S^2\times S^2, \mathrm{pr}_2} (\mathfrak{s}(m,n)) = \EuScript{L}_{D_+} ( U^{mn+m+n}\cdot \EuScript{L}_{D_-})= \EuScript{L}_{D_+} (U^n \cdot U^{m(n+1)}\cdot \EuScript{L}_{D_-}).  \]
Here $\EuScript{L}_{D_-}$ is the invariant associated with the southern hemisphere $D_-$,
\[ \EuScript{L}_{D_-} \in   HF_*(S^2\times S^1, \mathfrak{t}_{2 + 2n}),  \]
where $\langle c_1(\mathfrak{t}_k),  [S^2 \times \{\mathrm{pt}.\}] \rangle =k $. By the PSS isomorphism, 
\[\{ \EuScript{L}_{D_-}, U\cdot \EuScript{L}_{D_-},\dots, U^n\cdot \EuScript{L}_{D_-}\}\]
is a $\Z$-basis for $HF_*(S^2\times S^1, \mathfrak{t}_{2 + 2n})$. Because of the structure of the quantum cohomology just quoted, $U^{n+1}\cdot \EuScript{L}_{D_-}=\EuScript{L}_{D_-}$, and $U^n \cdot U^{m(n+1)}\cdot \EuScript{L}_{D_-} = U^n\cdot \EuScript{L}_{D_-} =\Phi( [\eta^n]) $ is a primitive element of Floer homology. Moreover, this element maps to $1\in \Z$ under  $\EuScript{L}_{D_+}\colon HF_*(S^2\times S^1, \mathfrak{t}_{2 + 2n})\to \Z$.

Hence  
\[   \EuScript{L}_{S^2\times S^2, \mathrm{pr}_2} (\mathfrak{s}(m,n)) =  \begin{cases} 
  U^{ (m+1)(n+1)-1}, &  n\geq 1, \, m\geq 0 \\
  0, & n\geq 1,\, m < 0.\end{cases}  \]
We can compare this result with the  Seiberg--Witten invariants. By the wall-crossing formula, we have
\[   \mathrm{SW}^\tau_{S^2\times S^2} (\mathfrak{s}(m,n)) =  \begin{cases} 
  U^{ (m+1)(n+1)-1}, &  n\geq 0, \, m\geq 0 \\
  0, & n\geq 1,\, m<0\end{cases}  \]
The superscript $\tau$ designates Taubes' chamber for a symplectic form of shape $\omega_{S^2}\oplus \omega_{S^2}$.

\subsubsection{A broken fibration on $(S^1\times S^3)\, \# \, (S^2\times S^2)$.} 
This broken fibration, $\pi\colon X \to S^2$, was introduced in Example \ref{torus ex}. Recall its structure:
\begin{itemize}
\item
Over the southern hemisphere, $X_-\to D_-$ is a trivial torus-bundle. The vanishing torus $Q\subset\partial X_-$ is $ l\times \partial D_-$, where $l \subset T^2$ is the projection of a line in $\R^2$ to $\R^2/\Z^2$, so $Q$ is Lagrangian with respect to the standard, product symplectic form.
\item
Over an equatorial annulus, $X_0\to A$ is the elementary broken fibration constructed from $Q$.
\item
Over the northern hemisphere, $X_+\to D_+$ is a trivial $S^2$-bundle. There are two possible ways to glue $\partial X_+$ to $\partial X_0$, corresponding to the two elements of $\pi_1(\SO(3))$. We choose the `untwisted' gluing, for which $X= X_+ \cup X_0\cup X_-$ is diffeomorphic to $(S^1\times S^3)\, \# \,(S^2\times S^2)$ (the other gluing gives $(S^1\times S^3)\,\#\, \PS^2 \,\#\, \overline{\PS^2}$; see \cite[Section 8.2]{ADK}). 
\end{itemize}
 For the standard complex structure $J_0$ on $T^2\times D_-$, the only holomorphic sections of $X_-\to D_-$ with boundary on $Q$ are the constant sections $c_x$, $x\in l$, and since $\ind(c_x)=1$, $J_0$ is regular. These constant sections represent a class $\beta \in H_2(X,Z;\Z)$ with $\beta\cdot [T^2]=1$. The corresponding $\spinc$-structure $\mathfrak{s}_\beta$ has $c_1(\mathfrak{s}_\beta)=(2,2)$ in the basis for $H_2(X;\Z)$ given by $(F,S)$, where $F$ is the fibre class and $S$ is represented by a square-zero section. Hence 
\[\EuScript{L}_{X,\pi}(\mathfrak{s}_\beta) = \pm 1\otimes \lambda \in \mathbb{A}(X) = \Z[U]\otimes \Lambda^*H^1(X;\Z) , \] 
where $\lambda$ is the generator for $H^1(X;\Z)=\Z$. This is in agreement with the Seiberg--Witten invariant $\mathrm{SW}_X(\mathfrak{s}_\beta)$, calculated in the Taubes chamber of a compatible near-symplectic form.

We can compute the invariants for more general $\spinc$-structures using our field theory. Note that the Lagrangian matching invariants for $n$-fold sections of $X_-$ and $(n-1)$-fold sections of $X_+$ are well-defined for any $n>0$. 

Write $\mathfrak{s}_{p,q}$ for the unique $\spinc$-structure on $X$ with $c_1(\mathfrak{s}_{p,q}) = (2p,2q)$. The relative homology class $\beta+ (n-1) S+ m F \in H_2(X,Z;\Z)$ has associated $\spinc$-structure $\mathfrak{s}_{\beta+(n-1)S+mF } = \mathfrak{s}_{m+1,n}$. The virtual dimension of the moduli space is $d(\mathfrak{s}_{m+1,n})= 2n(m+1)-1$. 

We shall calculate the invariant of $\mathfrak{s}_{m+1,n}$ ($m\geq 0$) by splitting $S^2$ into its three parts $D_-$, $A$ and $D_+$. By the gluing rule,
 \[ \EuScript{L}_{X,\pi}(\mathfrak{s}_{m+1,n}) =  \EuScript{L}_{D_+} \circ (U^{n(m+1)-1} \cdot\EuScript{L}_A) \circ  (l \cdot \EuScript{L}_{D_-}). \]
The relative invariant of $D_-$ is the element $\EuScript{L}_{D_-}=\Phi(1)\in HF_*(T^2\times S^1, \mathfrak{t}_{2n})$, and  $l \cdot \EuScript{L}_{D_-} =\Phi(l)$.
The relative invariant $\EuScript{L}_A$ for $X_0\to A$ is a map 
 \[  H_*(\sym^n(T^2);\Z) =  HF_*(T^2\times S^1,\mathfrak{t}_{2n}) \to HF_*(S^2\times S^1,\mathfrak{t}_{2n}) \cong H_*(\sym^{n-1}(S^2);\Z) . \] 
%Here we are using a notation introduced in Part I: when $\Sigma$ is a Riemann surface, $\mathbb{S}(\Sigma,n)$ is the graded abelian group $\bigoplus_{i=0}^{n-1}{ \Lambda^i H^1(\Sigma) \otimes_{\Z}\Z[U]/(U^{n+1-i})} = H_*(\sym^n(\Sigma);\Z)$ (here $ \Lambda^i H^1(\Sigma)\otimes \langle U^k \rangle$ is homogeneous of degree $i+2k$).
 
The map $\EuScript{L}_A$ is identified with the map on homology of symmetric products induced by the fundamental class of the correspondence $\widehat{V}$. Monotonicity ensures that there are no `quantum corrections': two holomorphic sections with the same asymptotic limits and the same index must also have the same energy, and the horizontal sections are exactly the ones with zero energy. Thus $\EuScript{L}_A$ sends $\Phi(l)$ to $\Phi(1)$. 

Next, $U^{nm+n-1} \cdot\Phi(1) = U^{n-1}\cdot\Phi(1) = \Phi (\eta^{n-1})$. Finally, $ \EuScript{L}_{D_+}$ maps $\Phi (\eta^{n-1})$ to $1$.  

We deduce that
\[ \EuScript{L}_{X,\pi} (\mathfrak{s}_{m+1,n}) = 
\begin{cases} \pm  U^{n-1}\otimes \lambda,& m\geq 0, \, n\geq 0; \\
0,  & m< 0, \, n\geq 0.
\end{cases}    \]
This is again in agreement with the Seiberg--Witten invariants, calculated in the Taubes chamber of a compatible near-symplectic form. 

The same method applies also to the broken fibration on $(S^1\times S^3)\# \PS^2 \# \overline{\PS^2}$ mentioned above.

\subsubsection{Connected sum with $\PS^2$} 
This is one of the topological examples from \cite{ADK}. Start with a broken fibration $\pi\colon X\to S$ with an isolated critical point $c$. Let $D\subset S$ be a small disc containing $\pi(c)$. There is a torus $Q\subset \pi^{-1}(\partial D)$, the union of the vanishing cycles for radial paths in $D$. Then there is a broken fibration $\pi'\colon X\to S'$ which coincides with $\pi$ outside $D$, and which contains a circle $Z$ of critical points mapping to a small circle in $D$. The manifold $X'$ is diffeomorphic to $X \, \#  \, \PS^2$. 

Let $\beta\in H_2(X' ,Z;\Z)$ be a class with $\beta\cdot [\Sigma]=1$, where $\Sigma$ is a fibre over a point in $\partial D$ (so $\beta \cdot [\bar{\Sigma}]=0$ for fibres $\bar{\Sigma}$ over the centre of $D$). 

\begin{Prop} We have $\EuScript{L}_{X',\pi'}(\mathfrak{s}_\beta)=0.$
\end{Prop}
\begin{pf}
Choose a near-symplectic form $\omega$ on $X'$ such that $Q$ is Lagrangian. The monodromy of $\pi'$ around $\partial D$ is a Dehn twist $\tau_L\in \aut(\Sigma)$ about the circle $L=Q\cap \Sigma$. We analyse the element $\EuScript{L}\in HF_*(Y,\mathfrak{t})$, where $Y$ is the mapping torus of the Dehn twist $\tau_L$ on $\Sigma$ and $\mathfrak{t}$ the restriction of $\mathfrak{s}_\beta$. 

Notice first that $\EuScript{L}\in HF_{\mathrm{odd}}(Y,\mathfrak{t})$. It is known by Seidel's work \cite{Se1} that $HF_{\mathrm{odd}}(Y,\mathfrak{t})\cong H_1(\Sigma,L ;\Z)$, by an isomorphism which is equivariant under the action by symplectic automorphisms of $\Sigma$ on the two sides. If $\delta \in \aut(\Sigma,\omega|\Sigma)$ acts as the identity near $L$, then $\delta_*(\EuScript{L})= \EuScript{L}$, because $\delta$ induces an automorphism of the LHF $Y\to S^1$ which acts trivially on $Q$. But we can obtain any element of $\Sp(H_1(\Sigma);\Z)$ which fixes $[L]$ from such a $\delta$, so there are no non-zero elements of $H_1(\Sigma,L;\Lambda_{\Z})$ fixed by all the automorphisms $\delta$. The result follows.
\end{pf}

\subsubsection{A broken fibration on $(S^1\times S^3)\,\#\, (S^2\times S^2) \,\# \, \overline{\mathbb{P}^2}$.}
In this example we take a step in the direction of an `isotropic blow-up' formula for Lagrangian matching invariants. This involves a model for the differentiable blow-up of a broken fibration described in \cite{ADK}. This model is unusual in that  there is a compatible near-symplectic form for which the exceptional sphere is \emph{isotropic}. The sphere lies over an arc in the base connecting an isolated critical value and a point on a circle of critical values.\footnote{In four-dimensional symplectic geometry, Lagrangian two-spheres have self-intersection $-2$, but in near-symplectic geometry there are other possibilities, including isotropic $(-1)$-spheres.}

Rather than studying the isotropic blow-up process in general, we apply it to our earlier example on $ (S^1\times S^3)\,\# \,(S^2\times S^2)$ and to a variant of that example on $(S^1\times S^3)\, \# \,\PS^2\, \# \,\overline{\PS^2}.$ Let $\widehat{\pi}\colon \widehat{X}\to S^2$ be either of the following two broken fibrations:
\begin{itemize}
\item
The set of critical values is the union of the equator and the point $0\in D_-$.
\item
The regular fibres over $D_-$ have genus 1. In the fibre $T^2=\widehat{\pi}^{-1}(1)$, the vanishing cycle for the path $[0,1]$ is a non-separating circle $L$. The vanishing surface $Q \subset \widehat{\pi}^{-1}(\partial D_-)$ has $Q\cap T^2 = L$.
\item
The fibres over $D_+$ have genus $0$, so $\widehat{X}_+ :=\widehat{\pi}^{-1}(D_+)$ is a trivial bundle $S^2\times D_+\to D_+$. 
\item
There are two possible ways of gluing $\partial \widehat{X}_+ $ to $\widehat{X}_0$. One of these has total space $\widehat{X}\cong (S^1\times S^3) \, \# \,(S^2\times S^2) \,\# \,\overline{\PS^2}$, the other $\widehat{X}\cong (S^1\times S^3)\,\# \,\PS^2 \, \# \, \overline{\PS^2} \, \# \, \overline{\PS^2}$).
\end{itemize} 
We will only consider those homology classes $\beta \in H_2(\widehat{X},Z;\Z)$ with $\beta\cdot [T^2] = 1$, so the fibration over $D_+$ is irrelevant. A useful local model is the fibration
\[q\colon (z,w) \mapsto z^2 + w^2  \]
on 
$ E^r  = \{ (z,w)\in \C^2: |q(z,q) | \leq r , (|z|^2+|w|^2)^2 - |q(z,w)|^2 \leq \epsilon \} $,
with its standard symplectic form $\omega_{\C^2}$ and Lagrangian boundary condition 
\[Q^r = \bigcup_{s\in \partial \bar{D}(0;r)} {Q^r_s},\quad Q^r_s= s^{1/2} S^1 \] 
(where $S^1$ is the unit circle in $\R^2\subset \C^2$). This model was carefully analysed by Seidel \cite[Section 2.3]{Se5}. For the standard almost complex structure, the holomorphic sections with boundary on $Q^r$ are given by
\[ u^r_{a,\pm} ( s) = ( a s+ \bar{a} , \pm \ii (a s- \bar{a}  )),\quad  |a|^2 = r/2.   \]
Moreover, all the sections $u^r_{a,\pm} $ are regular.

We can embed $E^1$ as a sub-fibration of $\widehat{X}_-\to D_-$ (its complement is then a trivial fibration). We are then free to shrink the base from $D_-= \bar{D}(0;1)$ to $\bar{D}(0;r)$, using the Lagrangian boundary condition $Q^r$. By the degeneration argument of \cite[Lemmas 2.14 and 2.15]{Se5} when $r$ is sufficiently small, every pseudo-holomorphic section of $\widehat{X}_- |\bar{D}(0;r) \to \bar{D}(0;r)$ lies inside $E^r$, and is therefore one of the $u^r_{a,\pm}$.

It is important to note that the circles in $Q^r$ traced out by $u^r_{a,+}$ and $u^r_{a,-}$ are \emph{not} mutually homologous, since they intersect transversely at a single point. Hence when we embed $E(r)$ into $X_-$,  $u^r_{a,+}$ and $u^r_{a,-}$ represent \emph{distinct} homology classes $\beta_+$, $\beta_-$. Thus we have moduli spaces of sections $\EuScript{M}^r_{\beta_+}=\{ u^r_{a,+} : |a|^2 = r/2\} $ and $\EuScript{M}^r_{\beta_-}=\{ u^r_{a,-}:|a|^2=r/2 \} $. Each is diffeomorphic to $S^1$ (by $u_{a,\pm}^r\mapsto \sqrt{2} r^{-1/2} a $) but precisely one of these diffeomorphisms is orientation-preserving (another consequence of Seidel's degeneration argument, which shows that $\EuScript{M}^r_{\beta_-}\cup \EuScript{M}^r_{\beta_+}$ is the boundary of a moduli space $\EuScript{M}'\cong S^1\times [-1,1]$). Let $\lambda \subset T^2$ be a loop representing the generator of $\pi_1(\widehat{X})$, and intersecting $L$ transversely at a point. Then the evaluation maps $\ev_+\colon \EuScript{M}^r_{\beta_+}\to L $ and $\ev_-\colon \EuScript{M}^r_{\beta_-}\to L$ are both transverse to $\lambda$. The fibre product moduli spaces, $\EuScript{M}^r_{\beta_\pm}\times_{ ev_{\pm}} \lambda$, are singletons, each consisting of one point, counted with the same signs. Hence
\[ \EuScript{L} _{\widehat{X},\widehat{\pi}}(\mathfrak{s}_{\beta_\pm}) =  \varepsilon [\lambda],  \] 
where $\varepsilon=\pm 1$ is a common sign. Moreover, $\EuScript{L} _{\widehat{X},\widehat{\pi}}(\mathfrak{s}_{\beta'})=0$ for every other class $\beta'$ with $\beta' \cdot T^2 = 1$.

To compare with Seiberg--Witten theory, let us take $\widehat{X}\cong (S^1\times S^3)\# (S^2\times S^2) \# \overline{\PS^2}$. We use the basis $(S,S',E)$ for $H_2(\widehat{X};\Z)$, where $S = [S^2\times \{pt.\}]$, $S' = [\{pt.\}\times S^2]$, $E$ the class of a $(-1)$-sphere in $\overline{\PS^2}$. For the $\spinc$-structure $\mathfrak{s}_{p,q,k}$ with $c_1(\mathfrak{s}_{p,q,k}) = 2 p S+ 2qS'+ (2k+1)E$, one has $d(\mathfrak{s}_{p,q,k}) = 2pq - k(k+1) -1$. For the classes $\beta_\pm$ one has $p=q =1$ and $d(\mathfrak{s}_{\beta_\pm})=1$, and it follows that $\{ \mathfrak{s}_{\beta_+} ,\mathfrak{s}_{\beta_-} \} = \{ \mathfrak{s}_{1,1,0}, \mathfrak{s}_{1,1,1}\}$. By the blow-up formula for Seiberg--Witten invariants, one has
\[ \mathrm{SW}^\tau_{\widehat{X}}(\mathfrak{s}_{1,1, k} ) = \begin{cases}
\pm [\lambda ] ,& k \in \{ 0,1\} , \\ 
0, & k \notin \{0,1 \}. \end{cases}\] 

\subsection{Separating model for $X_1\,\#\,X_2$}
There is a simple realisation of the connected sum (\emph{not} fibre sum!) of broken fibrations.

\begin{centering}
\begin{figure}[t!]
\includegraphics[width=13cm]{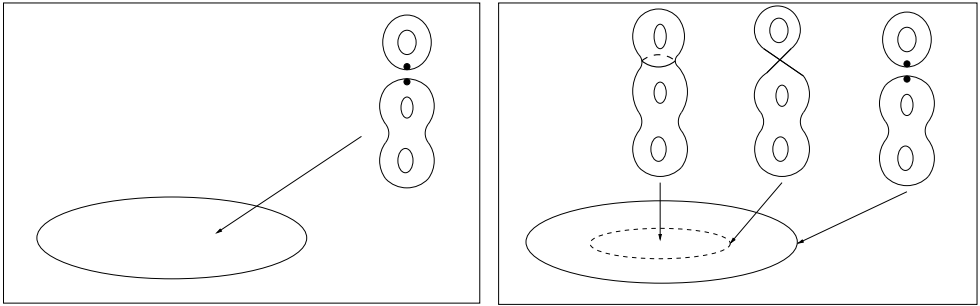}
\caption{\emph{Connected sum of broken fibrations. The first diagram shows the disjoint union of two surface-bundles over the disc. One converts them the disjoint union into a broken fibration as in the second diagram.}}\label{connsum}
\end{figure}\end{centering}

Suppose $\pi_i \colon X_i\to S$, $i=1,2$, are broken fibrations over the same base. Take a simultaneous regular value $s\in S$, and a small coordinate disc $(\Delta,0)\hookrightarrow (S,s)$. There is then a broken fibration $\pi\colon X\to S$ such that 
\begin{itemize}
\item
$\pi^{-1}(S\setminus \Delta)= \pi_1^{-1}(S\setminus \Delta) \sqcup \pi^{-1}_2(S\setminus \Delta)$, and on this region $\pi$ is the disjoint union of $\pi_1$ and $\pi_2$; 
\item
$X^{\crit} \cap \pi^{-1}(\Delta)$ is a circle;
\item $s$ is a regular value and the fibre $X_s=\pi^{-1}(s)$ is the connected sum $\pi_1^{-1}(s)\, \# \,\pi_2^{-1}(s) $.
\end{itemize}
The construction is best understood pictorially (see Figure \ref{connsum}). 

\begin{Lem}
$X\cong X_1\,\#\, X_2$.
\end{Lem}
\begin{pf}
We need to find a separating three-sphere. We work over the unit disc $\Delta$, taking the circle of critical values to be $\{|z|=1/2 \}$. As the connected sum of two surfaces, the central fibre $\pi^{-1}(0)$ contains a separating circle $L$ (drawn in the figure). Over each ray $\{re^{\ii\theta}: 0\leq r\leq 1/2 \}\subset \Delta$, there is a disc in $X$ bounded by $L$ and containing exactly one critical point of $\pi$. Let $Y$ be the union (over $\theta\in[0,2 \pi]$) of these discs. It is the union of a copy of $D^2\times S^1$ (lying over $\{|z| \leq 1/4\}$) and a copy of $S^1\times D^2$ (over $1/4 \leq |z| \leq 1/2$), joined in the obvious way along their common boundary $S^1\times S^1$. Hence $Y$ is diffeomorphic to $S^3$. The complement $X\setminus Y$ breaks into two components $X_i' \cong X_i \setminus B^4$, so $X\cong X_1 \, \# \, X_2$.
\end{pf}

\subsubsection{Vanishing of the invariant} 
The Lagrangian matching invariant of the connected sum of broken fibration vanishes:
\begin{TheoremF}
$\EuScript{L}_{(X_1\# X_2, \pi_1\# \pi_2)}(\mathfrak{s})=0$ for each admissible $\spinc$-structure $\mathfrak{s}$.
\end{TheoremF}
The admissible ones are essentially those corresponding to $n_1$ points on $\Sigma_1$ and $n_2$ points on $\Sigma_2$, where $n_i \leq \frac{1}{2} g(\Sigma_i) -1 $ for $i=1,2$ (though one could do a little better than this). The precise definition was given in Section \ref{definition}.

\begin{pf}
$\EuScript{L}_{(X_1\# X_2, \pi_1\# \pi_2)}(\mathfrak{s})$ factors through the map on homology induced by the fundamental class $[\widehat{V}_L]$: 
\[[\widehat{V}_L]\colon  H_*(\sym^n(\Sigma_1\cup\Sigma_2);\Z)\to H_*(\sym^{n-1}(\Sigma_1 \#\Sigma_2);\Z). \]
But here $L\subset\Sigma_1 \#\Sigma_2 $ is the loop separating the two connect summands, so the fundamental class map is zero (Part I, Lemma 3.18).
\end{pf}

\subsection{Dimensional reduction}

The cobordism category of broken fibrations contains a much simpler category $\mathcal{C}_{\text{Morse}}$ in which all the structures have $S^1$-symmetry. 

An object of $\mathcal{C}_{\text{Morse}}$ is a triple $(\Sigma,c,r)$ consisting of a closed, oriented surface with a locally constant function $c \colon \Sigma \to S^1$ and $r\in H_0(\Sigma;\Z)$. A morphism $(Y,f,a)$ from $(\Sigma_1,c_1,r_1)$ to $(\Sigma_2,c_2,r_2)$ is a cobordism $Y$, equipped with a circle-valued Morse function $f\colon Y\to S^1$ with no interior maximum or minimum, and a homology class $a\in H_1(Y, \partial Y \cup \crit(f) ;\Z)$. There are two requirements on $a$: its boundary must be
\[\partial a = r _2  - r_1 +  \sum_{x\in \crit(f)}{(-1)^{\ind(x)+1} [x]}\in  H_0( \partial Y \cup\crit(f);\Z), \]
and the number $a\cap [F] -\dim(H_1(F))/2 $ must be constant as $F$ ranges over regular fibres $F$ of $f$.

One obtains a broken fibration from a morphism in $\mathcal{C}_{\text{Morse}}$ by crossing with the identity map $\id\colon S^1\to S^1$.
 
\subsubsection{Elementary cobordisms} A basic case to consider is one where $Y$ is an elementary cobordism $C_{\mathrm{elem}}^-$ from $\Sigma$ to $\bar{\Sigma}$; the Morse function $f$ is real valued, and has a single critical point, of index $+1$. 

This gives rise to an invariant
\[ \EuScript{L}(\mathfrak{s}_d)  \colon HF_*(\Sigma\times S^1,\mathfrak{t}_d )\to HF_*(\bar{\Sigma}\times S^1,\mathfrak{t}_d ),  \]
where $\mathfrak{s}_d$ is the unique $\spinc$-structure which is pulled back from $C_{\mathrm{elem}}^-$ and which satisfies $\langle \mathfrak{s}_d, [\Sigma]\rangle = d := \chi(\Sigma)+2n$. 
Under the PSS isomorphisms of Floer and singular homology, this is the map
\[ [\widehat{V}_L ] \colon H_*(\sym^n(\Sigma);\Z) \to H_{*-1}(\sym^{n-1}(\bar{\Sigma});\Z). \]

We can also consider an elementary cobordism $C_{\mathrm{elem}}^+$ running in the opposite direction (so the Morse function has one critical point, of index 2). Its invariant is the adjoint homomorphism (with respect to Poincar\'e duality)
\[   H_*(\sym^{n-1}(\bar{\Sigma});\Z)  \to H_{*+1}(\sym^n(\Sigma);\Z) \]

It is not difficult to give completely explicit formulae for these maps, but what is more illuminating is to observe that they are the same as the maps that arise in a well-known (2+1)-dimensional TQFT, described in \cite{Seg,Don}, which we shall call the Segal-Donaldson model. This can be characterised as follows.

\begin{Defn}
The  {\bf admissible (2+1)-dimensional cobordism category} of degree $d$ is the category in which an object is an oriented 2-manifold, and a morphism is an `admissible cobordism' $Y$ equipped with a line bundle which has degree $d$ over the incoming and outgoing ends. Admissible means that either
\begin{align*}
&\partial Y=\emptyset,\quad H_1(Y;\Z)\cong\Z,\text{ or}\\  
&\partial Y\neq\emptyset,\quad H_1(\partial Y;\Z)\to H_1(Y;\Z)\text{ is onto.}
\end{align*}
\end{Defn}

The Segal-Donaldson model is a TQFT on the category of admissible $(2+1)$-dimensional cobordisms, i.e. a functor from this category to super $\Q$-vector spaces, preserving the multiplicative structures (disjoint union; super tensor product) and dualities. For our purposes,  the morphisms are defined only up to a sign $\pm 1$. 

\begin{Thm}[Donaldson \cite{Don}]
For any sequence of non-negative integers $n_0,n_1,\dots$, there is a unique TQFT $V_d$ on the degree-$d$ admissible cobordism category, defined over $\Q$, such that 
\[ V_d(\Sigma)=\bigoplus_{i\geq 0} {n_i\,\Lambda_i \Sigma}\]
for connected surfaces $\Sigma$, as representations of the mapping class group $\pi_0\diff^+(\Sigma)$. Here $\Lambda_i \Sigma$ denotes $\Lambda^{g(\Sigma)-i}H_1(\Sigma; \Q)$. If $Y^3$ is closed, so that there is a unique line bundle $L$ of degree $d$, then
\[      V_d(Y,L) =  \sum_i{n_i a_i} ,   \]
where $A(Y)=\pm (a_0 + \sum_i{a_i(t^i+t^{-i})})$ is the normalised Alexander polynomial of $Y$. 
\end{Thm}
(`Unique' means up to natural isomorphism.) We shall take $V_d$ to be the unique TQFT with
\[  V_d(\Sigma) = \bigoplus_{j\geq 1} {j \Lambda_{d+j} \Sigma}. \]
This is of interest because then
\[V_d (\Sigma) \cong  H_*(\sym^{n}(\Sigma);\Q) \] 
as representations of $\pi_0\diff^+(\Sigma)$, where $d=\chi(\Sigma)+2n$.

We do not know that Lagrangian matching invariants give a TQFT; they depend, in principle, on the fibrations. However, if we show that $\EuScript{L}(\mathfrak{s}_d)$ coincides with $V_d(C_{\mathrm{elem}}^-)$ then our field theory must coincide with $V_d$.
  
\begin{Lem}
Considered as a map $H_*(\sym^{n}(\Sigma);\Z)\to H_*(\sym^{n-1}(\bar{\Sigma}))$, $\EuScript{L}=\EuScript{L}(\mathfrak{s}_d)$ satisfies the following properties:
\begin{enumerate}
\item
Let $(\alpha,\bar{\alpha}) \in H_1(\Sigma;\Z)\times H_1(\bar{\Sigma};\Z)$ be such that the images of $\alpha$ and $\bar{\alpha}$ in $H_1(C^-_{\mathrm{elem}};\Z)$ are equal.
Then $ \bar{\alpha} \wedge\EuScript{L}(x) = \EuScript{L} (\alpha \wedge x)$. 
\item 
$U  \cdot  \EuScript{L}(x)  =\EuScript{L} ( U \cdot x)$ and $\theta_{\bar{\Sigma}}\cdot \EuScript{L}(x) = \EuScript{L}(\theta_\Sigma \cdot x)$;
\item
on the degree 1 part $H_1(\Sigma;\Z)=H^1(\sym^n(\Sigma);\Z)$, $\EuScript{L}$ is given, up to sign, by $x\mapsto  x\cap [L]$.
\end{enumerate}
Moreover, $\EuScript{L}$ is uniquely characterised (up to sign) by these properties. 
\end{Lem}
\begin{pf}
The first two properties are instances of (\ref{moving between components}). The third is another way of writing the equation $\mathrm{pr}_{1*}[\widehat{V}_L]=[\delta_L]$ for the fundamental class of the vanishing cycle $\widehat{V}_L$ (see Part I, Lemma 3.18).
 
For the uniqueness, we note that by (2) it suffices to show that $\EuScript{L}$ is uniquely determined on the images of the maps $\Lambda^k H_1(\Sigma;\Z)\to H^k\sym^n(\Sigma);\Z)$ sending $x_1\wedge \dots \wedge x_k$ to $\mu(x_1)\cup \dots \cup \mu(x_k)$; this is a simple consequence of (1) and (3). \end{pf}

One can easily verify from Donaldson's model that the Donaldson-Segal map $V_d(C^-_{\mathrm{elem}})$ satisfies these same properties. 
\begin{Cor}
$V_d(C^-_{\mathrm{elem}})= \pm \EuScript{L}$, and $V_d(C^+_{\mathrm{elem}})$ is its adjoint $\pm\EuScript{L}^\vee$ (the Lagrangian matching invariant of $C^+_{\mathrm{elem}}\times S^1$).
\end{Cor}

This has the following consequence.
\begin{TheoremE}
Let $Y$ be a closed 3-manifold with $H_1(Y;\Z)\cong \Z$, and $f\colon Y\to S^1$ a Morse function with connected fibres and no extrema. Let 
\[\id\times f \colon S^1\times Y\to S^1 \times S^1\] 
be the corresponding broken fibration, and let $\mathfrak{s}_d$ be an admissible $\spinc$-structure on $S^1\times Y$ which is pulled back from $Y$. Then 
\[\EuScript{L}_{S^1\times Y, \id\times f}(\mathfrak{s}_d) = \pm \sum_{i\geq 0}{ i a_{d+i}},\] where $A(Y)= \pm (a_0 +\sum_{i\geq 0} {a_i(t^i+ t^{-i})}$ is the Alexander polynomial.
\end{TheoremE}

The assumption $H_1(Y;\Z)\cong \Z$ implies that such Morse functions exist. It may be possible to push the argument through for a general $Y$ with $b_1>0$, replacing the Alexander polynomial by the Milnor--Turaev torsion.

\end{document}